\documentclass[reqno]{amsart}
\usepackage[utf8]{inputenc}
\usepackage[T1]{fontenc}
\usepackage{a4wide}
\usepackage{amsmath,amssymb,amsfonts,amsthm,mathtools,mathrsfs,mathtools,bbm}
\usepackage{graphicx,float,tikz}
\usepackage{hyperref}
\usepackage{aliascnt}
\usepackage[capitalise]{cleveref}
\hypersetup{hypertexnames=false,colorlinks=true,linkcolor=red,citecolor=red,urlcolor=red}

\usepackage[style=numeric,backend=biber,doi=false,isbn=false,url=false,eprint=false,
 giveninits=true,sorting=nyt]{biblatex}
\DeclareFieldFormat{title}{\mkbibemph{#1}}
\DeclareFieldFormat[article,inbook,incollection,inproceedings,thesis]{title}{\mkbibemph{#1}}
\AtEveryBibitem{\clearlist{language}\clearfield{note}\clearfield{month}%
 \clearfield{issn}\clearfield{isbn}\clearfield{doi}\clearfield{url}%
 \clearfield{urldate}\clearfield{eprint}\clearfield{eprinttype}}

\newcommand{\norm}[1]{\left\Vert#1\right\Vert}
\numberwithin{equation}{section}
\allowdisplaybreaks
\newtheorem{theorem}{Theorem}[section]

\newaliascnt{proposition}{theorem}
\newtheorem{proposition}[proposition]{Proposition}
\aliascntresetthe{proposition}
\crefname{proposition}{Proposition}{Propositions}

\newaliascnt{lemma}{theorem}
\newtheorem{lemma}[lemma]{Lemma}
\aliascntresetthe{lemma}
\crefname{lemma}{Lemma}{Lemmas}

\newaliascnt{corollary}{theorem}
\newtheorem{corollary}[corollary]{Corollary}
\aliascntresetthe{corollary}
\crefname{corollary}{Corollary}{Corollaries}

\newaliascnt{assumption}{theorem}
\newtheorem{assumption}[assumption]{Assumption}
\aliascntresetthe{assumption}
\crefname{assumption}{Assumption}{Assumptions}

\theoremstyle{definition}
\newaliascnt{definition}{theorem}
\newtheorem{definition}[definition]{Definition}
\aliascntresetthe{definition}
\crefname{definition}{Definition}{Definitions}

\theoremstyle{remark}
\newaliascnt{remark}{theorem}
\newtheorem{remark}[remark]{Remark}
\aliascntresetthe{remark}
\crefname{remark}{Remark}{Remarks}

\title[Scattering of water waves]{Scattering of water waves in presence\\ of a partially immersed  obstacle}
\thanks{Debayan Maity acknowledges the support of the Department of Atomic Energy, Government of India, under Project Identification No. RTI 40XX and ANRF MATRICS Research Grant ANRF/ARGM/2025/002006. Vicente Ocqueteau acknowledges funding from the European Union (Horizon Europe MSCA project ModConFlex, grant number 101073558).}
\author{Debayan Maity}
\address{TIFR Centre for Applicable Mathematics, Bangalore, Karnataka 560065, India}
\email{debayan@tifrbng.res.in}
\author{Vicente Ocqueteau}
\address{Univ. Bordeaux, CNRS, Bordeaux INP, IMB, UMR 5251, F-33400 Talence, France}
\email{vicente.ocqueteau@math.u-bordeaux.fr}

\author{Marius Tucsnak}
\address{Univ. Bordeaux, CNRS, Bordeaux INP, IMB, UMR 5251, F-33400 Talence, France and Institut Universitaire de France (IUF)}
\email{marius.tucsnak@math.u-bordeaux.fr}


\date{\today}

\begin{document}
\begin{abstract}
We study linear water waves in two dimensions in a fluid of infinite depth containing a fixed, partially immersed obstacle. Under natural geometric assumptions on the wetted boundary, we construct the harmonic lifting and prove that the associated Dirichlet-to-Neumann operator is an unbounded, self-adjoint and non-negative operator on a $L^2$ space. This yields global well-posedness and conservation of energy. We establish the absence of trapped modes and, by a compact perturbation argument, a limiting absorption principle at every positive spectral parameter. Consequently, for all finite energy initial data, the physical local energy tends to zero as $|t|\to\infty$, including in neighborhoods of the obstacle and the contact points. Quantitative resolvent estimates and fractional regularity in the spectral parameter yield algebraic local energy decay for Sobolev initial data in weighted $L^2$ spaces.
\end{abstract}

\maketitle

\section{Introduction}\label{sec:intro}

We study small-amplitude irrotational waves in an inviscid incompressible
fluid that occupies an unbounded domain $\Omega\subset\mathbb R^2$ of infinite
depth. The equilibrium free surface lies on $y=0$, with $y$ increasing downward. A single fixed obstacle intersects that surface. Its wetted boundary is denoted by $\Gamma$, and the exposed free surface is denoted by
$\mathcal E$. The normal $n$ points out of the fluid, so that  $n=-e_y$ on
$\mathcal E$. For a detailed description of the modeling and of the linearization process, we refer, for instance, to Lannes \cite{lannes2013water} and Stoker \cite{stoker1957water}.

\begin{figure}[H]
    \centering
    \begin{tikzpicture}[scale=0.8, >=stealth]
        
        \fill[cyan!15] (-4, -3) -- (-4, 0) -- (-1.5, 0)
            .. controls (-1.2, -2.8) and (-0.2, -1.2) .. (0.5, -1.6)
            .. controls (1.0, -1.9) and (1.3, -0.3) .. (1.5, 0)
            -- (4, 0) -- (4, -3) -- cycle;
        
        \draw[black, thick, dashed] (-4, 0) -- (-1.5, 0);
        \draw[black, thick, dashed] (1.5, 0) -- (4, 0);
        \node[above] at (-2.75, 0) {$\mathcal{E}$};
        \node[above] at (2.75, 0) {$\mathcal{E}$};
        
        \draw[black, thick] (-1.5, 0) 
            .. controls (-1.2, -2.8) and (-0.2, -1.2) .. (0.5, -1.6)
            .. controls (1.0, -1.9) and (1.3, -0.3) .. (1.5, 0);
        \node[below] at (-0.5, -1.15) {$\Gamma$};

        \node at (-2.8, -1.5) {\Large $\Omega$};

        \filldraw (-1.5, 0) circle (1.5pt) node[above=2pt] {$x_-$};
        \filldraw (1.5, 0) circle (1.5pt) node[above=2pt] {$x_+$};
        
        \node[left] at (-4, 0) {$-\infty$};
        \node[right] at (4, 0) {$\infty$};

        \draw[->, thick, gray] (0, 0) -- (0.7, 0) node[right, text=black] {\small $x$};
        \draw[->, thick, gray] (0, 0) -- (0, -0.7) node[right, text=black] {\small $y$};

    \end{tikzpicture}
    \vspace{2mm}
    \caption{Fluid domain $\Omega$, where $\mathcal{E}$ denotes the free surface, $\Gamma$ the wetted boundary of the obstacle, and $x_\pm$ the contact points.}
    \label{fig:physical_domain}
\end{figure}
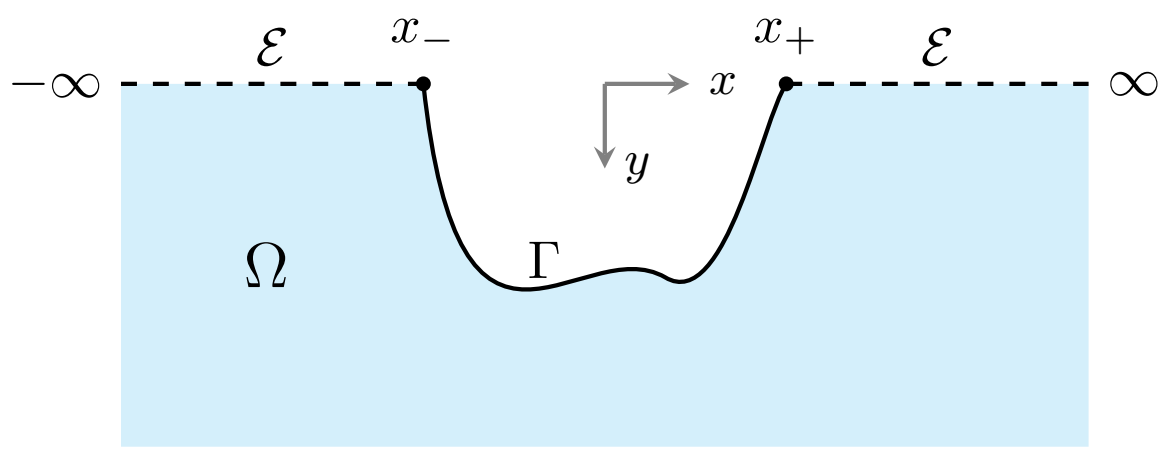

Our main results concern well-posedness and the long-time behavior of
these waves. We identify the operator and form domains of the
Dirichlet-to-Neumann operator, establish a limiting absorption principle,
and prove that the physical local energy tends to zero for every initial
datum of finite energy in the stated inhomogeneous spaces. For spatially
localized, sufficiently regular initial data, we also obtain an explicit,
nonoptimal algebraic decay rate. The precise statements are given in
Section~\ref{sec:main_results}.

Previous work on water waves in the presence of partially immersed
objects includes the recent study by Lannes and Ming
\cite{lannes_ming_2026well}. They establish well-posedness in finite-depth corner domains, allowing for several objects and a possibly emerging bottom.
Their approach defines the Dirichlet-to-Neumann operator variationally
on a homogeneous fractional Sobolev trace space and then studies its
regularity. The finite-depth geometry permits Poincar\'e-type estimates
that are not directly available in our infinite-depth setting.

Our construction follows a different route. We define the
Dirichlet-to-Neumann operator through the almost-everywhere normal trace
of the harmonic lifting and realize it as an unbounded operator on the
space of square-integrable functions on the free surface. We identify
its operator and form domains in standard inhomogeneous Sobolev spaces
and prove self-adjointness in this realization. The required boundary
regularity is obtained by localization and application of Brown's
bounded-domain mixed-problem theorem \cite{Brown1994} under its stated hypotheses.
This construction provides the basis for our spectral analysis and
local energy decay results.

For a fixed obstacle in infinite depth, Weber's thesis
\cite[Theorem~1]{Weber_1976} provides an earlier quantitative
decay result based on multiplier estimates. For symmetric motions and
sufficiently regular, localized data satisfying zero-mean conditions,
he obtains an $O(t^{-1})$ bound for the potential energy in bounded
parts of the free surface and for the kinetic energy in bounded fluid
regions separated from both the obstacle and the free surface.
His geometric hypotheses require symmetry about a vertical line,
star-shapedness of the fluid half-domain with respect to its contact
point, a strictly positive horizontal component of the oriented wetted
tangent, and a straight wetted segment near the contact point.
These conditions do not follow from our geometric assumptions, so his
decay theorem does not directly cover the class considered here.
Our estimates have a slower rate but allow local energy observations
up to the wetted boundary and the contact points, without symmetry,
local straightness, or zero-mean conditions on the data.

A substantial literature also concerns freely floating bodies, whose
motion is coupled to that of the fluid. The foundational papers of
Fritz John \cite{John1949OnTM,John1950OnTM} were followed, among
other contributions, by Ursell \cite{Ursell1949ONTH} and Hazard and Loret
\cite{HazardLoret2007}. More recently, Ocqueteau and Tucsnak
\cite{Ocqueteau2026} studied the initial value problem for a vertically
moving floating cylinder, and Lannes and Paulsen
\cite{lannes2025fjohnmodelcummins} investigated the freely floating
body problem and its relation to Cummins' equations.

The work \cite{HazardLoret2007} is particularly
close to ours: both studies consider infinite depth, and we follow
their perturbation approach to the limiting absorption principle.
Their paper already contains the water-wave LAP and the abstract
compact perturbation theorem underlying this approach. Our detailed
proofs of these results are included to make their application to the
operator realization used here explicit and self-contained.
In particular, their formulation does not identify the domain of the
Dirichlet-to-Neumann operator as an inhomogeneous Sobolev space on the
free surface. Here we construct this operator through its normal trace, identify its exact operator and form domains, and verify the hypotheses of the perturbation argument in the corresponding weighted spaces. The harmonic lifting, the mixed boundary estimates at the contact points, and the compact comparison of resolvents are thus established within a single functional framework.

The main new conclusions concern the physical local energy:
\cref{local_energy_decay} establishes its decay for every datum in the
stated energy space, and \cref{thm:full_polynomial_decay} gives an
explicit algebraic estimate under additional spatial localization and
Sobolev regularity. In particular, these observations include fluid
regions that meet the obstacle and the contact points. The quantitative result is the main contribution of this work; the precise
statements and an account of the proof strategy are gathered in
\cref{sec:main_results}.

The paper is organized as follows. In \cref{subse_not}, we introduce the geometry and notation and establish the basic properties of the energy spaces with boundary $L^2$ traces. We state our main results in \cref{sec:main_results}. In Sections~\ref{sec:spaces} and~\ref{sec_self}, we construct the obstacle operator and establish their self-adjointness and the well-posedness. In \cref{sec:free_waves}, we study the free problem and derive the resolvent estimates needed for the subsequent analysis. In \cref{sec:free}, we establish the limiting absorption principle and quantitative weighted resolvent estimates. In \cref{sec:local_energy_decay}, we prove \cref{local_energy_decay} by combining absolute continuity of the spectral measures with the compactness of the physical local-energy observations after spectral localization. In \cref{sec:quantitative_decay}, we establish \cref{thm:full_polynomial_decay,cor:polynomial_decay_consequences} using a Fourier decay estimate for H\"older continuous amplitudes and a decomposition into low-, intermediate-, and high-frequency components. Finally, in \cref{appendix:conformal_mapping}, we construct a conformal map used in the analysis of the Dirichlet map, while in \cref{appendix:lap}, we collect the abstract spectral results used throughout the paper.

\section{Notation and preliminaries}\label{subse_not}

\subsection{Geometry}

Throughout the paper we impose the following geometric assumptions.

\begin{assumption}\label{ass:geometry}\label{assumption_gamma}
Let $x_-<x_+$ and let $\gamma\in C^0([x_-,x_+])\cap C^2((x_-,x_+))$ satisfying
\[
 \gamma(x_\pm)=0,\qquad
 \gamma'(x_-)>0>\gamma'(x_+),
\]
\[
 \gamma(x)>0 \qquad\qquad(x_-<x<x_+),
\]
with $\gamma'(x_-)$ and $\gamma'(x_+)$ possibly equal to $\infty$ and $-\infty$, respectively.
We also assume that the wetted boundary is the graph $\Gamma$ of the function $\gamma$ defined above, i.e., that
 \begin{equation}\label{def_gamma_e}
\Gamma=\left\{\begin{bmatrix} x\\ \gamma(x)\end{bmatrix}\ \ |\ \ x\in[x_-,x_+]\right\},\qquad
 \mathcal E=\mathbb R\setminus[x_-,x_+].
\end{equation} 
Extending $\gamma$ by zero outside $[x_-,x_+]$ and denoting this extension by
$\overline\gamma$, the domain filled by the fluid is
\begin{equation}\label{def_omega_mare}
 \Omega=\left\{\begin{bmatrix} x\\ y\end{bmatrix}\in \mathbb{R}^2\ \ |\ \ y>\overline \gamma(x)\right\}.
 \end{equation}
\end{assumption}

The strict inequalities on the endpoint derivatives express
non-tangential contact.   

Denote by $\mathbb{H}$ the upper half-plane
\begin{equation}\label{defhmare}
\mathbb H=\left\{\begin{bmatrix} x\\ y\end{bmatrix}\in \mathbb{R}^2\ \ | \ \ y>0\right\}
\end{equation}
The bounded immersed part of the obstacle is
\begin{equation}\label{immersed}
\mathcal B=\mathbb H\setminus\overline\Omega.
\end{equation}
We denote by $B_r$ the open ball of $\mathbb R^2$ of radius $r>0$ centered at the origin. The coordinate $y$ increases downward, and the normal $n$ points out of
the fluid; in particular, $n=-e_y$ on $\mathcal E$.

\subsection{Function spaces and boundary traces}

For a linear operator $A$, its domain is denoted by $\mathcal D(A)$.
For normed spaces $X$ and $Y$, $\mathcal L(X,Y)$ denotes the space of
bounded linear operators from $X$ to $Y$, and $\mathcal L(X)=\mathcal L(X,X)$.
Hilbert-space inner products are taken to be linear in the first argument.

We denote by $\mathbb{N}$ the set of positive integers, i.e., 
\[
\mathbb{N}=\{1,2,\dots\}.
\]
For $d\in\mathbb N$ and an open set $\mathcal O\subset\mathbb R^d$, we denote by $C_c^\infty(\mathcal{O})$ the space of infinitely differentiable functions with compact support in $\mathcal{O}$, and by $C_c^\infty(\overline{\mathcal{O}})$ the space formed by the restrictions to $\overline{\mathcal{O}}$ of functions from $C_c^\infty(\mathbb{R}^d)$.
 For $m\in\mathbb N$, we denote by $H^m(\mathcal{O})$ the Sobolev space of square integrable functions in $\mathcal{O}$ whose distributional derivatives, up to order $m$, belong to $L^2(\mathcal{O})$. The space $H^m_{\rm loc}(\mathcal{O})$ consists of the functions
lying in $H^m(\widetilde{\mathcal{O}})$ for each open bounded subset $\widetilde{\mathcal{O}}$ satisfying  $\overline{\widetilde{\mathcal{O}}} \subset  \mathcal{O}$. 

\begin{remark}\label{rem:homogeneous_boundary_traces}
 If $\mathcal{O}$ has locally Lipschitz boundary, then any function $f \in L^2_{\rm loc}(\mathcal{O})$ with $\nabla f \in L^2(\mathcal{O})^2$ belongs to $H^1$ locally up to the boundary: for every $z\in\partial\mathcal O$ there is a neighbourhood
$U$ of $z$ such that $f|_{\mathcal O\cap U}\in H^1(\mathcal O\cap U)$, see for instance Galdi \cite[Lemma II.6.1]{galdinsbook}. 
Consequently, the standard Sobolev trace theorem yields a Dirichlet trace on each bounded Lipschitz boundary patch. These local traces agree on overlaps and therefore define a consistent local boundary trace. We use these traces throughout the paper, including on boundary patches that contain the contact points.
\end{remark}

For $r\geqslant 0$, $s\in \mathbb{R}$ and interval $I\subset \mathbb{R}$,  
we denote by $H^r(I)$ the space of restrictions to $I$ of functions from
$H^r(\mathbb R)$ and we set
\begin{equation}\label{not_cu_s}
 L_s^2(I)=\{f:\langle x\rangle^s f\in L^2(I)\},\qquad
 \langle x\rangle=(1+x^2)^{\frac{1}{2}}.
\end{equation}
The weighted norm is
\begin{equation}\label{eq:weighted_norm}
\|f\|_{L_s^2(I)}=\|\langle x\rangle^s f\|_{L^2(I)}.
\end{equation}
We note that for every $s>0$, the space $L_{-s}^2(I)$ is the dual of $L_s^2(I)$
with respect to the pivot space $L^2(I)$. 

For a Banach space $Y$, the space $C(\mathbb{R}; Y)$ denotes the space of continuous functions from $\mathbb{R}$ to $Y$. Note that this notation implies no boundedness at infinity, which distinguishes it from the space $C_b(\mathbb{R}; Y)$ of bounded continuous functions. On the disconnected free surface we use the direct sums
\[
H^r(\mathcal E)=H^r((-\infty,x_-))\oplus H^r((x_+,\infty)) \qquad\qquad(r>0),
\]
\[
L_s^2(\mathcal E)=L_s^2((-\infty,x_-))\oplus L_s^2((x_+,\infty)) \qquad\qquad(s\in \mathbb{R}),
\]
with the squared norms equal to the sums of the squared component norms.
In particular, no matching condition is imposed between the two contact points.

For a Banach space $Y$, an interval $J\subset\mathbb R$ and
$0<\alpha<1$, we denote by $C^{0,\alpha}(J;Y)$ the bounded
$Y$-valued H\"older continuous functions, with seminorm and norm
\begin{equation}\label{eq:holder_norms}
 \begin{aligned}
 [f]_{C^{0,\alpha}(J;Y)}
 &=\sup_{\substack{r,q\in J\\r\ne q}}
           \frac{\|f(r)-f(q)\|_Y}{|r-q|^\alpha},\\
 \|f\|_{C^{0,\alpha}(J;Y)}
 &=\sup_{r\in J}\|f(r)\|_Y+[f]_{C^{0,\alpha}(J;Y)}.
 \end{aligned}
\end{equation}
The target space $Y$ is omitted for scalar-valued functions. We write
\begin{equation}\label{eq:compact_holder_space}
 C_c^{0,\alpha}(\mathbb R;Y)
 =\{f\in C^{0,\alpha}(\mathbb R;Y):
                          \operatorname{supp}f\text{ is compact}\}.
\end{equation}

\subsection{Fourier transform}

For the one dimensional Fourier transform we use the unitary form
\begin{equation}\label{eq:fourier_convention}
\mathcal{F}(f)(\xi)=\hat f(\xi)=\frac{1}{\sqrt{2\pi}} \int_\mathbb{R} {\rm e}^{-i x \xi} f(x)\, {\rm d}x
\qquad\qquad(f\in L^1(\mathbb{R})).
\end{equation}
The transform extends by density to a unitary operator on $L^2(\mathbb R)$,
and Plancherel's identity reads
\begin{equation}\label{id_parseval}
    \|f\|_{L^2(\mathbb{R})}=\|\hat f\|_{L^2(\mathbb{R})} \qquad\qquad(f\in L^2(\mathbb{R})).
\end{equation}

\subsection{Energy spaces with boundary $L^2$ traces}
\label{subsec:energy_spaces}

For $G \in \{\Omega, \mathbb{H}\}$, we denote  $\Sigma_\Omega = \mathcal{E}$ and $\Sigma_{\mathbb{H}} = \mathbb{R}$, and we define the space
\begin{equation}\label{eq:robin_energy_space}
 \mathscr W_G=\left\{f\in L^2_{\mathrm{loc}}(G)\ \big|\ \nabla f\in L^2(G)^2,\quad f|_{\Sigma_G}\in L^2(\Sigma_G)\right\},
\end{equation}
equipped with the inner product and norm
\begin{equation}\label{eq:robin_energy_form}
 a_G(f,g)=\int_G\nabla f\cdot\overline{\nabla g} + \int_{\Sigma_G} f \,\overline g, \qquad \|f\|_{\mathscr W_G}^2=a_G(f,f).
\end{equation}
Here, the boundary trace is understood in the sense of Remark \ref{rem:homogeneous_boundary_traces}.

We define the closed subspace $\mathscr Z_G \subset \mathscr W_G$ of functions with zero trace on $\Sigma_G$:
\begin{equation}\label{eq:z_space}
 \mathscr Z_G = \left\{w\in \mathscr W_G \ \big|\ w|_{\Sigma_G}=0\right\},
\end{equation}
endowed with the inner product $\langle f,g \rangle_{\mathscr Z_G} = \int_G \nabla f \cdot \overline{\nabla g}$.

\begin{lemma}\label{lem:completeness_W_Z}
The spaces $(\mathscr W_G, \|\cdot\|_{\mathscr W_G})$ and $(\mathscr Z_G, \|\cdot\|_{\mathscr Z_G})$ are Hilbert spaces.
\end{lemma}
\begin{proof}
We only verify the completeness of $\mathscr W_G$, as the completeness of $\mathscr Z_G$ follows by the same argument combined with the continuity of the trace. Let $(f_n)$ be a Cauchy sequence in $\mathscr W_G$. Then $(\nabla f_n)$ converges to some $F$ in $L^2(G)^2$ and $(f_n|_{\Sigma_G})$ converges to some $h$ in $L^2(\Sigma_G)$. Let $\mathcal O$ be a bounded connected Lipschitz subdomain of $G$ such that $\Sigma_{\mathcal O} := \partial \mathcal O \cap \Sigma_G$ has positive boundary measure. The Friedrichs-Poincaré inequality with a boundary term (see, e.g., Toselli and Widlund \cite[Lemma A.14]{TW05}) yields
\begin{equation} \label{eq:poincare_friedrichs}
  \|{\rm u}\|_{L^2(\mathcal O)} \leqslant C_{\mathcal O} \left( \|\nabla {\rm u}\|_{L^2(\mathcal O)} + \|{\rm u}\|_{L^2(\Sigma_{\mathcal O})} \right) \qquad ({\rm u}\in H^1(\mathcal O)).  
\end{equation}
Applying this inequality to the differences $f_n - f_m$ shows that $(f_n)$ is a Cauchy sequence in $H^1(\mathcal O)$. In particular, we obtain $f_n \to f$ in $H^1_{\mathrm{loc}}(G)$, with $\nabla f = F$. By the continuity of the local trace operator, we also have $f|_{\Sigma_{\mathcal O}} = h|_{\Sigma_{\mathcal O}}$. From this, we infer that $f|_{\Sigma_G} = h$, and thus $f \in \mathscr W_G$. Finally, it is straightforward to check that $f_n \to f$ in the $\mathscr W_G$ norm. Therefore, $\mathscr W_G$ is complete.
\end{proof}

\begin{lemma} \label{lem:dense}
   $C_c^\infty(\overline{\Omega})$ is dense in $\mathscr{W}_{\Omega}.$ 
\end{lemma}

\begin{proof}
It suffices to consider real-valued functions, since the real and
imaginary parts can be approximated separately. For $f\in\mathscr W_\Omega$,
set $f_M=\max(-M,\min(f,M))$. Using the dominated convergence it follows that
\[
\|f_M-f\|_{\mathscr W_\Omega}\longrightarrow0
\qquad\qquad(M\to\infty).
\]
Indeed, $|\nabla(f_M-f)|\leqslant
\mathbf 1_{\{|f|\geqslant M\}}|\nabla f|$, while the difference
of the boundary traces tends to zero pointwise and is bounded in
modulus by the modulus of the trace of $f$. Thus we may assume that $|f|\leqslant M$.

We next use a classical logarithmic cutoff construction;, see, for instance, the proof of Lemma 3.2 in Eichmair and Metzger  \cite{EichmairMetzger2012}. More precisely,
for $R>1$, we define the radial Lipschitz cutoff
\[
\chi_R(z)=
\begin{cases}
1,& |z|\leqslant R,\\
\displaystyle\frac{\log(R^2/|z|)}{\log R},&R<|z|<R^2,\\
0,&|z|\geqslant R^2.
\end{cases}
\]
A direct calculation yields
\[
\|\nabla\chi_R\|_{L^2(\mathbb R^2)}^2=\frac{2\pi}{\log R}.
\]
Consequently,
\[
\|\nabla(\chi_Rf-f)\|_{L^2(\Omega)}
\leqslant\|(1-\chi_R)\nabla f\|_{L^2(\Omega)}
+M\left(\frac{2\pi}{\log R}\right)^{1/2}
\longrightarrow0.
\]
Dominated convergence on $\mathcal E$ also gives
$\|(\chi_R-1)f\|_{L^2(\mathcal E)}\to0$.
Hence $\chi_Rf\to f$ in $\mathscr W_\Omega$.

For fixed $R$, the function $\chi_Rf$ belongs to $H^1(\Omega)$
and has bounded support. Choose a sufficiently large bounded
Lipschitz truncation of $\Omega$ containing this support in its
interior relative to $\overline\Omega$. The Sobolev extension theorem,
followed by multiplication by a smooth cutoff and mollification,
provides a sequence of functions $(\varphi_n)\subset C_c^\infty(\mathbb R^2)$, supported in a common bounded set, such that
\[
\varphi_n|_\Omega\longrightarrow\chi_Rf
\quad\hbox{in }H^1(\Omega).
\]
The trace theorem on a bounded Lipschitz truncation containing
their supports also gives convergence in $L^2(\mathcal E)$.
Thus the convergence holds in $\mathscr W_\Omega$.
Combining the three approximations proves the assertion.
\end{proof}

\section{Statement of the main results}\label{sec:main_results}

\subsection{The evolution problem}

If $\zeta$ denotes downward displacement of the free surface and $\varphi$ the velocity potential, the linearized kinematic and dynamic conditions are $\zeta_t=\partial_y\varphi$ and $\varphi_t=g\zeta$ on $\mathcal{E}$, where $g$ is the acceleration due to gravity. Combining these conditions on $\mathcal{E}$ with the impermeability condition on $\Gamma$ and the incompressibility and irrotationality conditions on $\Omega$, we get the following mixed problem for $\varphi$:
\begin{equation}\label{phi_coupled}
\left\{
 \begin{aligned}
 &\Delta\varphi(t,x,y)=0& \left(\begin{bmatrix}x \\ y \end{bmatrix}\in\Omega\right),\\
 &\frac{\partial\varphi}{\partial n}(t,x,y)=0& \left(\begin{bmatrix} x \\  y\end{bmatrix}\in \Gamma\right),\\
  &\varphi_{tt}(t,x,0)-g\frac{\partial \varphi}{\partial y}(t,x,0)=0& \left(x\in \mathcal E\right) .
 \end{aligned}
 \right.
\end{equation}
Denoting $v=\varphi|_{\mathcal E}$, the general approach in water waves theory is to rephrase the system \eqref{phi_coupled} as an evolution 
equation for $v$. This approach is well established in the absence of partially immersed obstacle or walls, but it has been much less studied
 in their presence. Reducing the system to an evolution equation involving just $v$ requires several ingredients, more difficult to 
 construct in the presence of partially immersed obstacles.  
The first one, denoted $D_\Omega$, is called {\em Dirichlet map} on $\Omega$ and it associates to each ${\rm v}\in H^1(\mathcal{E})$ the function $D_\Omega {\rm v}=\phi$, where
$\phi$ is the unique solution in $\mathscr{W}_\Omega$ of
\begin{equation}\label{Dirichlet_op_pb}
\left\{
\begin{alignedat}{2}
    &\Delta \phi(x,y)=0 && \qquad \qquad \left(\begin{bmatrix}x\\ y\end{bmatrix}\in\Omega\right),\\
    &\frac{\partial \phi}{\partial n}(x,y)=0 && \qquad \qquad \left(\begin{bmatrix}x\\ y\end{bmatrix}\in\Gamma\right),\\
    &\phi(x,0)={\rm v}(x) && \qquad \qquad\left(x\in\mathcal{E}\right).
\end{alignedat}
\right.
\end{equation}
The notion of solution in $\mathscr{W}_\Omega$ of \eqref{Dirichlet_op_pb}  is the variational
one, which means that the Dirichlet trace of $\phi$ on $\mathcal{E}$,
understood as in Remark~\ref{rem:homogeneous_boundary_traces}, equals to ${\rm v}$ and
\begin{equation}\label{Dirichlet_op_pb_var}
\int_\Omega \nabla\phi\cdot \nabla\eta \, {\rm d}x\, {\rm d}y = 0,
\end{equation}
for every $\eta\in \mathscr{Z}_\Omega$.

Moreover, as shown in Proposition \ref{prop:h1lifting} below,  for every ${\rm v}\in H^1(\mathcal{E})$, the function
$D_\Omega {\rm v}$ is continuous and bounded on $\overline\Omega$.

The second ingredient, designed as {\em Dirichlet-to-Neumann map} in the remaining part of the paper, is the
operator associating to each ${\rm v}\in H^1(\mathcal{E})$ the trace on $\mathcal{E}$ of the normal derivative of $D_\Omega {\rm v}$. 
Equivalently, we can write 
\begin{equation}
(\Lambda_\Omega {\rm v})(x)=-\frac{\partial (D_\Omega {\rm v})}{\partial y}(x,0)\qquad\qquad(x\in\mathcal{E}).
\end{equation}
In  Section \ref{sec_self} it is proved  that $\Lambda_\Omega:H^1(\mathcal{E})\to L^2(\mathcal{E})$ is a self-adjoint, positive, one-to-one operator on $L^2(\mathcal{E})$ and that $\mathcal{D}(\Lambda_\Omega^\frac12)=H^\frac12(\mathcal{E})$. The proof of the above properties of $\Lambda_\Omega$ is a point where the presence of the obstacle (through possible singularities of $D_\Omega \rm v$ around the corners) and the unboundedness (both in the $x$ and $y$ directions) of $\Omega$ considerably complicates the analysis.

With the Dirichlet-to-Neumann map $\Lambda_\Omega$ at hand, the system of equations for the unknown $v:\mathbb{R}\times \mathcal{E} \to \mathbb{R}$, describing the small oscillations of the free surface in the presence of the partially immersed obstacle $\mathcal{B}$, can be written as
\begin{equation}\label{v_eq}
\left\{
\begin{aligned}
 &v_{tt}(t,x)+g(\Lambda_\Omega v)(t,x)=0\qquad\qquad(t\in \mathbb{R}, x\in \mathcal{E}),\\ &v(0,x)=v_0(x),\quad v_t(0,x)=v_1(x) \qquad\qquad (x\in \mathcal{E}).
 \end{aligned}
 \right.
\end{equation}
For comparison, in the absence of an obstacle the free surface is
$\mathbb R$ and the Dirichlet-to-Neumann operator is the Fourier
multiplier
\[
 \widehat{\Lambda_{\mathbb H}{\rm v}}(\xi)
   =|\xi|\widehat{{\rm v}}(\xi),
 \qquad\qquad({\rm v}\in H^1(\mathbb R)).
\]
The evolution equation becomes $v_{tt}+g\Lambda_{\mathbb H}v=0$,
with dispersion relation $\omega^2=g|\xi|$. This explicit free
operator serves as the reference for the spectral analysis of
$\Lambda_\Omega$. Its properties and the required resolvent estimates
are developed in Section~\ref{sec:free_waves}.

\subsection{Well-posedness and limiting absorption principle}

The remaining main results tackle the initial value problem \eqref{v_eq}.
We first remark that this problem is well-posed, in the sense of the result below, which will be proved in Section \ref{sec_self}.

\begin{proposition}\label{thm:well_posedness_energy}
With the above notation and under Assumption~\ref{ass:geometry}, for every
$v_0 \in H^{\frac12}(\mathcal E),\  v_1\in L^2(\mathcal E)$ there exists a
unique solution of \eqref{v_eq}, satisfying
\begin{equation}\label{reg_slaba}
 v\in C(\mathbb R;H^{\frac12}(\mathcal E))\cap C^1(\mathbb R;L^2(\mathcal E)),
\end{equation}
 and this solution depends
continuously on the initial data on bounded time intervals.
If $v_0\in H^1(\mathcal E),\ v_1\in H^{\frac{1}{2}}(\mathcal E)$, then
\begin{equation}\label{regtare}
 v\in C(\mathbb R;H^1(\mathcal E))\cap C^1(\mathbb R;H^{\frac{1}{2}}(\mathcal E))
                         \cap C^2(\mathbb R;L^2(\mathcal E)).
\end{equation}
The conserved energy is
\begin{equation}\label{eq:totalenergy}
 \begin{aligned}E(t)&=\frac12\|v_t(t,\cdot)\|_{L^2(\mathcal{E})}^2+
           \frac g2\|(\Lambda_\Omega^{\frac{1}{2}}v)(t,\cdot)\|_{L^2(\mathcal{E})}^2\\
     &=\frac12\|v_t(t,\cdot)\|_{L^2(\mathcal{E})}^2+\frac g2\|\nabla (D_\Omega v)(t,\cdot)\|_{L^2(\Omega)}^2.
\end{aligned}
\end{equation}
\end{proposition}

\begin{remark}\label{rem:energy_norm}
The conserved energy provides quantitative control over $v_t$ and $\Lambda_\Omega^{\frac{1}{2}}v$ in $L^2(\mathcal{E})$. However, it does not prevent the $L^2$ norm of $v(t, \cdot)$ from growing linearly in time, which is why the solution space $C(\mathbb{R}; H^{\frac{1}{2}}(\mathcal{E}))$ does not imply uniform boundedness for all $t \in \mathbb{R}$. Furthermore, as established in \cref{lem:quarter_plane_comparison} below, the coercivity of the semi-norm $\|\Lambda_\Omega^{\frac{1}{2}}v\|_{L^2(\mathcal{E})}$ can be understood by comparing it to the Neumann Laplacian on the free surface components. This non-local effect of the energy semi-norm is a characteristic feature of the problem, analogous to the difference of constants term in the finite-depth framework of \cite{lannes_ming_2026well}.
\end{remark}

The sense in which function $v$ with the regularity \eqref{reg_slaba} satisfies \eqref{v_eq}
is the variational one, i.e.,
\[
\frac{{\rm d}}{{\rm d}t} \langle v_t(t,\cdot),\eta\rangle_{L^2(\mathcal{E})}+g\langle (\Lambda_\Omega^\frac12 v)(t,\cdot),
\Lambda_\Omega^\frac12 \eta \rangle_{L^2(\mathcal{E})} = 0 \qquad\qquad \left(t\in\mathbb{R} \ {\rm a.e.}, \eta\in H^\frac12(\mathcal{E})\right).
\]
Our second main result shows that $\Lambda_\Omega$ satisfies a limiting absorption principle, designed as LAP in the remaining part of this work. More precisely, we have:

\begin{proposition}\label{thm:mainlap}
With the notation and under the assumptions in \cref{thm:well_posedness_energy}, for every
$s\in \left[1,\frac32\right)$ and $\lambda>0$ we have
\[
 \lim_{\varepsilon \to 0+} (\Lambda_\Omega-\lambda\mp i \varepsilon)^{-1} =
 (\Lambda_\Omega-\lambda\mp i0)^{-1}
 \in\mathcal L(L_s^2(\mathcal E),L_{-s}^2(\mathcal E)),
\]
exist in operator norm, locally uniformly in $\lambda$.
\end{proposition}

\subsection{Local energy decay}

The LAP, the absence of eigenvalues of $\Lambda_\Omega$ and compactness
of the spectrally localized energy observations yield the following
decay result for every datum in the stated energy space.

\begin{theorem}\label{local_energy_decay}
Under the hypotheses of \cref{thm:mainlap}, let $v$ be the solution of \eqref{v_eq}
with $v_0\in H^{\frac{1}{2}}(\mathcal E)$ and $v_1\in L^2(\mathcal E)$. Then for $|t|\to\infty$
we have the weak convergence
\[
 (\Lambda_\Omega^{\frac{1}{2}}v)(t,\cdot) \to 0 \quad \text{and} \quad v_t(t,\cdot) \to 0 \qquad\text{weakly in } L^2(\mathcal E).
\]
Moreover, setting $\Omega_R=\Omega\cap B_R$,
$\mathcal E_R=\mathcal E\cap(-R,R)$ for every $R>0$ such that $[x_-,x_+]\subset(-R,R)$, we have
\begin{equation}\label{eq:physical_local_energy_decay}
 E_R(t)=\frac{g}{2}\int_{\Omega_R}|\nabla (D_\Omega v)(t,\cdot)|^2
              +\frac{1}{2}\int_{\mathcal E_R}|v_t(t,\cdot)|^2
                   \longrightarrow0 \text{ as }|t|\to\infty.
\end{equation}
\end{theorem}
Since the downward surface displacement is $\zeta=v_t/g$, this is the physical local
mechanical energy multiplied by the constant $g/\rho$, where $\rho$ is
the fluid density. 
Finally and most importantly, we obtain the following local energy decay result:

\begin{theorem}\label{thm:full_polynomial_decay}
With the notation and under the assumptions in \cref{local_energy_decay}, fix
$1\leqslant s<\frac{3}{2}$ and $0<\beta<s-\frac{1}{2}$. Suppose that
\[
 v_0\in H^1(\mathcal E)\cap L_s^2(\mathcal E),\qquad
 v_1\in H^{\frac12}(\mathcal E)\cap L_s^2(\mathcal E).
\]
The weighted spaces and norms are defined in
\eqref{not_cu_s}--\eqref{eq:weighted_norm}. Let $v$ be the solution of
\eqref{v_eq}, and put
\begin{equation}\label{eq:data_norm_decay}
 \mathcal N_s(v_0,v_1)=
 \|v_0\|_{L_s^2(\mathcal{E})}+\|v_1\|_{L_s^2(\mathcal{E})}
 +\| v_0\|_{H^1(\mathcal E)}+\|v_1\|_{H^\frac12(\mathcal E)}.
\end{equation}
For every fixed $R > 0$ such that $[x_-,x_+]\subset(-R,R),$ the physical local energy $E_R(t)$ satisfies
\begin{equation}\label{eq:optimized_logarithmic_decay}
 E_R(t)\leqslant C_{R,s,\beta}\mathcal N_s(v_0,v_1)^2
        (1+|t|)^{-\gamma}\log(2+|t|)
        \qquad(t\in\mathbb R),
\end{equation}
where 
\begin{equation}
\gamma  =  \frac{2\beta}{29+3\beta}.
\end{equation}
\end{theorem}

The exponent $\gamma$ arises from balancing the decay contribution on an intermediate spectral interval with the low- and high-frequency tails. The denominator reflects the losses incurred in the weighted resolvent estimates as the intermediate spectral interval expands. We do not claim that this exponent represents an intrinsic or optimal decay rate. The precise choice of the cutoff parameters is explained in \cref{sec:quantitative_decay}.

The following consequence applies to spatially localized initial data.

\begin{corollary}\label{cor:polynomial_decay_consequences}
Under the geometric assumptions of \cref{thm:full_polynomial_decay},
let $M>0$ be such that $[x_-,x_+]\subset(-M,M).$ 
Suppose that $v_0\in H^1(\mathcal E)$ and
$v_1\in H^{\frac12}(\mathcal E)$ 
satisfy 
\begin{equation}
    v_0(x) =  v_1(x) = 0 \mbox{ for } |x| > M.
\end{equation}
Then, for every $0<\delta<1/16$, the solution of \eqref{v_eq}
satisfies
\begin{equation}\label{eq:bounded_support_decay}
 \begin{split}
 E_R(t)\leqslant C_{R,M,\delta}(1+|t|)^{-\delta}
 \bigl(\|v_0\|_{H^1(\mathcal E)}^2+
             \|v_1\|_{H^{\frac{1}{2}}(\mathcal E)}^2\bigr).
 \end{split}
\end{equation}
\end{corollary}

For bounded supported initial data, the decay threshold $\frac{1}{16}$ follows from
$\frac{2\beta}{29+3\beta}\longrightarrow\frac1{16}$
as $\beta\to1^-$. The logarithmic factor can be absorbed by taking any decay exponent strictly below $\frac{1}{16}$. Thus, the endpoint exponent $\frac{1}{16}$ is not claimed. The bounded-support assumption ensures that the required weighted norms are controlled by the corresponding unweighted norms.

The decay rate can be improved by imposing additional regularity on the initial data with respect to the operator $\Lambda_\Omega$. More precisely, \cref{prop:higher_regularity_decay} in \cref{subsec:higher_regularity_decay} establishes stronger decay estimates for smooth, compactly supported initial data belonging to the prescribed domains of higher powers of $\Lambda_\Omega$.

\subsection{Contributions and main ideas of the proofs}
\label{subsec:proof_strategy}

The main new contribution of this work is the quantitative decay estimate in \cref{thm:full_polynomial_decay}, together with the qualitative decay of the physical local energy established in \cref{local_energy_decay} for the present geometry and energy space. The well-posedness and limiting absorption principle (LAP) results in \cref{thm:well_posedness_energy,thm:mainlap} build on the theory developed in \cite{HazardLoret2007} and, in the case of well-posedness for a different finite-depth setting, in \cite{lannes_ming_2026well}. We do not claim novelty for the abstract compact perturbation principle. Rather, our contribution lies in adapting these methods to the strong $L^2(\mathcal E)$ realization of $\Lambda_\Omega$, with particular attention to the precise characterization of its operator and form domains and the estimates near the contact points.

The proof of \cref{local_energy_decay} relies on the spectral properties established through the LAP. The absence of eigenvalues provides the injectivity required in the Fredholm argument used to prove the LAP.  The Riemann--Lebesgue lemma then yields weak decay of the energy variables. On compact spectral intervals contained in $(0,\infty)$, the compactness of the physical local-energy observation operators allows this weak convergence to be strengthened to strong local decay. Finally, conservation of the total energy permits the removal of the spectral cutoffs for arbitrary fixed initial data in the energy space. This establishes local energy decay as $|t|\to\infty$, rather than merely decay in a time-averaged sense. We emphasize that the absence of eigenvalues alone is insufficient to obtain this conclusion.

The proof of \cref{thm:full_polynomial_decay} provides a quantitative version of this spectral argument. The proof relies on a smooth spectral decomposition of the initial data into low-, intermediate-, and high-frequency components. For $0<\varepsilon\leqslant\frac12$ and $N\geqslant1$, these components are spectrally supported in $[0,2\varepsilon]$, $[\varepsilon,2N]$, and $[N,\infty)$, respectively, where frequency refers to the spectral parameter of $\Lambda_\Omega$. The overlapping intervals allow us to choose smooth spectral cutoffs forming a partition of unity.

The main step is to establish decay on the intermediate-frequency interval. We first obtain a weighted resolvent estimate at every positive frequency using an outgoing decomposition of the harmonic potential and a vertical integration argument (\cref{prop-lambda4}). We refer to \cite{KM13, MR1925354} for a similar outgoing decomposition and related estimates in the study of the water-wave problem in the presence of a submerged body. Combining this estimate with the free resolvent bounds and the perturbation formula yields fractional H\"older continuity of the resolvent with respect to the spectral parameter (\cref{prop-Holder-ROmega}). Stone's formula then represents the intermediate-frequency contribution as an oscillatory integral. After the change of variables $r=\sqrt{g\lambda}$, applying a Fourier decay lemma to the H\"older continuous amplitude gives the intermediate-frequency energy bound
\begin{equation*}
   C(1+|t|)^{-2\beta}\varepsilon^{-(16+2\beta)}N^{12+\beta}
\mathcal N_s(v_0,v_1)^2,
\end{equation*}
see, \cref{prop:local-spectral-estimate}.

The low- and high-frequency contributions are controlled uniformly in time using conservation of the total energy. At low frequencies, a comparison with the square roots of the Neumann Laplacians on the half-lines yields the low-spectral-mass estimate \cref{lem:low_spectral_mass}, providing a bound of order $\varepsilon\log(2/\varepsilon)$. At high frequencies, the additional Sobolev regularity of the initial data gives a bound of order $N^{-1}$. Finally, choosing $\varepsilon$ and $N$ appropriately as functions of time balances these contributions with the intermediate-frequency decay estimate and yields \eqref{eq:optimized_logarithmic_decay}. The spectral cutoffs are used only in the proof and do not impose any spectral localization assumptions on the initial data.

\section{The Dirichlet and the Dirichlet-to-Neumann operator maps on $\Omega$}\label{sec:spaces}

The aim of this section is to rigorously construct the Dirichlet map $D_\Omega$ and the Dirichlet-to-Neumann map, both formally introduced
in Section \ref{sec:main_results}, and to prove their main properties. 
Some of the arguments below are based on our previous work \cite{Ocqueteau2026}
where these operators have been studied in the case when $\Omega=\Omega_0$, $\Gamma=\Gamma_0$,
$\mathcal{E}=\mathcal{E}_0$, where
\begin{equation}\label{def_cu_zero_1}
 \Omega_0=\left\{\begin{bmatrix} x_0\\ y_0\end{bmatrix}\  |\  y_0>0, x_0^2+y_0^2>1\right\},\quad
 \Gamma_0=\left\{\begin{bmatrix} x_0\\ y_0\end{bmatrix}\  |\  y_0>0, x_0^2+y_0^2=1\right\},
 \end{equation}
 \begin{equation}\label{def_cu_zero_2}
 \mathcal E_0=(-\infty,-1)\cup(1,\infty),
\end{equation}
so that the wetted boundary $\Gamma_0$ is the unit semicircle in the upper half-plane $\mathbb{H}$.

We first recall a result which follows from \cite[Proposition 4.1]{Ocqueteau2026}. 

\begin{proposition}\label{prop:reference_solution}
For  $\Omega=\Omega_0$ and $\tilde {\rm v}\in C(\overline{\mathcal{E}_0})$, $\tilde {\rm v}$ bounded, the problem \eqref{Dirichlet_op_pb} has a unique solution $\phi_0=D_{\Omega_0} {\tilde{\rm v}}$ in the class of functions which are continuous and bounded on $\overline\Omega_0$. 
Moreover, $D_{\Omega_0}$ is a linear bounded operator from $C_b(\overline{\mathcal{E}_0})$ to $C_b(\overline\Omega_0)$, where $C_b(\overline{\Omega_0})$
stands for the space of continuous bounded functions on $\overline{\Omega_0}$, endowed with the $L^\infty$ norm. Finally, $\phi_0\in \mathscr{W}_{\Omega_0}$ whenever $\tilde{\rm v}\in H^1(\mathcal{E}_0)$.
\end{proposition}

\begin{remark}\label{rem_extension}
Examining the proof of Proposition 4.1 in \cite{Ocqueteau2026}, it is easily seen that
the solution provided by \Cref{prop:reference_solution} is a classical one. More precisely,
for each $\tilde{\rm v}\in C_b(\overline{\mathcal{E}_0})$ the function $\Phi_0=D_{\Omega_0} \tilde {\rm v}$ extends to a function which is harmonic in the half-plane $\mathbb{H}$  and satisfies, in the classical sense, a homogeneous Neumann boundary condition in $\Gamma_0$.  
\end{remark}

The construction uses the spaces $\mathscr W_\Omega$ and
$\mathscr Z_\Omega$ defined in \eqref{eq:robin_energy_space} and
\eqref{eq:z_space}. Their completeness, established in
\cref{lem:completeness_W_Z}, allows us to apply the Riesz representation theorem.

The first main result of this section is:

\begin{proposition}\label{prop:h1lifting}
For ${\rm v}\in H^1(\mathcal E)$ the boundary value problem \eqref{Dirichlet_op_pb}
admits a unique solution $\phi=D_\Omega {\rm v}$ (in the sense of \eqref{Dirichlet_op_pb_var}) 
in $\mathscr{W}_\Omega$.  This solution is harmonic in $\Omega$,  satisfies the zero Neumann condition on $\Gamma$ in a classical sense, is continuous and bounded on
$\overline\Omega$ and the restriction of $\phi$ to $\mathcal E$ coincides with ${\rm v}$. Moreover, there exists a constant $C>0$ such that
\begin{equation}\label{eq:lifting_bound}
 \|\nabla D_\Omega {\rm v}\|_{L^2(\Omega)}
       +\|D_\Omega {\rm v}\|_{L^\infty(\Omega)}
                    \leqslant C\|{\rm v}\|_{H^1(\mathcal E)} \qquad\qquad({\rm v}\in H^1(\mathcal E)).
\end{equation}
\end{proposition}

\begin{proof}
We first note that there exists $C>0$ such that for every ${\rm v}\in H^1(\mathcal{E})$ there exists ${\rm v}_e\in H^1(\mathbb R)$ with
${\rm v}_e|_{\mathcal E}={\rm v}$ and $\|{\rm v}_e\|_{H^1(\mathbb R)}\leqslant C\|{\rm v}\|_{H^1(\mathcal E)}$.
Choose $\chi_0\in C_c^{\infty}(\mathbb{R})$ with $\chi_0(0)=1$, and set
${\rm V}(x,y)={\rm v}_e(x)\chi_0(y)$ for every $\begin{bmatrix} x\\ y\end{bmatrix}\in \overline{\Omega}$.  Then ${\rm V}\in H^1(\Omega)$ and it has a trace equal to ${\rm v}$ on $\mathcal E$.  

We next remark that $\phi\in \mathscr{W}_\Omega$ is a solution of \eqref{Dirichlet_op_pb}
(in the sense of \eqref{Dirichlet_op_pb_var})
if and only if ${\rm w}=\phi-{\rm V}$ satisfies
${\rm w}\in\mathscr Z_\Omega$ (introduced in \eqref{eq:z_space}) and
\[
 \int_\Omega\nabla {\rm w}\cdot\nabla\eta \, {\rm d}x\, {\rm d}y
       =-\int_\Omega\nabla {\rm V}\cdot\nabla\eta\, {\rm d}x\, {\rm d}y
       \qquad\qquad(\eta\in\mathscr Z_\Omega).
\]
By the Riesz representation theorem there exists a unique ${\rm w}\in\mathscr Z_\Omega$ satisfying the above condition, which implies that indeed \eqref{Dirichlet_op_pb}
admits a unique solution $\phi=D_\Omega {\rm v}$  in $\mathscr{W}_\Omega$.

To prove the other asserted properties of $\phi$ we use Proposition \ref{prop:reference_solution}
and the conformal map $\Theta:\Omega\to\Omega_0$   constructed in \Cref{lemma:conformal_mapping}.
More precisely, we first remark that
\begin{equation}\label{o_compunere}
 \varphi_0:=\phi\circ\Theta^{-1}
\end{equation}
is harmonic in $\Omega_0$. By conformal invariance 
$\nabla\varphi_0\in L^2(\Omega_0)$; its Dirichlet trace is well defined by
Remark~\ref{rem:homogeneous_boundary_traces}, and on $\mathcal E_0$ it coincides with
$\tilde{\rm v}:={\rm v}\circ\Theta^{-1}$. Consequently, the function $\varphi_0-D_{\Omega_0} \tilde {\rm v}$,
where $D_{\Omega_0}$ is the operator in  Proposition \ref{prop:reference_solution}, is harmonic on $\Omega_0$
and lies in  $\mathscr{Z}_{\Omega_0}$, so it vanishes. It follows that 
\[
\varphi_0=\phi_0=D_{\Omega_0} \tilde {\rm v},
\]
where $D_{\Omega_0}$ is the operator introduced in Proposition \ref{prop:reference_solution}. 
From this and \eqref{o_compunere} it follows that
\[
\phi=\phi_0 \circ \Theta.
\]
All the the other asserted properties of $\phi$ follow now from the properties of $\phi_0$ listed in Proposition \ref{prop:reference_solution}, Remark \ref{rem_extension} and from the properties of $\Theta$ stated and proved  in \Cref{lemma:conformal_mapping}. 
\end{proof}

We next construct the Dirichlet-to-Neumann operator $\Lambda_\Omega$.  We first record a consequence of  a result of Brown \cite[Theorem 2.1]{Brown1994}.

\begin{lemma} \label{lem:brown_traces_green}
Let $G\subset\mathbb R^2$ be a bounded connected Lipschitz domain
whose boundary consists of finitely many regular $C^1$ arcs, with
one-sided tangents at their endpoints. Let $D$ and $N$ be disjoint
relatively open unions of boundary arcs, with $D\ne\varnothing$,
whose closures cover $\partial G$. Assume that every junction between
$D$ and $N$ has an interior angle strictly smaller than $\pi$.
For $f_D\in H^1(D)$ and $f_N\in L^2(N)$, let $h\in H^1(G)$ be the
variational solution of
\begin{equation}
\left\{
\begin{aligned}
 &\Delta h(x,y)=0 && \left(\begin{bmatrix}x\\ y\end{bmatrix} \in G\right)\\
 &h(x,y)=f_D(x,y) && \left(\begin{bmatrix}x\\ y\end{bmatrix} \in D\right)\\
 &\partial_n h(x,y)=f_N(x,y) && \left(\begin{bmatrix}x\\ y\end{bmatrix} \in N\right).
 \end{aligned}
 \right.
\end{equation}
Then, the gradient $\nabla h$ admits an almost-everywhere non-tangential boundary trace
$q\in L^2(\partial G)^2$, and the outward normal derivative $\partial_n h=q\cdot n$
satisfies
\begin{equation}\label{eq:brown_trace_bound}
 \|\partial_n h\|_{L^2(\partial G)}
 \leqslant\|q\|_{L^2(\partial G)}
 \leqslant C\bigl(\|f_D\|_{H^1(D)}+\|f_N\|_{L^2(N)}\bigr),
\end{equation}
where $C>0$ depends only on $G$ and its boundary partition.
Moreover, $\partial_n h$ coincides with the weak normal trace of $h$,
and the following generalized Green's formula holds:
\begin{equation}\label{eq:brown_green_formula}
 \int_G\nabla h\cdot\overline{\nabla w}\,{\rm d}x\,{\rm d}y
 =\int_{\partial G}\partial_n h\,
                       \overline{w}\,{\rm d}\sigma
 \qquad(w\in H^1(G)).
\end{equation}
\end{lemma}

\begin{proof}
We verify that $G$, $D$, and $N$ satisfy the geometric hypotheses required in \cite{Brown1994}, namely, that $\partial G$ is locally represented as a Lipschitz graph separated into Dirichlet and Neumann regions by a Lipschitz interface $\psi$, with uniform $L^\infty$ bounds on the gradients of both functions and with the slope of the graph changing sign across each junction. At any such junction, let $\omega\in(0,\pi)$ denote the interior angle. By rotating and translating coordinates so that the vertical axis bisects the corner, the tangent wedge is given by the epigraph of $x_1\mapsto \cot(\omega/2)|x_1|$. Since $\omega<\pi$, we have $\cot(\omega/2)>0$, and by the continuity of the one-sided tangents, the bounding Lipschitz function has a strictly positive derivative on the Neumann side and a strictly negative derivative on the Dirichlet side. 

Consequently, \cite[Theorem 2.1]{Brown1994} applies and guarantees that $\nabla h$ possesses an almost-everywhere non-tangential boundary limit $q\in L^2(\partial G)^2$, and that the non-tangential maximal function $(\nabla h)^*$ satisfies
\[
 \|(\nabla h)^*\|_{L^2(\partial G)} \leqslant C\bigl(\|f_D\|_{H^1(D)}+\|f_N\|_{L^2(N)}\bigr).
\]
Since $|q|\leqslant (\nabla h)^*$ almost everywhere on $\partial G$, the estimate \eqref{eq:brown_trace_bound} follows. Finally, the generalized Green's formula \eqref{eq:brown_green_formula} follows from the result of Dautray and Lions \cite[Chapter VII, Corollary 1, p. 382]{dlv2}.
\end{proof}

Fix $R>0$ so large that the set $\mathcal B$ defined in \eqref{immersed} satisfies $\overline{\mathcal B}\subset B_R$,
and set
\begin{equation}\label{eq:truncated_fluid_domain}
 G_R=\Omega\cap B_{2R}.
\end{equation}
The boundary of $G_R$ consists of the free-surface portions, the artificial
circular arc, and $\Gamma$.

\begin{lemma}
\label{lem:dauge}
With $G_R$ defined in \eqref{eq:truncated_fluid_domain}, let $F\in L^2(G_R)$ and let $\psi$ solve
\begin{equation}
\left\{
\begin{aligned}
 &\Delta\psi(x,y)=F(x,y) && \left(\begin{bmatrix}x\\ y\end{bmatrix} \in G_R\right),\\
 &\psi(x,y)=0 &&\left(\begin{bmatrix}x\\ y\end{bmatrix} \in (\mathcal E\cap B_{2R})\cup(\Omega\cap\partial B_{2R})\right),\\
 &\partial_n\psi(x,y)=0 &&\left(\begin{bmatrix}x\\ y\end{bmatrix} \in \Gamma\right).
\end{aligned}
\right.
\end{equation}
Then there exists $\varepsilon>0$, depending only on the geometry of $\Omega$ near $x_\pm$, such that $\psi\in H^{\frac{3}{2}+\varepsilon}(G_R)$ and
\begin{equation}\label{eq:dauge_estimate}
 \|\psi\|_{H^{\frac{3}{2}+\varepsilon}(G_R)}\leqslant C_R\|F\|_{L^2(G_R)}.
\end{equation}
In particular $\partial_n\psi|_{\mathcal E\cap B_{2R}}\in H^\varepsilon(\mathcal E\cap B_{2R})
\hookrightarrow L^2(\mathcal E\cap B_{2R})$, and
\begin{equation}\label{eq:dauge_normal_estimate}
 \|\partial_n\psi\|_{L^2(\mathcal E\cap B_{2R})}
   \leqslant C_R\|F\|_{L^2(G_R)}.
\end{equation}
\end{lemma}

\begin{proof}
The boundary conditions remain mixed at $x_\pm$, even though their
data are homogeneous. Let $\omega_\pm$ denote the  angles at
the contact points. The first positive singular exponent for a
Dirichlet--Neumann sector of opening $\omega_\pm$ is
$\pi/(2\omega_\pm)>\frac{1}{2}$, since $\omega_\pm<\pi$. By Dauge \cite[Corollary~23.5 and Remark~23.6]{Dauge1988}, there exists
$\varepsilon>0$ such that $\psi\in H^{\frac{3}{2}+\varepsilon}(G_R)$ with estimate
\eqref{eq:dauge_estimate}.  The trace theorem on Lipschitz boundaries, as stated in
McLean \cite[Theorem~3.38]{mclean2000strongly}, then gives
$\partial_n\psi\in H^\varepsilon(\mathcal E\cap B_{2R})$, and
since $\varepsilon>0$ this embeds continuously into $L^2$.
\end{proof}

The second main result in this section is:

\begin{proposition}\label{prop:DtN_L2}
Let $\Omega$, $\Gamma$ and $\mathcal{E}$ be as defined in \eqref{def_omega_mare}, \eqref{def_gamma_e} and assume that
they satisfy Assumption \ref{assumption_gamma}.
For ${\rm v}\in H^1(\mathcal E)$ let $\phi$ be the solution of \eqref{Dirichlet_op_pb}
constructed in Proposition \ref{prop:h1lifting}. 
Then the partial derivative
$\partial_y \phi$ has an almost-everywhere non-tangential trace on
$\mathcal E$ that belongs to $L^2(\mathcal E)$ and satisfies
\begin{equation}\label{eq:mixedestimate}
 \left\|\frac{\partial \phi}{\partial y}(\cdot,0)\right\|_{L^2(\mathcal E)}
   \leqslant C\|{\rm v}\|_{H^1(\mathcal E)}.
\end{equation}
Consequently, the operator $\Lambda_\Omega$ defined by
\begin{equation}\label{eq:dtn_pointwise}
 \Lambda_\Omega {\rm v}=-\left.\frac{\partial( D_\Omega {\rm v})}{\partial y}\right|_{\mathcal E}
   \quad\text{a.e.\ on }\mathcal E \qquad\qquad({\rm v}\in H^1(\mathcal{E})),
\end{equation}
lies in $\mathcal{L}\left(H^1(\mathcal{E}),L^2(\mathcal{E})\right)$. 
\end{proposition}

\begin{proof}
 Recall from \Cref{prop:h1lifting}  
that $\phi\in C_b(\overline\Omega)\cap \mathscr{W}_\Omega$ and that $\frac{\partial\phi}{\partial n}=0$, in a classical sense, on $\Gamma$.

Let $R$ and $G_R$ be as in \eqref{eq:truncated_fluid_domain}.
Choose a smooth cutoff function $\chi\in C_c^\infty(\mathbb R^2)$ with
$\chi=1$ on $B_R$ and $\operatorname{supp}\chi\subset B_{2R}$, chosen so that
$\partial_y\chi(x,0)=0$ for every $x\in \mathbb{R}$. We write
\begin{equation}\label{desci_phi}
\phi=\phi_1+\phi_2,
\end{equation}
\begin{equation}\label{form_desco}
\phi_1=\chi\phi,\qquad \phi_2=(1-\chi)\phi,
\end{equation}
so that
\begin{equation}\label{eq:phi1_system}
\left\{
\begin{aligned}
&\Delta\phi_1(x,y)=2\nabla\chi(x,y)\cdot\nabla\phi(x,y)+(\Delta\chi(x,y))\phi(x,y) && \left(\begin{bmatrix}x\\ y \end{bmatrix} \in G_R\right)\\
& \phi_1(x,0)=\chi(x,0) {\rm v}(x) && (x\in \mathcal E\cap B_{2R})\\
 &\phi_1(x,y)=0 &&  \left(\begin{bmatrix}x\\ y \end{bmatrix} \in \Omega\cap\partial B_{2R}\right)\\
 &\partial_n\phi_1(x,y)=0 &&  \left(\begin{bmatrix}x\\ y \end{bmatrix} \in \Gamma\right).
\end{aligned}
\right.
\end{equation}

We next remark that we can write 
\begin{equation}\label{descophi1}
\phi_1=\phi_{11}+\phi_{12},
\end{equation}
where $\phi_{11}$ satisfies
\begin{equation}\label{eq:phi11_system}
\left\{
\begin{aligned}
 &\Delta\phi_{11}(x,y)=0 &&  \left(\begin{bmatrix}x\\ y \end{bmatrix} \in  G_R\right),\\
 &\phi_{11}(x,0)=\chi(x,0) {\rm v}(x)&&  \left(\begin{bmatrix}x\\ y \end{bmatrix} \in \mathcal E\cap B_{2R}\right)\\
 &\phi_{11}(x,y)=0 &&  \left(\begin{bmatrix}x\\ y \end{bmatrix} \in \Omega\cap\partial B_{2R}\right)\\
 &\partial_n\phi_{11}(x,y)=0 &&  \left(\begin{bmatrix}x\\ y \end{bmatrix} \in \Gamma\right).
\end{aligned}
\right.
\end{equation}
Applying \cref{lem:brown_traces_green} with $G=G_R$ and,
up to boundary endpoints,
\[
 D_R=(\mathcal E\cap B_{2R})\cup(\Omega\cap\partial B_{2R}),
 \qquad N_R=\Gamma.
\]
it follows that
\begin{equation}\label{eq:local_homogeneous_estimate}
 \|\partial_n\phi_{11}\|_{L^2(\mathcal E\cap B_{2R})}
 \leqslant C_R\|\chi{\rm v}\|_{H^1(\mathcal{E}\cap B_{2R})}
 \leqslant C_R\|{\rm v}\|_{H^1(\mathcal E)}.
\end{equation}

The remainder $\phi_{12}=\phi_1-\phi_{11}$ satisfies
\begin{equation}\label{eq:phi12_system}
\left\{
\begin{aligned}
 &\Delta\phi_{12}(x,y)=F_1(x,y)
&& \left(\begin{bmatrix}x\\ y \end{bmatrix} \in G_R\right),\\
 &\phi_{12}(x,y)=0 && \left(\begin{bmatrix}x\\ y \end{bmatrix} \in (\mathcal E\cap B_{2R}) \cup(\Omega\cap\partial B_{2R})\right)\\
& \partial_n\phi_{12}(x,y)=0 && \left(\begin{bmatrix}x\\ y \end{bmatrix} \in \Gamma\right),
\end{aligned}
\right.
\end{equation}
where
\begin{equation}\label{DEFF1}
    F_1(x,y)=2\nabla\chi(x,y)\cdot\nabla\phi(x,y)+(\Delta\chi(x,y))\phi(x,y).
\end{equation}
Since $\phi=D_\Omega {\rm v}\in C_b(\overline\Omega)$ and
$\nabla\phi\in L^2(\Omega)$ (\cref{prop:h1lifting}), and
$\chi\in C_c^\infty(\mathbb R^2)$, the source term satisfies
\begin{equation}\label{eq:source_bound}
 \|F_1\|_{L^2(G_R)}\leqslant C_R\|{\rm v}\|_{H^1(\mathcal E)}.
\end{equation}
The system \eqref{eq:phi12_system} has zero Dirichlet data on both
parts of the Dirichlet boundary, so \cref{lem:dauge} applies and
gives
\begin{equation}\label{eq:local_remainder_estimate}
 \|\partial_n\phi_{12}\|_{L^2(\mathcal E\cap B_{2R})}
   \leqslant C_R\|F_1\|_{L^2(G_R)}
   \leqslant C_R\|{\rm v}\|_{H^1(\mathcal E)}.
\end{equation}
Adding \eqref{eq:local_homogeneous_estimate} and
\eqref{eq:local_remainder_estimate}, and using the fact $\phi_1$ vanishes outside 
$B_{2R}$, it follows that
\begin{equation}\label{eq:local_total}
 \|\partial_n\phi_1\|_{L^2(\mathcal E)}\leqslant C_R\|{\rm v}\|_{H^1(\mathcal E)}.
\end{equation}

\smallskip
 Next, we estimate the normal derivative on $\mathcal{E}$ of the function  $\phi_2$
defined in \eqref{form_desco}.
Extending $\phi_2$ by zero inside the obstacle the obstacle $\mathcal{B}$ we have that
\begin{equation}
\left\{
\begin{aligned}
 &\Delta\phi_2(x,y)=-F_1(x,y) &&\left(\begin{bmatrix}x\\ y \end{bmatrix} \in\mathbb H\right)\\
& \phi_2(x,0)={\rm v}_2(x) && (x\in\mathbb{R})
 \end{aligned}
 \right.
\end{equation}
where ${\rm v}_2:=(1-\chi){\rm v}\in H^1(\mathbb R)$, $F_1$ has been defined in \eqref{DEFF1} and we recall that $\mathbb{H}$ is the upper half plane.
The source term $F_1$ has support in $\{R\leqslant|z|\leqslant 2R\}$,
in particular in the strip $0\leqslant y\leqslant 2R$. Since $\phi_2 \in L^\infty(\mathbb{H}),$ 
Fourier transformation in $x$ gives, for $\xi\ne 0$,
\begin{align*}
 \widehat\phi_2(\xi,y)
 ={}&e^{-|\xi|y}\widehat {\rm v}_2(\xi)\\
 &+\frac1{2|\xi|}\int_0^\infty
 \left(e^{-|\xi||y-\eta|}-e^{-|\xi|(y+\eta)}\right)
                           \widehat F_1(\xi,\eta)\,{\rm d}\eta.
\end{align*}
Taking the outward normal derivative $-\partial_y$ at $y=0$ yields
\begin{equation}\label{eq:far_normal_formula}
 \widehat{-\partial_y\phi_2}(\xi,0)
       =|\xi|\widehat {\rm v}_2(\xi)
          -\int_0^{2R}e^{-|\xi|\eta}\widehat F_1(\xi,\eta)\,{\rm d}\eta.
\end{equation}
Since $F_1$ is supported in $y\leqslant 2R$, Cauchy--Schwarz and
Plancherel give
\[
 \biggl\|\int_0^{2R}e^{-|\cdot|\eta}\widehat F_1(\cdot,\eta)\,{\rm d}\eta
         \biggr\|_{L^2(\mathbb{R})}
   \leqslant\sqrt{2R}\|F_1\|_{L^2(\mathbb H)},
\]
so
\begin{equation}\label{eq:far_total}
 \|-\partial_y\phi_2(\cdot,0)\|_{L^2(\mathbb R)}
   \leqslant\|\partial_x {\rm v}_2\|_{L^2(\mathbb{R})}+\sqrt{2R}\|F_1\|_{L^2(\mathbb H)}
   \leqslant C_R\|{\rm v}\|_{H^1(\mathcal E)}.
\end{equation}

\smallskip
Combining \eqref{descophi1}, \eqref{desci_phi}, \eqref{eq:local_total} and \eqref{eq:far_total}
yields the conclusion \eqref{eq:mixedestimate}.

\end{proof}

\section{Self-adjointness and proof of \cref{thm:well_posedness_energy}}\label{sec_self}

In this section, we  continue the study of the Dirichlet-to-Neumann operator which we
introduced in Proposition \ref{prop:DtN_L2} and apply them to obtain well-posedness results for $\eqref{v_eq}$.
We continue to use the notation in the previous sections for the sets $\Omega$, $\Gamma$, $\mathcal{E}$
and for the operators $D_\Omega$ and $\Lambda_\Omega$. 

The following lemma establishes a useful Green type formula, involving functions defined on unbounded domains (more precisely, in $\mathscr W_\Omega$ introduced in \eqref{eq:robin_energy_space} and $H^1(\mathcal E)$).

\begin{lemma}\label{lem:operatorgreen}
For every ${\rm u}\in H^1(\mathcal E)$ and $g\in\mathscr W_\Omega$, we have
\begin{equation}\label{eq:operatorgreen}
 \int_\Omega\nabla D_\Omega {\rm u}\cdot\overline{\nabla g}\, {\rm d}x\, {\rm d}y
   =\int_{\mathcal E}\Lambda_\Omega {\rm u}\,\overline g\, {\rm d}\sigma.
\end{equation}
In particular, setting $g=D_\Omega {\rm v}$ with ${\rm v}\in H^1(\mathcal E)$, 
\begin{equation}\label{eq:form_non-neg}
 \langle\Lambda_\Omega {\rm u},{\rm v}\rangle_{L^2(\mathcal E)}
 =\int_\Omega\nabla D_\Omega {\rm u}\cdot\overline{\nabla D_\Omega {\rm v}}\, {\rm d}x\, {\rm d}y,
\end{equation}
so that $\Lambda_\Omega$ is a symmetric and non-negative operator on $L^2(\mathcal{E})$, with domain $H^1(\mathcal E)$.
\end{lemma}

\begin{proof}
Let $R$ be large enough to contain the obstacle $\mathcal B$ and choose a radial
$\chi\in C_c^\infty(B_{2R})$, equal to one on $B_R$, with $|\nabla\chi|\leqslant C/R$. Using \eqref{eq:brown_green_formula} and setting $\phi=D_\Omega {\rm u}$ it follows that
\[
 \int_{\Omega\cap B_{2R}}\chi \, \nabla \phi\cdot\overline{\nabla g}
 +\int_{\Omega\cap B_{2R}}(\nabla\phi\cdot\nabla\chi)\, \overline g
 =\int_{\mathcal E}\chi \Lambda_\Omega {\rm u}\,
                                      \overline g \, {\rm d}\sigma.
\]
Using the facts that  $\nabla g\in L^2(\Omega)$ and that the trace of $g$ on $\mathcal{E}$, together with $\Lambda_\Omega {\rm u}$, lie in $L^2(\mathcal E)$,
it follows that the first and last integrals converge to the corresponding terms in
\eqref{eq:operatorgreen}. The conclusion of the lemma follows if we show that
 the middle term tends to zero when $R\to \infty$. To check this fact we first introduce
$\mathcal A_R=\Omega\cap(B_{2R}\setminus\overline{B_R})$.
For such $R$, the obstacle lies inside $B_R$, hence
$\mathcal A_R=R\mathcal A$, where
\[
 \mathcal A=\{z\in\mathbb H:1<|z|<2\},\qquad
 S=\left\{\begin{bmatrix}x\\ 0\end{bmatrix}\ \big|\ 1<|x|<2\right\}\subset\partial\mathcal A.
\]
We apply the inequality \eqref{eq:poincare_friedrichs}  to $\tilde g(z)=g(Rz)$ on $\mathcal{A}.$ In dimension two,
\begin{align*}
 \|\tilde g\|_{L^2(\mathcal A)}&=R^{-1}\|g\|_{L^2(\mathcal A_R)},\\
 \|\nabla \tilde g\|_{L^2(\mathcal A)}&=\|\nabla g\|_{L^2(\mathcal A_R)},\\
 \|\tilde g\|_{L^2(S)}&=R^{-\frac{1}{2}}\|g|_{RS}\|_{L^2(RS)}.
\end{align*}
Since $RS=\mathcal E\cap\{R<|x|<2R\}$, this proves
\begin{equation}\label{eq:exterior_poincare}
 \frac1R\|g\|_{L^2(\mathcal A_R)}
 \leqslant C\left(\|\nabla g\|_{L^2(\mathcal A_R)}
       +\frac1{\sqrt R}\|g|\|_{L^2(\mathcal E\cap\{R<|x|<2R\})}\right),
\end{equation}
with $C$ independent of $R$.
Consequently,
\[
 \left|\int_{\mathcal A_R}(\nabla\phi\cdot\nabla\chi)\overline g\right|
 \leqslant C\|\nabla\phi\|_{L^2(\mathcal A_R)}
                         \frac1R\|g\|_{L^2(\mathcal A_R)}\longrightarrow0 \qquad(R\to \infty).
\]
This proves \eqref{eq:operatorgreen}.

Setting $g=D_\Omega {\rm v}$ (which belongs to $\mathscr W_\Omega$ by
\cref{prop:h1lifting}) in \eqref{eq:operatorgreen} directly
recovers \eqref{eq:form_non-neg}, confirming symmetry and non-negativity.
\end{proof}

\begin{remark}\label{late_remark}
For later use, we record the  Poincar\'e inequality used in the above proof in the
following form: if $G\in\{\Omega,\mathbb H\}$ and $R$ is sufficiently
large, then
\begin{equation}\label{eq:energy_controls_local_H1}
 \|f\|_{H^1(G\cap B_R)}\leqslant C_R\|f\|_{\mathscr W_G}
 \qquad\qquad(f\in\mathscr W_G).
\end{equation}
Here, $C_R$ may depend on $R$, unlike the constant in
\eqref{eq:exterior_poincare}.
\end{remark}

We next prove a result playing an essential role in the proof of
the self-adjointness of $\Lambda_\Omega$.

\begin{lemma}\label{lemma:I_Lambda_onto}
    The operator  $\mathbb{I}+\Lambda_\Omega$ maps $H^1(\mathcal{E})$ onto $L^2(\mathcal{E})$.  
\end{lemma}

\begin{proof}
    We first note that, by the Riesz representation theorem, for every $f\in L^2(\mathcal E)$ there exists $\Phi\in \mathscr W_\Omega$ such that
\begin{equation}\label{ne_egale}
a_\Omega (\Phi,\Psi)= \langle f, \Psi(\cdot,0)\rangle_{L^2(\mathcal{E})}
\qquad \qquad(\Psi\in \mathscr W_\Omega),
\end{equation}
where $a_\Omega(\cdot,\cdot)$ is given by \eqref{eq:robin_energy_form}.

On the other hand, since $\Phi$ is clearly harmonic in $\Omega$,  
the   trace of its exterior normal derivative $\frac{\partial\Phi}{\partial n}$
 is defined as an element of $\left[H^1(\partial\Omega)\right]'$ (the dual of
$H^1(\partial\Omega)$ with respect to the pivot space $L^2(\partial\Omega)$)
by the formula 
$$
\left\langle \frac{\partial\Phi}{\partial n},{\rm v}\right\rangle_{\left[H^1(\partial\Omega)\right]',H^1(\partial\Omega)}=
\langle\nabla\Phi,\nabla(D_\Omega {\rm v})\rangle_{L^2(\Omega)} \qquad\qquad({\rm v}\in H^1(\mathcal{E})).
$$
Moreover, by Proposition \ref{prop:h1lifting}, we have $D_\Omega {\rm v}\in 
\mathscr W_\Omega$ for every ${\rm v}\in H^1(\Omega)$.
Comparing the above formula and \eqref{ne_egale} (with $\psi=D_\Omega {\rm v}$) it follows that $\frac{\partial\Phi}{\partial n}\in L^2(\partial\Omega)$, that this trace vanishes on $\partial\Omega\setminus \mathcal{E}$ and that 
$$
\frac{\partial \Phi}{\partial y}(\cdot,0)=f-\Phi(\cdot,0) \qquad\qquad {\rm in}\quad L^2(\mathcal{E}).
$$
It follows that the map $x\mapsto \frac{\partial \Phi}{\partial y}(x,0)$ lies in $L^2(\mathcal{E})$ and
$$
\int_\mathcal{E} \left| \frac{\partial \Phi}{\partial y}(x,0) \right|^2\, {\rm d}x \leqslant
2\|f\|^2_{L^2(\mathcal{E})}.
$$
By \cite[Theorem~A]{Brown1994} we obtain that ${\rm u}=
\Phi(\cdot,0)\in H^1(\mathcal E)$.
Thus ${\rm u}\in \mathcal{D}(\Lambda_\Omega)$ and combining 
\eqref{eq:form_non-neg} and \eqref{ne_egale} it follows that
\[
{\rm u}+\Lambda_\Omega{\rm u}=f,
\]
which ends the proof.
\end{proof}

The main result in this section is:

\begin{theorem}\label{thm:self-adjoint}
Let $\Lambda_\Omega$ be the operator on $L^2(\mathcal E)$, with domain 
$\mathcal{D}(\Lambda_\Omega)=H^1(\mathcal E)$, introduced in Proposition \ref{prop:DtN_L2}.
Then $\Lambda_\Omega$ is self-adjoint, non-negative, and injective. Finally, its
graph norm is equivalent to the $H^1(\mathcal E)$ norm, whereas
\begin{equation}\label{eq:domains}
 \mathcal{D}(\Lambda_\Omega^{\frac{1}{2}})=H^{\frac{1}{2}}(\mathcal E),\qquad
 \|\Lambda_\Omega^{\frac{1}{2}} {\rm v}\|_{L^2(\mathcal{E})}^2
   =\|\nabla (D_\Omega {\rm v})\|_{L^2(\Omega)}^2 \qquad\qquad({\rm v}\in H^1(\mathcal E)).
\end{equation}
\end{theorem}

\begin{proof}
We have already seen in Lemma \ref{lem:operatorgreen} that $\Lambda_\Omega$ is 
symmetric. Furthermore, according to Lemma  \ref{lemma:I_Lambda_onto},
the operator $\mathbb{I}+\Lambda_\Omega$ is onto.  Consequently, see, for instance Tucsnak and Weiss \cite[Proposition 3.2.4]{Obs_Book}, the operator $\Lambda_\Omega$ is self-adjoint.

The non-negativity of $\Lambda_\Omega$ follows from \eqref{eq:form_non-neg}. 
Its injectivity follows from \eqref{eq:form_non-neg} and \eqref{eq:dtn_pointwise}.
Indeed, if $\Lambda_\Omega {\rm u}=0$ then $\|\nabla D_\Omega {\rm u}\|_{L^2(\Omega)}^2
=\langle\Lambda_\Omega {\rm u},{\rm u}\rangle=0$, so $D_\Omega {\rm u}$ is constant. The fact
that its trace on $\mathcal E$ lies in $L^2(\mathcal E)$  forces that
constant to be zero, hence ${\rm u}=0$.

Since $\Lambda_\Omega$ is closed and has domain $H^1(\mathcal E)$,
the closed graph theorem shows that its graph norm and the
$H^1(\mathcal E)$ norm are equivalent.
We denote by $(X_0,X_1)_{\theta,2}$ the real interpolation space
obtained by the $K$-method; for Hilbert couples it coincides, up to
equivalent norms, with the complex interpolation space $[X_0,X_1]_\theta$.
Since $\mathbb I+\Lambda_\Omega$ is self-adjoint and bounded below
by $\mathbb I$, we may apply the interpolation theorem of
Lunardi \cite[Theorem~4.36]{lunardi2018interpolation}. Using
$\mathcal{D}(\mathbb I+\Lambda_\Omega)=\mathcal{D}(\Lambda_\Omega)=H^1(\mathcal E)$,
with equivalent graph norms, we obtain
\[
 \begin{aligned}
 \mathcal{D}((\mathbb I+\Lambda_\Omega)^{\frac12})
   &=(L^2(\mathcal E),\mathcal{D}(\Lambda_\Omega))_{\frac12,2}\\
   &=(L^2(\mathcal E),H^1(\mathcal E))_{\frac12,2}
    =H^{\frac12}(\mathcal E),
 \end{aligned}
\]
with equivalence of norms. The last equality is the usual Sobolev
interpolation identity on each of the two half-lines forming
$\mathcal E$, obtained by even reflection and restriction from $\mathbb R$.
To identify the square-root domains, use the spectral calculus notation
of Appendix~B, with $A=\Lambda_\Omega$.
For ${\rm v}\in L^2(\mathcal E)$, the notation
${\rm d}\|E_{\Lambda_\Omega}(\lambda){\rm v}\|_{L^2(\mathcal{E})}^2$ denotes the
non-negative scalar measure determined by the non-decreasing function
$\lambda\mapsto\|E_{\Lambda_\Omega}([0,\lambda]){\rm v}\|_{L^2(\mathcal{E})}^2$.
Its total mass is $\|{\rm v}\|_{L^2(\mathcal{E})}^2$, so
\begin{equation*}
  \int_{[0,\infty)}(1+\lambda)\,
       {\rm d}\|E_{\Lambda_\Omega}(\lambda){\rm v}\|_{L^2(\mathcal{E})}^2 = \|{\rm v}\|_{L^2(\mathcal{E})}^2+
       \int_{[0,\infty)}\lambda\,
       {\rm d}\|E_{\Lambda_\Omega}(\lambda){\rm v}\|_{L^2(\mathcal{E})}^2.   
\end{equation*}
By \eqref{def_domain_fA}, applied to $f(\lambda)=\sqrt{1+\lambda}$
and $f(\lambda)=\sqrt\lambda$, respectively, finiteness of the
left-hand side characterizes $\mathcal{D}((\mathbb I+\Lambda_\Omega)^{\frac12})$,
whereas finiteness of the last integral characterizes
$\mathcal{D}(\Lambda_\Omega^{\frac12})$. Since ${\rm v}\in L^2(\mathcal E)$,
these conditions are equivalent, and both domains equal
$H^{\frac12}(\mathcal E)$.
Finally, the second identity in \eqref{eq:domains} follows directly
from \eqref{eq:form_non-neg} and
$\langle\Lambda_\Omega {\rm v},{\rm v}\rangle
 =\|\Lambda_\Omega^{\frac12}{\rm v}\|_{L^2(\mathcal{E})}^2$ for
${\rm v}\in \mathcal{D}(\Lambda_\Omega)=H^1(\mathcal E)$.
\end{proof}

We are now in a position to prove the well-posedness of \eqref{v_eq}.

\begin{proof}[Proof of \cref{thm:well_posedness_energy}]
By \cref{thm:self-adjoint}, the operator $A=g\Lambda_\Omega$
is non-negative and self-adjoint on $L^2(\mathcal E)$, with
$\mathcal D(A)=H^1(\mathcal E)$ and $\mathcal D(A^{\frac12})=H^{\frac12}(\mathcal E)$,
with equivalent graph and Sobolev norms.
\cref{prop_abstract_eq_spectral} therefore gives existence,
uniqueness, continuous dependence on the initial data, and the
regularity assertions \eqref{reg_slaba}--\eqref{regtare}.

For initial data in $H^1(\mathcal E)\times H^{\frac12}(\mathcal E)$,
taking the real part of the $L^2(\mathcal E)$ inner product of
\eqref{v_eq} with $v_t$ gives
\begin{multline*}
 \frac{\mathrm d}{\mathrm dt}
 \left(\frac12\|v_t(t,\cdot)\|_{L^2(\mathcal{E})}^2
       +\frac g2\|(\Lambda_\Omega^{\frac12}v)(t,\cdot)\|_{L^2(\mathcal{E})}^2\right)\\
 =\operatorname{Re}\langle v_{tt}(t,\cdot)+g(\Lambda_\Omega v)(t,\cdot) , v_t(t,\cdot)\rangle_{L^2(\mathcal{E})} =0 \qquad\qquad(t\in \mathbb{R}).
\end{multline*}
The regularity in \eqref{regtare} justifies this computation.
Density of $H^1(\mathcal E)\times H^{\frac12}(\mathcal E)$ in
$H^{\frac12}(\mathcal E)\times L^2(\mathcal E)$ and continuous dependence
then give energy conservation for all the stated initial data.
The expression of the energy in terms of the harmonic lifting in
\eqref{eq:totalenergy} follows from \eqref{eq:domains}, whose extension
to $H^{\frac12}(\mathcal E)$ by continuity is detailed in
\cref{lem:local_energy_compactness} below.
\end{proof}

The following lemma is a formulation, in our functional setting, of
the classical absence-of-trapped-modes result under the geometric
condition of \cite{John1950OnTM}. We include a proof because its exponentially weighted vertical integral and the resulting
one-dimensional equation will also be used in
\cref{prop-lambda4}.

\begin{lemma}\label{lem:no_fixed_eigenvalues}
Under Assumption~\ref{ass:geometry},
\begin{equation*}
   \ker(\Lambda_\Omega-\lambda I)=\{0\},
\qquad (\lambda>0). 
\end{equation*}
Consequently, since $\Lambda_\Omega$ is self-adjoint, non-negative, and injective (\cref{thm:self-adjoint}) it has no eigenvalues.
\end{lemma}

\begin{proof}
Fix $\lambda>0$ and suppose that ${\rm u}\in \mathcal{D}(\Lambda_\Omega)$ satisfies $\Lambda_\Omega {\rm u}=\lambda {\rm u}$. Set $\phi=D_\Omega {\rm u}$. Then ${\rm u}\in H^1(\mathcal E)$ and  $\phi\in \mathscr{W}_\Omega$. Moreover, from \eqref{eq:operatorgreen}, we deduce 
\begin{equation}\label{eq:417}
 \int_\Omega\nabla\phi\cdot\overline{\nabla g} \ {\rm d}x\,{\rm d}y
       =\lambda\int_{\mathcal E}{\rm u}\overline{g} \ {\rm d} \sigma.
       \qquad(g \in\mathscr W_\Omega).
\end{equation}
By taking $g = \phi$ yields 
\begin{equation} \label{eq-517}
    \|\nabla \phi\|^2_{L^2(\Omega)} = \lambda \|{\rm u}\|^2_{L^2(\mathcal{E})}.
\end{equation}
Define, for $x \in \mathcal{E}$,
\begin{equation} \label{def-F0}
    F(x):=\frac1\lambda\left(
{\rm u}(x)+\int_0^\infty e^{-\lambda y}\partial_y\phi(x,y)\,{\rm d}y
\right).
\end{equation}
Note that, under \cref{ass:geometry}, $F$ is well defined. Moreover, 
by Cauchy--Schwarz we have $F\in L^2(\mathcal E)$. Integrating by parts, we get 
\begin{equation} \label{def-F}
    F(x)=\int_0^\infty e^{-\lambda y}\phi(x,y)\,{\rm d}y \quad (x \in \mathcal{E}).
\end{equation}
Additionally, differentiation with respect to $x,$ we obtain 
\begin{equation}
    F'(x)=\int_0^\infty e^{-\lambda y}\partial_x \phi(x,y)\,{\rm d}y \qquad (x \in \mathcal{E}).
\end{equation}
Therefore $F \in H^1(\mathcal{E}).$ For $\vartheta\in C_c^\infty(\mathcal E)$, we extended it by zero on $(x_-, x_+).$ We denote this extension by $\vartheta_{\textnormal e}.$ Define $g(x, y) = \vartheta_{\textnormal e}(x)e^{-\lambda y}.$ Note that $g\in\mathscr W_\Omega$, with $\mathscr W_\Omega$ defined in \eqref{eq:robin_energy_space}. Taking  this $g$ in \eqref{eq:417}, yields 
\begin{equation}
 \int_{\mathcal E}\int_0^\infty
e^{-\lambda y}\partial_x\phi(x,y)\overline{\vartheta'(x)}\,{\rm d}y\,{\rm d}x
-\lambda\int_{\mathcal E}\int_0^\infty
e^{-\lambda y}\partial_y\phi(x,y)\overline{\vartheta(x)}\,{\rm d}\,y{\rm d}x
=\lambda\int_{\mathcal E}{\rm u}\,\overline\vartheta\,{\rm d}x.   
\end{equation}
By Fubini's theorem, \eqref{def-F0} and \cref{def-F}, from the above identity we deduce that 
\begin{equation}
\int_{\mathcal E}F'\overline{\vartheta'}
-\lambda^2\int_{\mathcal E}F\overline\vartheta=0
\qquad\mbox{ for all } \vartheta\in C_c^\infty(\mathcal E).
\end{equation}
Consequently, $F''(x)+\lambda^2F(x) = 0$ in the sense of distribution. Since $F \in L^2(\mathcal{E}),$ it follows that $F = 0.$ From the definition of $F$ in \eqref{def-F0}, we get 
\begin{equation*}
    {\rm u}(x)=-\int_0^\infty e^{-\lambda y} \partial_y\phi(x,y)\,{\rm d}y.
\end{equation*}
Applying, Cauchy--Schwarz and \eqref{eq-517}, we deduce 
\begin{equation*}
 \lambda\|{\rm u}\|_{L^2(\mathcal E)}^2 \leqslant  \frac12\|\nabla\phi\|_{L^2(\Omega)}^2
=\frac{\lambda}{2}\|{\rm u}\|_{L^2(\mathcal E)}^2.
\end{equation*}
Thus ${\rm u}=0$, as required.
\end{proof}

The absence of eigenvalues proved above is qualitative. For the
quantitative decay argument, we also need to control how much of a
spatially localized datum can lie near the bottom of the continuous
spectrum. More precisely, we shall estimate the norm of
$\mathbf1_{[0,\varepsilon]}(\Lambda_\Omega)f$ in terms of $\varepsilon$
and the $L^1$ and $L^2$ norms of $f$. Such an estimate is useful because
the very low frequencies cannot be treated by oscillatory estimates
uniformly down to zero.

The idea is to retain only the two quarter-planes lying beneath the
free surface. Removing the remaining part of the fluid and allowing
independent extensions in these quarter-planes can only decrease the
minimum Dirichlet energy. The resulting comparison operator is
explicit: on each half-line it is the square root of the Neumann
Laplacian. The Neumann condition occurs on the artificial vertical
edge of each quarter-plane; it imposes no condition on the value of
the free-surface datum at a contact point. We isolate this geometric
comparison first, and then deduce the low-frequency estimate in
\cref{lem:low_spectral_mass}.

Let $A_N$ be the Neumann Laplacian on $L^2(\mathcal E)$, defined by
\[
 A_N{\rm u}=-{\rm u}'',\qquad
 \mathcal D(A_N)=\left\{{\rm u}\in H^2(\mathcal E):
             {\rm u }'(x_-)={\rm u}'(x_+)=0\right\},
\]
where the endpoint derivatives are taken from within $\mathcal E$.
This operator is non-negative and self-adjoint. We define
\begin{equation}\label{eq:LN_definition}
 L_N=A_N^{\frac12},\qquad \mathcal D(L_N)=H^1(\mathcal E).
\end{equation}
There are no endpoint conditions in $\mathcal D(L_N)$.
It is not difficult to check that $L_N$ is non-negative and self-adjoint
and that $\mathcal D(L_N^{\frac12})=H^{\frac12}(\mathcal E)$.

We denote by $q_N$ and $q_\Omega$ the closed quadratic forms
associated with $L_N$ and $\Lambda_\Omega$, respectively:
\begin{equation}\label{eq:comparison_forms}
 q_N({\rm u})=\|L_N^{\frac12}{\rm u}\|_{L^2(\mathcal E)}^2,
 \qquad
 q_\Omega({\rm u})=\|\Lambda_\Omega^{\frac12}{\rm u}\|_{L^2(\mathcal E)}^2
 \qquad\qquad({\rm u}\in H^{\frac12}(\mathcal E)).
\end{equation}
Both forms have domain $H^{\frac12}(\mathcal E)$, by the preceding
observation and \cref{thm:self-adjoint}.

\begin{lemma}\label{lem:quarter_plane_comparison}
The operator $L_N$ defined by \eqref{eq:LN_definition} satisfies
\[
 \|L_N^{\frac12}{\rm u}\|_{L^2(\mathcal E)}^2
 \leqslant\|\Lambda_\Omega^{\frac12}{\rm u}\|_{L^2(\mathcal E)}^2
 \qquad\qquad({\rm u}\in H^{\frac12}(\mathcal E)).
\]
Consequently, for every $\varepsilon>0$ and $f\in L^2(\mathcal E)$,
\begin{equation}\label{eq:neumann_resolvent_comparison}
 \langle(\Lambda_\Omega+\varepsilon\mathbb I)^{-1}f,f\rangle_{L^2(\mathcal E)}
 \leqslant
 \langle(L_N+\varepsilon\mathbb I)^{-1}f,f\rangle_{L^2(\mathcal E)}.
\end{equation}
\end{lemma}

\begin{proof}
For ${\rm u}\in H^{\frac12}(\mathcal E)$, consider the
$L^2(\mathcal E)$-valued function
\[
 U(\cdot,y)=e^{-yL_N}{\rm u}\qquad\qquad(y\geqslant0).
\]
It converges to ${\rm u}$ in $L^2(\mathcal E)$ as $y\downarrow0$.
Since $L_N^2=A_N$, it satisfies
$U_{yy}=A_NU=-U_{xx}$ for $y>0$, together with the Neumann conditions
on the two artificial vertical edges. These assertions follow from
the spectral calculus, since $e^{-yL_N}{\rm u}\in\mathcal D(A_N)$ for $y>0$.
Moreover, the quadratic form of the Neumann Laplacian gives
\[
 \|L_Na\|_{L^2(\mathcal E)}^2
 =\|a_x\|_{L^2(\mathcal E)}^2
 \qquad\qquad(a\in H^1(\mathcal E)).
\]
We now regard $U$ as a function on the union of the two quarter-planes
\[
 Q=Q_-\cup Q_+=\mathcal E\times(0,\infty),\qquad
 Q_-=(-\infty,x_-)\times(0,\infty),\quad
 Q_+=(x_+,\infty)\times(0,\infty).
\]
Its horizontal boundary is $\mathcal E\times\{0\}$, on which its
trace is ${\rm u}$. The preceding identity and the spectral theorem give,
with $q_N$ as defined in \eqref{eq:comparison_forms},
\begin{equation}\label{eq:LN_form}
 \begin{aligned}
 \int_Q|\nabla U|^2\,{\rm d}x\,{\rm d}y
 &=\int_0^\infty\left(
      \|U_y(\cdot,y)\|_{L^2(\mathcal E)}^2
     +\|L_NU(\cdot,y)\|_{L^2(\mathcal E)}^2\right)\,{\rm d}y\\
 &=2\int_0^\infty\|L_Ne^{-yL_N}{\rm u}\|_{L^2(\mathcal E)}^2\,{\rm d}y
 =\|L_N^{\frac12}{\rm u}\|_{L^2(\mathcal E)}^2=q_N({\rm u}).
 \end{aligned}
\end{equation}
The spectral-calculus computation uses the scalar identity
$2\int_0^\infty t^2e^{-2yt}\,{\rm d}y=t$ for $t\geqslant0$,
integrated against the spectral measure of ${\rm u}$. In particular, $U$ has finite
Dirichlet energy and belongs to $H^1$ on every bounded truncation of $Q$.

We next show that $U$ has the least Dirichlet energy among all
extensions of the boundary datum ${\rm u}$ from $\mathcal E\times\{0\}$
to $Q$. Precisely, we claim that
\begin{equation}\label{eq:quarter_plane_minimum}
 q_N({\rm u})=\int_Q|\nabla U|^2
 =\min_{\substack{\nabla V\in L^2(Q)^2\\V(\cdot,0)={\rm u}}}
       \int_Q|\nabla V|^2.
\end{equation}
Here $V$ ranges over functions on $Q$ with square-integrable gradient
and prescribed horizontal trace ${\rm u}$; the traces are understood
as in \cref{rem:homogeneous_boundary_traces}.
No condition is imposed on $V$ on the vertical edges, and $V$ need
not be harmonic. To prove the claim, take any such $V$ and put $Z=V-U$.
The harmonic equation and the homogeneous Neumann conditions give
\[
 \int_Q\nabla U\cdot\overline{\nabla Z}\,{\rm d}x\,{\rm d}y=0.
\]
For completeness, this identity is valid first for tests with compact
support in $y$ and zero trace at $y=0$, by integration by parts in
$L^2(\mathcal E)$ using the form of $A_N$.
It extends to $Z$ by cutting off at large $y$:
the zero horizontal trace and the one-dimensional Hardy inequality give
\[
 \int_Q\frac{|Z(x,y)|^2}{y^2}\,{\rm d}x\,{\rm d}y
 \leqslant4\int_Q|Z_y(x,y)|^2\,{\rm d}x\,{\rm d}y.
\]
Thus, if $\eta_R(y)=1$ for $y\leqslant R$, vanishes for $y\geqslant2R$,
and satisfies $|\eta_R'|\leqslant C/R$, then
$\nabla(\eta_R Z)\to\nabla Z$ in $L^2(Q)$.
Indeed, the derivative of the cutoff is controlled by the tail of
$\int_Q|Z|^2/y^2$, while the other terms converge by dominated
convergence. The cutoff tests belong to the required form space:
$Z(\cdot,y)=\int_0^y Z_y(\cdot,t)\,{\rm d}t$ in $L^2(\mathcal E)$,
so they are square integrable on every strip of finite height.
It follows that
\[
 \int_Q|\nabla V|^2
 =\int_Q|\nabla U|^2+\int_Q|\nabla Z|^2
 \geqslant q_N({\rm u}).
\]

Together with \eqref{eq:LN_form}, this proves the minimum identity
\eqref{eq:quarter_plane_minimum}.

By Assumption~\ref{ass:geometry}, $Q\subset\Omega$.
For ${\rm u}\in H^1(\mathcal E)$, take $V=(D_\Omega{\rm u})|_Q$.
Its gradient is square integrable on $Q$, and its trace on
$\mathcal E\times\{0\}$ is ${\rm u}$, so it belongs to the class of
functions in \eqref{eq:quarter_plane_minimum}. In particular, its
Dirichlet energy on $Q$ is at least that of $U$; no boundary condition
on the artificial vertical edges needs to be checked.
Using \eqref{eq:domains} and the definitions of $q_N$ and $q_\Omega$
in \eqref{eq:comparison_forms}, we obtain
\begin{equation}\label{eq:quarter_plane_form_comparison}
 q_N({\rm u})\leqslant\int_Q|\nabla D_\Omega {\rm u}|^2
 \leqslant\int_\Omega|\nabla D_\Omega {\rm u}|^2
 =\|\Lambda_\Omega^{\frac12}{\rm u}\|_{L^2(\mathcal E)}^2=q_\Omega({\rm u}).
\end{equation}
The form domains of both operators are $H^{\frac12}(\mathcal E)$,
with equivalent form norms. Hence density of $H^1(\mathcal E)$
in $H^{\frac12}(\mathcal E)$ extends this inequality to all
${\rm u}\in H^{\frac12}(\mathcal E)$.

Finally, for any non-negative self-adjoint operator $T$ with closed
quadratic form $q_T$, completing the square in the form norm gives
\[
 \langle(T+\varepsilon\mathbb I)^{-1}f,f\rangle
 =\sup_{{\rm u}\in\mathcal D(q_T)}
 \left\{2\operatorname{Re}\langle f,{\rm u}\rangle
               -q_T({\rm u})-\varepsilon\|{\rm u}\|_{L^2(\mathcal{E})}^2\right\}.
\]
The supremum is attained at $(T+\varepsilon\mathbb I)^{-1}f$.
Applying this formula to $L_N$ and $\Lambda_\Omega$, whose form
domains coincide, proves \eqref{eq:neumann_resolvent_comparison}.
\end{proof}

The next lemma estimates the part of a datum
$f\in L^1(\mathcal E)\cap L^2(\mathcal E)$ lying in the spectral
interval $[0,\varepsilon]$ of $\Lambda_\Omega$, and quantifies how
its $L^2$ norm tends to zero as $\varepsilon\to0^+$.

\begin{lemma}\label{lem:low_spectral_mass}
Let $P_\varepsilon=\mathbf1_{[0,\varepsilon]}(\Lambda_\Omega)$,
$f\in L^1(\mathcal E)\cap L^2(\mathcal E)$ and $0<\varepsilon\leqslant1$. Then
\begin{equation}\label{eq:low_spectral_mass}
 \|P_\varepsilon f\|_{L^2(\mathcal{E})}^2
  \leqslant C\varepsilon\left(
         \log\left(\frac{2}{\varepsilon}\right)\|f\|_{L^1(\mathcal{E})}^2+\|f\|_{L^2(\mathcal{E})}^2\right).
\end{equation}
\end{lemma}
\begin{proof}
We first estimate the scalar product
$\langle(L_N+\varepsilon\mathbb I)^{-1}f,f\rangle_{L^2(\mathcal E)}$. Set
\[
 f_-(r)=f(x_--r),\qquad f_+(r)=f(x_++r)\qquad\qquad(r>0),
\]
and denote by $\mathcal C$ the unitary cosine transform on
$L^2(0,\infty)$, defined initially on $L^1\cap L^2$ by
\[
 (\mathcal Ca)(\xi)=\sqrt{\frac2\pi}
       \int_0^\infty\cos(r\xi)a(r)\,{\rm d}r
 \qquad\qquad(\xi>0).
\]
On each half-line, the cosine transform diagonalizes $A_N$ as multiplication
by $\xi^2$ in the frequency domain, and $L_N=A_N^{\frac12}$ acts as multiplication by $\xi$.
Consequently, by Plancherel's theorem,
\[
 \langle(L_N+\varepsilon\mathbb I)^{-1}f,f\rangle_{L^2(\mathcal E)}
 =\int_0^\infty
       \frac{|\mathcal Cf_-(\xi)|^2}{\xi+\varepsilon}\,{\rm d}\xi + \int_0^\infty
       \frac{|\mathcal Cf_+(\xi)|^2}{\xi+\varepsilon}\,{\rm d}\xi.
\]
On $(0,1)$, we use
$|\mathcal Cf_\pm(\xi)|\leqslant\sqrt{2/\pi}
\|f_\pm\|_{L^1(0,\infty)}$ and
\[
 \int_0^1\frac{{\rm d}\xi}{\xi+\varepsilon}
 =\log\left(\frac{1+\varepsilon}{\varepsilon}\right)
 \leqslant\log\left(\frac2\varepsilon\right).
\]
On $(1,\infty)$, the bound $(\xi+\varepsilon)^{-1}\leqslant1$
and unitarity of $\mathcal C$ give an upper bound
$\|f_\pm\|_{L^2(0,\infty)}^2$.
Finally,
\begin{gather*}
\|f_-\|_{L^1(0,\infty)}^2+\|f_+\|_{L^1(0,\infty)}^2
 \leqslant\|f\|_{L^1(\mathcal E)}^2,
 \\
 \|f_-\|_{L^2(0,\infty)}^2+\|f_+\|_{L^2(0,\infty)}^2
 =\|f\|_{L^2(\mathcal E)}^2.
\end{gather*}
Consequently,
\[
 \langle(L_N+\varepsilon\mathbb I)^{-1}f,f\rangle_{L^2(\mathcal E)}
 \leqslant \frac2\pi\log\left(\frac2\varepsilon\right)
                  \|f\|_{L^1(\mathcal E)}^2
            +\|f\|_{L^2(\mathcal E)}^2.
\]

We now pass from this resolvent estimate to the spectral projection.
The elementary inequality
\[
 \mathbf1_{[0,\varepsilon]}(\lambda)
 \leqslant\frac{2\varepsilon}{\lambda+\varepsilon}\qquad\qquad(\lambda\geqslant0),
\]
combined with the spectral theorem, the identity
$\|P_\varepsilon f\|_{L^2(\mathcal{E})}^2=\langle P_\varepsilon f,f\rangle$
and \cref{lem:quarter_plane_comparison}, gives
\begin{equation}\label{eq:low_projection_resolvent}
 \begin{aligned}
 \|P_\varepsilon f\|_{L^2(\mathcal E)}^2
 &=\langle P_\varepsilon f,f\rangle_{L^2(\mathcal E)}\\
 &\leqslant2\varepsilon
       \langle(\Lambda_\Omega+\varepsilon\mathbb I)^{-1}f,f\rangle_{L^2(\mathcal E)}\\
 &\leqslant2\varepsilon
       \langle(L_N+\varepsilon\mathbb I)^{-1}f,f\rangle_{L^2(\mathcal E)}.
 \end{aligned}
\end{equation}
Combining this with the preceding bound proves \eqref{eq:low_spectral_mass}.
\end{proof}

\section{The free problem and its resolvent estimates}\label{sec:free_waves}

Having established the basic properties of the obstacle operator, we now study the free problem, which provides the reference operator for the limiting absorption principle in Section~\ref{sec:free}. We first recall its governing equations and classical properties, and then establish the resolvent estimates needed below. In this case the domain $\Omega$ is replaced by  the  upper half-plane \(\mathbb{H}\) defined in \eqref{defhmare},
whereas the velocity potential $\varphi(t,x,y)$ and the free surface downward displacement $\zeta(t,x)$ satisfy:
\begin{equation}\label{eq:free_coupled}
\left\{
\begin{aligned}
&\Delta\varphi(t,x,y)=0 & \qquad\left(\begin{bmatrix}x \\ y \end{bmatrix}\in\mathbb{H}\right),\\
&\zeta_t(t,x) = \frac{\partial\varphi}{\partial y}(t,x,0) & \qquad(x\in \mathbb{R}),\\
&\varphi_t(t,x,0) = g\zeta(t,x) & \qquad(x\in \mathbb{R}).
\end{aligned}
\right.
\end{equation}
As recalled in Section \ref{sec:main_results} , the above system reduces to a single second order evolution equation on the free surface of unknown $v(t,x) = \varphi(t,x,0)$. 
More precisely, the linear free oscillations of the fluid are decribed by the initial value problem
\begin{equation}\label{v_eq_plan}
\left\{
\begin{aligned}
 &v_{tt}(t,x)+g(\Lambda_{\mathbb{H}} v)(t,x)=0&&(t\in\mathbb{R}, x\in \mathbb{R}),\\ &v(0,x)=v_0(x),\quad v_t(0,x)=v_1(x) && (x\in \mathbb{R}),
 \end{aligned}
 \right.
\end{equation}
where $\Lambda_{\mathbb{H}}$, designed as the {\em half-plane Dirichlet to Neumann operator} in the remaining part of this work, is defined by 
\begin{equation}\label{eq:free_dtn_fourier}
    (\Lambda_{\mathbb{H}}{\rm v})(x) = \mathcal{F}_{\xi}^{-1}\left(|\xi| \hat{\rm v}(\xi)\right)(x) \qquad\qquad ({\rm v}\in H^1(\mathbb{R}),\ x\in\mathbb{R}).
\end{equation}
Here $\mathcal F$ is the Fourier transform normalized as in \eqref{eq:fourier_convention}, and $\mathcal{F}_{\xi}^{-1}$ denotes its inverse in the variable $\xi$.
Alternatively (see \cite[Proposition 3.1]{Ocqueteau2026}), $\Lambda_{\mathbb{H}}$ can be defined by
\begin{equation}
    (\Lambda_{\mathbb{H}}{\rm v})(x) = -\frac{\partial (D_{\mathbb{H}}{\rm v})}{\partial y}(x,0)\qquad(x\in\mathbb{R}),
\end{equation}
where $D_{\mathbb{H}}$ is the operator associating to every ${\rm v} \in C_b(\mathbb{R})$
the solution $\Phi=D_{\mathbb{H}}{\rm v}$ of
\begin{equation}
\left\{
\begin{aligned}
&\Delta\Phi(x,y) = 0 && \left(\begin{bmatrix} x \\ y\end{bmatrix} \in \mathbb{H}\right)\\
&\Phi(x,0) = {\rm v} && (x\in \mathbb{R}).
\end{aligned}
\right.
\end{equation}
Using the Fourier transform with respect to $x$ the operator $D_\mathbb{H}$ can be alternatively written:
\begin{equation}\label{def_dir_h}
    \left(\mathcal{F} {D_\mathbb{H}} \eta\right)(\xi,y)=\exp(-|\xi|y) \hat\eta(\xi) \qquad\qquad(\eta\in H^1(\mathbb{R}), \xi\in \mathbb{R},\ y\geqslant 0).
\end{equation}

From \eqref{eq:free_dtn_fourier} it is easily proved that:

\begin{proposition}\label{thm:self-adjoint_plane}
$\Lambda_\mathbb{H}$  defined in \eqref{eq:free_dtn_fourier}  is a self-adjoint, non-negative, and injective  operator on $L^2(\mathbb R )$, with domain 
$\mathcal{D}(\Lambda_\mathbb{H})=H^1(\mathbb R)$. Moreover, its
graph norm is equivalent to the $H^1(\mathbb R)$ norm, whereas
\begin{equation}\label{eq:domains_plan}
 \mathcal{D}(\Lambda_\mathbb{H}^{\frac{1}{2}})=H^{\frac{1}{2}}(\mathbb R).
\end{equation}
\end{proposition}

By combining \Cref{thm:self-adjoint_plane} and \Cref{prop_abstract_eq_spectral} we obtain:
\begin{proposition}\label{thm:well_posedness_energy_plan}
For every
$v_0 \in H^{\frac12}(\mathbb R),\  v_1\in L^2(\mathbb R)$ there exists a
unique solution of \eqref{v_eq_plan}, satisfying
\begin{equation}\label{reg_slaba_plan}
 v\in C(\mathbb R;H^{\frac12}(\mathbb R))\cap C^1(\mathbb R;L^2(\mathbb R)),
\end{equation}
 and this solution depends
continuously on the initial data on bounded time intervals.
If $v_0\in H^1(\mathbb R),\ v_1\in H^{\frac{1}{2}}(\mathbb R)$, then
\begin{equation}\label{regtare_plan}
 v\in C(\mathbb R;H^1(\mathbb R))\cap C^1(\mathbb R;H^{\frac{1}{2}}(\mathbb R))
                         \cap C^2(\mathbb R;L^2(\mathbb R)).
\end{equation}
The conserved energy is
\begin{equation}\label{eq:totalenergy_plan}
 \begin{aligned}
 E(t)&=\frac12\|v_t(t,\cdot)\|_{L^2(\mathbb R)}^2+
           \frac g2\|(\Lambda_{\mathbb{H}}^{\frac{1}{2}}v)(t,\cdot)\|_{L^2(\mathbb R)}^2\\
     &=\frac12\|v_t(t,\cdot)\|_{L^2(\mathbb R)}^2+\frac g2\int_{\mathbb{H}}|\nabla (D_{\mathbb{H}} v)(t,\cdot)|^2.
\end{aligned}
\end{equation}
\end{proposition}

Taking the Fourier transform with respect to $x$, the solution of \eqref{v_eq_plan} writes
\begin{equation}\label{eq:free_solution_fourier}
    v(t,x) = \frac{1}{\sqrt{2\pi}} \int_{\mathbb{R}} e^{ix\xi} \left( \hat{v}_0(\xi)\cos(t\sqrt{g|\xi|}) + \hat{v}_1(\xi)\frac{\sin(t\sqrt{g|\xi|})}{\sqrt{g|\xi|}} \right) {\rm d}\xi.
\end{equation}

By separating the positive and negative frequencies and introducing the change of variables $\kappa = \sqrt{|\xi|}$, we have:
\begin{multline}\label{eq:free_solution_kappa}
    v(t,x) = \sqrt{\frac{2}{\pi}} \int_0^\infty \Bigg[ \kappa \left( e^{ix\kappa^2}\hat{v}_0(\kappa^2) + e^{-ix\kappa^2}\hat{v}_0(-\kappa^2) \right) \cos(t\sqrt{g}\kappa) \\
    + \frac{1}{\sqrt{g}}\left( e^{ix\kappa^2}\hat{v}_1(\kappa^2) + e^{-ix\kappa^2}\hat{v}_1(-\kappa^2) \right) \sin(t\sqrt{g}\kappa) \Bigg] {\rm d}\kappa.
\end{multline}
Assuming that $v_0$ and $v_1$ have piecewise continuous third derivatives in $L^1(\mathbb{R})$ and
integrating by parts in the right-hand side of \eqref{eq:free_solution_kappa}, it has been shown in
\cite[Appendix A]{Weber_1976} that $v(t,x) = O(t^{-1})$, uniformly for $x$ in compact intervals.

The remaining part of this section provides several useful properties of the operator $\Lambda_\mathbb{H}$ defined in \eqref{eq:free_dtn_fourier}. We first note the following result:

\begin{lemma} \label{lem-RHz}
The resolvent set of $\Lambda_\mathbb{H}$ contains \ $\mathbb{C}\setminus [0,\infty)$, so that the resolvents of $\Lambda_\mathbb{H}$ defined by 
\begin{equation}\label{def_r_h}
     \mathcal{R}_{\mathbb{H}}(z) = (\Lambda_{\mathbb{H}} - z)^{-1}\qquad(z\in\mathbb{C}\setminus [0,\infty)),
\end{equation} 
lie in $\mathcal{L}(L^2(\mathbb R))$.
Moreover, for every
    $z\in\mathbb{C}\setminus[0,\infty)$ and every $f\in L^2(\mathbb{R})$
\begin{equation} \label{RHz-formula}
    \mathcal{R}_{\mathbb{H}}(z) f  = \mathcal{F}^{-1}_{\xi}
        \left(\frac{\hat f(\xi)}{|\xi|-z}\right) = T(z)f +  \mathcal{R}_{\mathbb{H}}(-z)f,
\end{equation}
where the second equality holds for $z\in\mathbb{C}\setminus\mathbb{R}$, and for such $z$ the operator $T(z)$ is defined by
\begin{equation}\label{def_op_t}
  \left(T(z)f\right)(x)  = \begin{dcases}
   i\int_{\mathbb{R}} e^{ iz|x - r|} f(r)\, {\rm d}r & \qquad  \mathrm{Im} \;z > 0,\\
   -i\int_{\mathbb{R}} e^{- iz|x - r|} f(r)\, {\rm d}r & \qquad \mathrm{Im} \;z < 0,
  \end{dcases},
\end{equation}
for \(f\in L^2(\mathbb{R})\) and \(x \in \mathbb{R}\).
\end{lemma}

\begin{proof}
The fact that for \ $z\in \mathbb{C}\setminus [0,\infty)$\ the operators defined in \eqref{RHz-formula}
lie in $\mathcal{L}(L^2(\mathbb{R}))$ follows from  Proposition \ref{thm:self-adjoint_plane} and basic spectral theory of self-adjoint operators.

The fact that $\mathcal{R}_{\mathbb{H}}(z)$ is well-defined follows from the fact that $\Lambda_{\mathbb{H}}$ is a self-adjoint positive operator.

The first equality in \eqref{RHz-formula} follows by applying the one dimensional Fourier transform to the equation
\[
\Lambda_\mathbb{H}{\rm v}- z{\rm v}=f\in L^2(\mathbb{R}).
\]
To prove the second equality in \eqref{RHz-formula}, we note that
for every $z\in \mathbb{C}\setminus \mathbb{R}$ we have
    \begin{equation*}
        \frac{1}{|\xi| -z} = \frac{2z}{\xi^2 - z^2} + \frac{1}{|\xi| +z},  
    \end{equation*}
and that for every $x\in \mathbb R$ we have
\begin{equation}
 \mathcal{F}_\xi^{-1} \left ( \frac{2z}{\xi^2 - z^2}\right)(x) =
 \begin{dcases}
 \sqrt{2\pi} \, i e^{i z |x|} & \qquad \qquad{\rm Im}\, z> 0,\\
 -\sqrt{2\pi} \, i e^{ - i z |x|} & \qquad \qquad{\rm Im}\, z < 0.
 \end{dcases}
\end{equation}
\end{proof}

\begin{remark}\label{rem_op_t}
We need below the operators $T^\pm(\lambda)$, with $\lambda\in [0,\infty)$, formally obtained
by passing  to the limit for $z\to \lambda$ in \eqref{def_op_t}. These operators are no longer bounded on 
$L^2(\mathbb{R})$, so that we restrict them to $L_s^2(\mathbb{R})$, with $s>\frac12$. More precisely, for every $\lambda\geqslant 0$ we set 
\begin{equation}\label{def_pm_op}
(T^\pm(\lambda)f)(x)= \pm i\int_{\mathbb{R}} e^{\pm  i\lambda |x - r|} f(r)\, {\rm d}r
\qquad\qquad(f\in L^2_s(\mathbb{R}),\ x\in \mathbb R).
\end{equation}
Since for $s>\frac12$ we have the continuous inclusion \ $L^2_s(\mathbb{R}) \hookrightarrow L^1(\mathbb{R})$ \ it follows that
\begin{equation} \label{eq:RH-sup2}
        \left|\left(T^{\pm}(\lambda)f\right)(x)\right| \leqslant \int_{\mathbb{R}} |f(r)| \ {\rm d}r \leqslant C \|f\|_{L^2_s(\mathbb{R})} \qquad\qquad (f\in L^2_s(\mathbb{R}),x\in \mathbb R).
\end{equation}
From the above estimate it follows, in particular, that $T^{\pm}(\lambda) \in \mathcal{L}\left(L^2_s(\mathbb{R}),L^2_{-s}(\mathbb{R})\right)$ for every $\lambda>0$.
\end{remark}

The next result describes the behaviour of $\mathcal{R}_{\mathbb{H}}(z)$ when $z$ approaches the positive real axis.

\begin{proposition}\label{prop:free-lap}
     For each $s > \frac12$ and $\lambda > 0$, let  $T^{\pm}(\lambda) \in \mathcal{L}\left(L^2_s(\mathbb{R}),L^2_{-s}(\mathbb{R})\right)$ the operators introduced 
    in Remark \ref{rem_op_t}. For each $\lambda > 0$, let $\mathcal{R}^{\pm}_{\mathbb{H}}(\lambda)
    \in \mathcal{L}\left(L^2_s(\mathbb{R}),L^2_{-s}(\mathbb{R})\right)$ be defined by
\begin{equation} \label{eq:RH-boundary}
        \mathcal{R}^{\pm}_{\mathbb{H}}(\lambda) =  T^{\pm}(\lambda) + \mathcal{R}_{\mathbb{H}}(-\lambda).
    \end{equation}
Then
\begin{equation}\label{R_H_limit}
        \lim_{\varepsilon \to 0^+}
        \mathcal{R}_{\mathbb{H}}(\lambda \pm i\varepsilon)=\mathcal{R}^{\pm}_{\mathbb{H}}(\lambda)
        \qquad\qquad(\lambda > 0)
    \end{equation}
    in the norm of
    $\mathcal{L}(L^2_s(\mathbb{R}),L^2_{-s}(\mathbb{R}))$.
    More precisely, for every compact interval $[a,b]\subset(0,\infty)$ and
    every $0<\beta<\min\{s-\frac12,1\}$, there exists
    $C_{s,\beta}>0$ such that
    \begin{equation}\label{eq:free_lap_norm_rate}
        \sup_{\lambda\in [a,b]}
        \big\|\mathcal{R}_{\mathbb{H}}(\lambda\pm i\varepsilon)
        -\mathcal{R}^{\pm}_{\mathbb{H}}(\lambda)\big\|
        _{\mathcal{L}(L^2_s(\mathbb{R}),L^2_{-s}(\mathbb{R}))}
        \leqslant (C_{s,\beta} + a^{-2})\varepsilon^\beta
        \qquad(0<\varepsilon\leqslant1).
    \end{equation}
    
\end{proposition}

\begin{proof}
    
Let $z = \lambda + i  \varepsilon$, with $\lambda,\ \varepsilon > 0$. Fix $\beta$ as in the statement of the proposition. Then for every $x\in \mathbb{R}$ the operators $T(z)$ and $T^+(\lambda)$ defined in 
\eqref{def_op_t} and \eqref{def_pm_op}, respectively, satisfy
\begin{equation}\label{tdeasupra}
        \left| \left(T(z)f\right)(x) - \left(T^+(\lambda) f\right)(x)\right| \leqslant \varepsilon^\beta \int_{\mathbb{R}} |x-r|^\beta|f(r)| \ {\rm d}r.
\end{equation}
    In the above estimate we have used the fact that
\[
0\leqslant 1-e^{-r}\leqslant \min\{r,1\}\leqslant r^\beta \qquad\qquad(r\geqslant 0).
\]
Consequently, for every $s>\frac12$, $z = \lambda + i  \varepsilon$, with $\lambda,\ \varepsilon > 0$ we have 
    \begin{multline} \label{eq-459}
         \left\| T(z) - T^+(\lambda) \right\|^2_{ _{\mathcal{L}(L^2_s(\mathbb{R}),L^2_{-s}(\mathbb{R}))}} \leqslant \varepsilon^{2\beta} \int_{\mathbb{R}^2} \langle x\rangle^{-2s} |x-r|^{2\beta} \langle r\rangle^{-2s} {\rm d}x\, {\rm d}r \\
         \leqslant C_\beta\, \varepsilon^{2\beta} \int_{\mathbb{R}^2} \langle x\rangle^{-2s} \left(\langle x\rangle^{2\beta} + \langle r\rangle^{2\beta}\right) \langle r\rangle^{-2s} {\rm d}x\, {\rm d}r \leqslant C_{s, \beta}\, \varepsilon^{2\beta}.
     \end{multline}
     Similarly we obtain that for $s>\frac12$,  $z = \lambda - i  \varepsilon$, with $\lambda,\ \varepsilon > 0$ and $\beta$ as above we have
\begin{equation} \label{eq-459-jos}
         \left\| T(z) - T^-(\lambda) \right\|^2_{ _{\mathcal{L}(L^2_s(\mathbb{R}),L^2_{-s}(\mathbb{R}))}}  \leqslant C_{s, \beta}\, \varepsilon^{2\beta},
\end{equation}
where $T^-(\lambda)$ has been defied in \eqref{def_pm_op}.

Let $a,\ b\in (0,\infty)$ with $a<b$ and let $\lambda \in [a, b]$ Then, using the fact that $\Lambda_\mathbb{H}$ is self-adjoint and non-negative, combined with the resolvent identity yields
\begin{equation*}
    \|\mathcal{R}_{\mathbb{H}}(-\lambda \pm i \varepsilon) - \mathcal{R}_{\mathbb{H}}(-\lambda) \|_{\mathcal{L}(L^2(\mathbb{R}))} \leqslant a^{-2}\varepsilon \qquad\qquad(\lambda,\varepsilon>0), 
\end{equation*}
so that, 
\begin{equation} \label{eq-460}
     \|\mathcal{R}_{\mathbb{H}}(-\lambda \pm i \varepsilon) - \mathcal{R}_{\mathbb{H}}(-\lambda) \|_{ _{\mathcal{L}(L^2_s(\mathbb{R}),L^2_{-s}(\mathbb{R}))}} \leqslant a^{-2}\varepsilon
     \qquad\qquad\left(s>\frac12, \lambda,\varepsilon>0\right). 
\end{equation}
Combining \eqref{eq:RH-boundary} and \eqref{eq-459}-\eqref{eq-460}, we obtain the conclusion \eqref{eq:free_lap_norm_rate}.
\end{proof}

Using the terms introduced in \Cref{def_lap}, the above proposition (specifically the existence of the limits \eqref{R_H_limit}) is summarized as follows:

\begin{corollary}\label{cor_lap}
For every $s>\frac12$, the operator $\Lambda_\mathbb{H}$ introduced in \eqref{eq:free_dtn_fourier}
 satisfies the limiting absorption principle  with respect to the scale of spaces 
\begin{equation}
L^2_s(\mathbb{R}) \hookrightarrow L^2(\mathbb{R}) \hookrightarrow L^2_{-s}(\mathbb{R}).
\end{equation}
\end{corollary}

The following result establishes the \textit{stable LAP} property for the free operator (see \Cref{def:stable_lap}).

\begin{proposition}\label{prop:stable_lap_free}
Let $s\geqslant1$, $\lambda>0$ and let 
$\mathcal{R}^{\pm}_{\mathbb{H}}(\lambda)
    \in \mathcal{L}\left(L^2_s(\mathbb{R}),L^2_{-s}(\mathbb{R})\right)$ be the operators defined in 
\eqref{eq:RH-boundary}. If $f\in L_s^2(\mathbb R)$ is such that ${\rm v} := \mathcal{R}_\mathbb{H}^\pm(\lambda)f \in L^2(\mathbb{R})$, then ${\rm v} \in H^1(\mathbb{R})$ and $(\Lambda_{\mathbb{H}} - \lambda){\rm v} = f$. Moreover, if
\[
 \mathrm{Im}\langle \mathcal{R}_\mathbb{H}^\pm(\lambda)f,f\rangle_{L_{-s}^2(\mathbb{R}),L_s^2(\mathbb{R})}=0,
\]
then $\mathcal{R}_\mathbb{H}^\pm(\lambda)f\in L^2(\mathbb{R})$. 
\end{proposition}

\begin{proof}
 If ${\rm v} = \mathcal{R}_{\mathbb{H}}^\pm(\lambda)f \in L^2(\mathbb{R})$, then $f + (\lambda+1){\rm v} \in L^2(\mathbb{R})$, so the identity ${\rm v} = \mathcal{R}_{\mathbb{H}}(-1)[f + (\lambda+1){\rm v}]$ (justified by \cref{prop:free-lap} as in \cref{lem:RH-boundary-extension}) implies ${\rm v} \in \mathcal{D}(\Lambda_{\mathbb{H}}) = H^1(\mathbb{R})$ and $(\Lambda_{\mathbb{H}} - \lambda){\rm v} = f$. Next, we remark that for every $\lambda>0$ and $f\in L^2_s(\mathbb{R})$ we have
\begin{equation}\label{eq:free_imag_formula}
 \mathrm{Im}\langle \mathcal{R}_{\mathbb{H}}^\pm(\lambda)f,f\rangle_{L_{-s}^2(\mathbb{R}),L_s^2(\mathbb{R})} 
   =\pm \frac{1}{2} \left(  \left|\int_{\mathbb R}e^{-i\lambda r}f(r)\,{\rm d}r\right|^2
+\left|\int_{\mathbb R}e^{i\lambda r}f(r)\,{\rm d}r\right|^2 \right).
\end{equation}
Indeed, this follows from two facts. The first one is that 
$\langle \mathcal{R}_{\mathbb{H}}(-\lambda)f, f \rangle_{L^2(\mathbb{R})}$ is real since $\Lambda_{\mathbb{H}}$ is self-adjoint. 
The second one is that for every $\lambda>0$ the operators $T^\pm(\lambda)$ defined
in \eqref{def_pm_op} satisfy
\begin{multline*}
    \mathrm{Im}\langle T^\pm(\lambda)f,f\rangle_{L_{-s}^2(\mathbb{R}),L_s^2(\mathbb{R})} = \pm \int_{\mathbb{R}^2} \cos(\lambda(x-r)) f(r) \overline{f(x)} \; {\rm d}r\, {\rm d}x \\
    = \pm \frac{1}{2} \int_{\mathbb{R}^2} \left(e^{i\lambda(x-r)}
+ e^{-i\lambda(x-r)} \right) f(r) \overline{f(x)} \; {\rm d}r \,{\rm d}x  
= \pm \frac{1}{2} \left(  \left|\int_{\mathbb R}e^{-i\lambda r}f(r)\,{\rm d}r\right|^2
+\left|\int_{\mathbb R}e^{i\lambda r}f(r)\,{\rm d}r\right|^2 \right).
\end{multline*}
 so that we have \eqref{eq:free_imag_formula}. Consequently, the conclusion follows if we show that 
 \begin{equation}\label{conc_follows}
     \int_{\mathbb{R}} e^{\pm i \lambda r} f(r) \ {\rm d}r = 0 \implies T^{\pm}(\lambda)f \in L^2(\mathbb{R}). 
 \end{equation}
 Take $x\geqslant 0.$ Assuming that $f$ satisfies \eqref{conc_follows} and using \eqref{def_pm_op}, it follows that 
 \begin{multline*}
(T^{\pm}(\lambda)f )(x) 
=\pm i e^{\pm i\lambda x}
 \int_{-\infty}^{x}e^{\mp i\lambda r}f(r)\,{\rm d}r
\pm i e^{\mp i\lambda x}
 \int_x^\infty e^{\pm i\lambda r}f(r)\,{\rm d}r \\
 = \mp i e^{\pm i\lambda x} \int_x^\infty e^{\mp i\lambda r}f(r)\,{\rm d}r  \pm i e^{\mp i\lambda x}
 \int_x^\infty e^{\pm i\lambda r}f(r)\,{\rm d}r. 
 \end{multline*}
 Consequently, 
 \begin{equation*}
     \left|(T^{\pm}(\lambda)f )(x) \right| \leqslant 2 \int_x^\infty |f(r)| \; {\rm d} r, \qquad (x \geqslant 0).
 \end{equation*}
 In a similar manner, we also have 
  \begin{equation*}
     \left|(T^{\pm}(\lambda)f )(x) \right| \leqslant 2 \int_{-\infty}^x |f(r)| \; {\rm d} r, \qquad (x \leqslant 0).
 \end{equation*}
 We use the dual form of Hardy's inequality, see for instance
 Hardy, Littlewood and P\'{o}lya \cite[Theorem 328]{Hardy}. More precisely, for every $a\in\mathbb R$
 and every measurable function $h$ for which the corresponding right-hand
 side is finite,
 \begin{equation}\label{eq:hardy-tail}
 \begin{aligned}
 \left\|\int_x^\infty |h(r)|\,{\rm d}r\right\|_{L^2(a,\infty)}
 &\leqslant 2\|(x-a)h\|_{L^2(a,\infty)},\\
 \left\|\int_{-\infty}^x |h(r)|\,{\rm d}r\right\|_{L^2(-\infty,a)}
 &\leqslant 2\|(a-x)h\|_{L^2(-\infty,a)}.
 \end{aligned}
 \end{equation}
 The second inequality follows from the first by reflection.
 Applying \eqref{eq:hardy-tail} with $a=0$, we obtain 
 \begin{multline*}
     \left\| T^{\pm} f\right\|_{L^2(\mathbb{R})}^2 \leqslant 4 
\int_0^\infty\left(\int_x^\infty|f(r)|\,{\rm d}r\right)^2 {\rm d}x 
\quad+
4\int_{-\infty}^0\left(\int_{-\infty}^x|f(r)|\,{\rm d}r\right)^2 {\rm d}x\\
\leqslant  16 \int_0^\infty x^2 |f(x)|^2 {\rm d}x + 16 \int_{-\infty}^{0} x^2 |f(x)|^2 {\rm d}x \leqslant C \|f\|^2_{L^2_s(\mathbb{R})} \quad (s \geqslant 1).
 \end{multline*}
 This completes the proof of the proposition. 
\end{proof}

\begin{lemma} \label{prop-holder-RH}
Let $s>\frac12$ and $0<\beta<\min\{s-\frac12,1\}$.
There is a constant $C_{s,\beta}>0$, independent of $0<a<b$, such that
\begin{equation}\label{eq:sup-RH}
 \sup_{\lambda\in[a,b]}
 \|\mathcal R_{\mathbb H}^{\pm}(\lambda)\|
 _{\mathcal L(L_s^2(\mathbb R),L_{-s}^2(\mathbb R))}
 \leqslant C_{s,\beta}\bigl(1+\log(1+a^{-1})\bigr),
\end{equation}
and, for every $\lambda,\mu\in[a,b]$,
\begin{equation}\label{eq:holder-RH}
\|\mathcal R_{\mathbb H}^{\pm}(\lambda)
       -\mathcal R_{\mathbb H}^{\pm}(\mu)\|
 _{\mathcal L(L_s^2(\mathbb R),L_{-s}^2(\mathbb R))}\leqslant C_{s,\beta}
 \bigl(1+\min\{a^{-\beta},a^{-1-\beta}\}\bigr)
 |\lambda-\mu|^\beta.
\end{equation}
\end{lemma}
\begin{proof}
We use the decomposition \eqref{eq:RH-boundary}, namely
\[
 \mathcal R_{\mathbb H}^{\pm}(\lambda)
 =T^\pm(\lambda)+\mathcal R_{\mathbb H}(-\lambda),
\]
where $T^\pm(\lambda)$ is defined by \eqref{def_pm_op} and
$\mathcal R_{\mathbb H}(z)$ by \eqref{def_r_h}.

We first estimate $T^\pm$. Since $s>\frac12$, Cauchy--Schwarz gives
\[
 \|f\|_{L^1(\mathbb{R})}\leqslant
 \left(\int_{\mathbb R}\langle r\rangle^{-2s}\,{\rm d}r\right)^{\frac{1}{2}}
 \|f\|_{L_s^2(\mathbb{R})},\qquad
 \|h\|_{L_{-s}^2(\mathbb{R})}\leqslant
 \left(\int_{\mathbb R}\langle x\rangle^{-2s}\,{\rm d}x\right)^{\frac{1}{2}}
 \|h\|_{L^\infty(\mathbb{R})}.
\]
Thus \eqref{eq:RH-sup2} implies the uniform estimate
\[
 \|T^\pm(\lambda)\|_{\mathcal L(L_s^2(\mathbb{R}),L_{-s}^2(\mathbb{R}))}\leqslant C_s
 \qquad(\lambda>0).
\]
For $\lambda,\mu>0$, the elementary inequality
$|e^{iu}-e^{iv}|\leqslant 2^{1-\beta}|u-v|^\beta$ yields
\[
 |((T^\pm(\lambda)-T^\pm(\mu))f)(x)|
 \leqslant 2^{1-\beta}|\lambda-\mu|^\beta
 \int_{\mathbb R}|x-r|^\beta|f(r)|\,{\rm d}r.
\]
Applying Cauchy--Schwarz in $r$ and then integrating in $x$, we obtain
\[
 \begin{aligned}
 &\|T^\pm(\lambda)-T^\pm(\mu)\|_{\mathcal L(L_s^2(\mathbb{R}),L_{-s}^2(\mathbb{R}))}^2\\
 &\quad\leqslant C_\beta|\lambda-\mu|^{2\beta}
 \int_{\mathbb R^2}\langle x\rangle^{-2s}|x-r|^{2\beta}
                  \langle r\rangle^{-2s}\,{\rm d}x\,{\rm d}r\\
 &\quad\leqslant C_{s,\beta}|\lambda-\mu|^{2\beta}.
 \end{aligned}
\]
Indeed, $|x-r|^{2\beta}\leqslant
C_\beta(\langle x\rangle^{2\beta}+\langle r\rangle^{2\beta})$,
and the resulting integrals are finite because $s>\frac12$ and
$s-\beta>\frac12$.

We next estimate the regular part $\mathcal R_{\mathbb H}(-\lambda)$.
It is differentiable in $\mathcal L(L^2(\mathbb{R}))$ for $\lambda>0$, with
\[
 \frac{{\rm d}}{{\rm d}\lambda}\mathcal R_{\mathbb H}(-\lambda)
 =-\mathcal R_{\mathbb H}(-\lambda)^2.
\]
By the Fourier representation \eqref{RHz-formula}, the operator
$\mathcal R_{\mathbb H}(-\lambda)^2$ is convolution with
\[
 k_\lambda(x)=\frac1{2\pi}\int_{\mathbb R}
          \frac{e^{ix\xi}}{(|\xi|+\lambda)^2}\,{\rm d}\xi.
\]
Since
\[
 \|k_\lambda\|_{L^\infty(\mathbb{R})}
 \leqslant\frac1{2\pi}\int_{\mathbb R}
                   \frac{{\rm d}\xi}{(|\xi|+\lambda)^2}
 =\frac1{\pi\lambda},
\]
the two embeddings established above imply
\[
 \left\|\frac{{\rm d}}{{\rm d}\lambda}
          \mathcal R_{\mathbb H}(-\lambda)\right\|
 _{\mathcal L(L_s^2(\mathbb{R}),L_{-s}^2(\mathbb{R}))}
 \leqslant\frac{C_s}{\lambda}.
\]
For $0<\lambda\leqslant1$, integration from $\lambda$ to $1$ gives
\[
 \|\mathcal R_{\mathbb H}(-\lambda)\|_{\mathcal L(L_s^2(\mathbb{R}),L_{-s}^2(\mathbb{R}))}
 \leqslant 1+C_s\log\left(\frac{1}{\lambda}\right).
\]
Here we used $\|\mathcal R_{\mathbb H}(-1)\|_{\mathcal L(L^2(\mathbb{R}))}\leqslant1$
and the norm-one inclusions $L_s^2(\mathbb{R})\hookrightarrow L^2(\mathbb{R})\hookrightarrow L_{-s}^2(\mathbb{R})$.
For $\lambda\geqslant1$, non-negativity and self-adjointness give
$\|\mathcal R_{\mathbb H}(-\lambda)\|_{\mathcal L(L^2(\mathbb{R}))}
\leqslant\lambda^{-1}\leqslant1$.
Combining these estimates with the uniform bound for $T^\pm$ proves
\eqref{eq:sup-RH}.

To estimate differences, we may assume $a\leqslant\lambda\leqslant\mu$
and set $h=\mu-\lambda$. Integrating the same derivative estimate yields
\[
 \begin{aligned}
 \|\mathcal R_{\mathbb H}(-\lambda)
        -\mathcal R_{\mathbb H}(-\mu)\|_{\mathcal L(L_s^2(\mathbb{R}),L_{-s}^2(\mathbb{R}))}
 &\leqslant C_s\log\left(\frac{\mu}{\lambda}\right)\\
 &\leqslant C_s\log\left(1+\frac{h}{a}\right)
 \leqslant C_{s,\beta}a^{-\beta}h^\beta.
 \end{aligned}
\]
The last inequality follows from
$\log(1+t)\leqslant C_\beta t^\beta$ for $t\geqslant0$.
We also retain the estimate obtained directly from the Fourier multiplier:
\[
 \begin{aligned}
 \|\mathcal R_{\mathbb H}(-\lambda)
        -\mathcal R_{\mathbb H}(-\mu)\|_{\mathcal L(L_s^2(\mathbb{R}),L_{-s}^2(\mathbb{R}))}
 &\leqslant\sup_{\xi\in\mathbb R}
       \frac{h}{(|\xi|+\lambda)(|\xi|+\mu)}\\
 &=\frac{h}{\lambda\mu}
 \leqslant\min\{ha^{-2},a^{-1}\}
 \leqslant a^{-1-\beta}h^\beta.
 \end{aligned}
\]
Here $h/(\lambda\mu)=\lambda^{-1}-\mu^{-1}\leqslant a^{-1}$,
and we used $\min\{u,v\}\leqslant u^\beta v^{1-\beta}$ for $u,v\geqslant0$.
Taking the better of these two bounds for the regular part and adding
the estimate for $T^\pm(\lambda)-T^\pm(\mu)$ proves
\eqref{eq:holder-RH}.
\end{proof}

\begin{lemma} \label{lem-local-estimate-free}
      Let $\frac12<s<\frac32$, and let $D_\mathbb{H}$ be the Dirichlet map defined in \eqref{def_dir_h}.  Then for any $R > 0$ there exists a constant $C_{R, s} > 0$ such that 
      \begin{equation} \label{eq:1205}
          \left\| D_{\mathbb{H}}\mathcal{R}_{\mathbb{H}}(-1) f \right\|_{H^1(\mathbb{H} \cap B_R)} \leqslant C_{R,s}\left\| f\right\|_{L^2_{-s}(\mathbb{R})} \qquad\qquad(f\in L^2(\mathbb{R})).
      \end{equation}
     Moreover, there exists a constant $C_s > 0$ such that 
     \begin{equation} \label{eq:1206}
         \left\| \mathcal{R}_{\mathbb{H}}(-1) f \right\|_{L^2_{-s}(\mathbb{R})} \leqslant C_{s}\left\| f\right\|_{L^2_{-s}(\mathbb{R})} \qquad\qquad(f\in L^2(\mathbb{R})).
     \end{equation} 
\end{lemma}
\begin{proof}
    Let us set $\psi = D_{\mathbb{H}} \mathcal{R}_{\mathbb{H}}(-1) f$. By combining \eqref{def_dir_h}
and \eqref{RHz-formula} it follows that
\begin{equation*}
    \hat \psi(\xi,y)=\frac{\exp(-|\xi|y)}{1+|\xi|} \hat f(\xi) \qquad\qquad(\xi\in \mathbb{R},\ y\geqslant 0).
\end{equation*}
Consequently, $\psi \in \mathscr{W}_{\mathbb{H}}$ and there exists an absolute constant $C>0$
such that
    \begin{equation} \label{eq-920}
        \|\nabla \psi\|_{L^2(\mathbb{H})} \leqslant C \|f\|_{L^2(\mathbb{R})} \qquad\qquad(f\in L^2(\mathbb{R})). 
    \end{equation}
    Moreover, $\psi$ satisfies 
    \begin{equation}  \label{eq:DHRH-1}
    \left\{
    \begin{aligned}
     &\Delta \psi(x,y) = 0 &&\left(\begin{bmatrix}x \\ y \end{bmatrix} \in \mathbb{H}\right)\\
     &-\partial_y \psi(x,0) + \psi(x,0) = f(x) &&(x\in \mathbb{R}).    
    \end{aligned}
    \right.
    \end{equation}
    Then, by using the Fourier transform, $\psi$ admits the following integral representation
     \begin{equation} \label{eq-formula-u}
        \psi(x,y) = \int_{\mathbb{R}} G(x-r, y) f(r) \; {\rm d}r,
    \end{equation}
    where
    \begin{equation} \label{def-G}
        G(x, y) = \frac{1}{\pi} \int_{0}^{\infty} e^{-r} \frac{y+r}{x^2 + (y+r)^2} \; {\rm d}r.
    \end{equation}
    For $|x|\geqslant 1$ and $0\leqslant y \leqslant R$, direct differentiation shows that
\begin{equation} \label{decay-G}
  |G(x,y)|+|\partial_x G(x,y)| + |\partial_y G(x,y)| \leqslant C_R|x|^{-2}.  
\end{equation}
Indeed,
\begin{equation*}
    |G(x,y)| \leqslant \frac{1}{\pi |x|^2} \int_0^\infty e^{-r} (y+r) \, {\rm d} r \leqslant \frac{R+1}{\pi |x|^2} \quad (  |x|\geqslant 1, \; 0\leqslant y \leqslant R)
\end{equation*}
Differentiating $G$ with respect to $x$, we obtain 
\begin{equation*}
    \partial_x G(x,y) =  -\frac2\pi\int_0^\infty
e^{-r}\frac{x(y+r)} {[x^2+(y+r)^2]^2}\,{\rm d}r,
\end{equation*}
so that 
\begin{equation*}
   | \partial_x G(x,y)| \leqslant \frac{2}{\pi |x|^3} \int_0^\infty e^{-r} (y+r)\, {\rm d} r \leqslant \frac{2(R+1)}{\pi |x|^2} \quad (|x|\geqslant 1, \; 0\leqslant y \leqslant R).
\end{equation*}
We can estimate $\partial_yG(x,y)$ in a similar manner. Write
\begin{equation*}
    f = f_1 + f_2, \quad f_1 = f \mathbf{1}_{\{|x|\leqslant 2R +2\}}, \quad f_2 = f \mathbf{1}_{\{|x| > 2R +2\}},
\end{equation*}
and denote the corresponding solution by $\psi_1$ and $\psi_2$. Since $f_1$ is compactly supported, applying \eqref{eq-920} to $\psi_1$ we get 
 \begin{equation*} 
        \|\psi_1\|_{H^1(\mathbb{H} \cap B_R)} \leqslant C \|f_1\|_{L^2(\mathbb{R})} \leqslant C_{R,s} \left\| f\right\|_{L^2_{-s}(\mathbb{R})}. 
    \end{equation*}
To estimate $\psi_2,$ we use \eqref{eq-formula-u}. If $[x,y]^\top \in \mathbb{H} \cap B_R$ and $|r| > 2R+2,$ then $\displaystyle |x-r| > \frac{|r|}{2}.$ By \eqref{decay-G} and Cauchy-Schwarz we get, for any $[x,y]^\top \in \mathbb{H} \cap B_R$
\begin{multline*}
    |\psi_2(x,y)| + |\nabla \psi_2(x,y)| \leqslant C_R \int_{|r| > 2R+2} |r|^{-2} |f(r)| \; {\rm d}r \\
    \leqslant C_{R} \|f\|_{L^2_{-s}(\mathbb{R})} \left( \int_{|r| > 2R+2} |r|^{-4} \langle r\rangle^{2s} \; {\rm d}r \right)^{\frac{1}{2}} \leqslant C_{R,s} \|f\|_{L^2_{-s}(\mathbb{R})}. 
\end{multline*}
In the above estimate we have used $s < \frac{3}2$ so that $2s-4 <-1.$ The above estimate together with the estimate of $\psi_1$ gives \eqref{eq:1205}. 

We now  prove \eqref{eq:1206}. Note that, from \eqref{eq-formula-u}, for any $x \in \mathbb{R}$
\begin{multline*}
    \langle x \rangle^{-s} \mathcal{R}_{\mathbb{H}}(-1) f(x) = \int_{\mathbb{R}} G(x-r, 0) \frac{\langle r \rangle^{s}}{\langle x \rangle^{s}} \langle r \rangle^{-s} f(r) \;{\rm d}r \\
    = \int_{\langle r \rangle \leqslant 2 \langle x \rangle} G(x-r, 0) \frac{\langle r \rangle^{s}}{\langle x \rangle^{s}} \langle r \rangle^{-s} f(r) \;{\rm d}r + \int_{\langle r \rangle > 2 \langle x \rangle} G(x-r, 0) \frac{\langle r \rangle^{s}}{\langle x \rangle^{s}} \langle r \rangle^{-s} f(r) \;{\rm d}r : = I_1(x) + I_2(x).
\end{multline*}
It is easy to see that 
\begin{equation*}
    |I_1(x)| \leqslant 2^s \int_{\mathbb{R}} G(x-r, 0) \langle r \rangle^{-s} |f(r)| \; {\rm d} r.
\end{equation*}
Young inequality for convolution gives, 
\begin{equation*}
    \|I_1\|_{L^2(\mathbb{R})} \leqslant 2^s \|G(\cdot, 0)\|_{L^1(\mathbb{R})} \|\langle \cdot \rangle^{-s} f \|_{L^2(\mathbb{R})} \leqslant 2^s \|f\|_{L^2_{-s}(\mathbb{R})}.
\end{equation*}
When $\langle r \rangle > 2 \langle x \rangle,$ we have $|x -r| \geqslant \frac{1}{2} \langle r \rangle.$ This yields, whenever $\langle r \rangle > 2 \langle x \rangle,$
\begin{equation*}
    |G(x-r,0)| \leqslant \frac{1}{\pi} \int_0^\infty \frac{ze^{-z}}{(x-r)^2 + z^2} \; {\rm d}z \leqslant \frac{1}{\pi (x-r)^2} \int_0^\infty z e^{-z} \;{\rm d}z \leqslant \frac{4}{\pi} \langle r \rangle^{-2}.
\end{equation*}
Therefore, by Cauchy-Schwarz
\begin{multline*}
    |I_2(x)|^2 \leqslant C \langle x \rangle^{-2s} \left| \int_{\mathbb{R}} \langle r \rangle^{s-2} \langle r \rangle^{-s} f(r) \; {\rm d} r \right|^2 \\ \leqslant C \langle x \rangle^{-2s} \|f\|_{L^2_{-s}(\mathbb{R})}^2 \int_{\mathbb{R}} \langle r \rangle^{2s-4} \;{\rm d}r \leqslant C  \langle x \rangle^{-2s} \|f\|_{L^2_{-s}(\mathbb{R})}^2. 
\end{multline*}
Consequently, 
\begin{equation*}
    \|I_2\|^2_{L^2(\mathbb{R})} \leqslant C\|f\|_{L^2_{-s}(\mathbb{R})}^2 \int_{\mathbb{R}} \langle x \rangle^{-2s} \; {\rm d}x \leqslant  C\|f\|_{L^2_{-s}(\mathbb{R})}^2.
\end{equation*}
The above estimate together with the estimate of $I_1$ gives \eqref{eq:1206}. 
\end{proof}

\begin{lemma} \label{lem:RH-boundary-extension}
Let $\frac12<s<\frac32$, and let $\iota:L^2_s(\mathbb{R})\hookrightarrow L^2_{-s}(\mathbb{R})$ denote  the natural inclusion map.  Then 
\begin{equation}
\mathcal R_{\mathbb H}^{\pm}(\lambda) = \mathcal{R}_{\mathbb H}(-1)
\left[\iota+(\lambda+1)\mathcal R_{\mathbb H}^{\pm}(\lambda)\right] \quad (\lambda > 0)
\end{equation}
holds in $\mathcal{L}(L^2_s(\mathbb{R}), L^2_{-s}(\mathbb{R})).$ 
\end{lemma}

\begin{proof}
   We first note that, for any  $z \in \mathbb{C} \setminus [0, \infty),$ 
    \begin{equation*}
\mathcal{R}_{\mathbb{H}}(-1)\left[ \mathbb{I} + (z+1) \mathcal{R}_{\mathbb{H}}(z)\right] = \mathcal{R}_{\mathbb{H}}(z),
    \end{equation*}
    holds in $\mathcal{L}(L^2(\mathbb{R})).$ The conclusion of the lemma follows from \cref{prop:free-lap}  and \eqref{eq:1206}. 
\end{proof}

More precisely, we give below another representation of $\mathcal{R}_{\mathbb{H}}^{\pm}(\lambda)$, which is used afterwards to obtain decay estimates. We show that $\mathcal{R}_{\mathbb{H}}^{+}(\lambda)$ and $\mathcal{R}_{\mathbb{H}}^{-}(\lambda)$ admit a decomposition into outgoing and incoming waves respectively on the two half-lines and an
$H^1(\mathbb R)$ remainder. The wave amplitudes are determined
by the Fourier transform of the datum at $\pm\lambda$.

\begin{lemma}\label{lem:repRH-boundary}
Let $\lambda>0,$ $1\leqslant s < \frac32,$ and $f\in L_s^2(\mathbb{R})$. Fix $r_0>0$ and smooth cutoffs $\chi_{+}$ and $\chi_{-}$ with
\begin{equation}
\begin{dcases}
     &  \chi_+(x)=0 \mbox{ for }x\leqslant r_0, \quad \chi_+(x)=1\mbox{  for }x\geqslant 2r_0,\\
 & \chi_-(x)=1 \mbox{ for }x\leqslant -2r_0, \quad \chi_-(x)=0\mbox{  for }x\geqslant -r_0.  
\end{dcases}
\end{equation}
Then, for every $\lambda > 0$ and $x \in \mathbb{R}$ we have
\begin{equation}\label{eq:outgoing}
 (\mathcal{R}_{\mathbb{H}}^+(\lambda)f)(x)
 =d_+(\lambda) \chi_+(x)e^{i\lambda x}+d_-(\lambda) \chi_-(x)e^{-i\lambda x} + g_\lambda(x) ,
\end{equation}
with $ g_\lambda \in H^1(\mathbb{R})$,
where $d_{\pm}(\lambda)$ are uniquely determined by
\begin{equation}\label{eq:dpm}
 d_+(\lambda)=i\sqrt{2\pi}\,\widehat{f}(\lambda),\qquad
 d_-(\lambda)=i\sqrt{2\pi}\,\widehat{f}(-\lambda) \qquad\qquad(\lambda > 0).
 \end{equation}
 A similar decomposition also holds for $\mathcal{R}_{\mathbb{H}}^{-}(\lambda)f$  from the relation 
 \begin{equation}
     \mathcal R_{\mathbb H}^{-}(\lambda)f=\overline{\mathcal R_{\mathbb H}^{+}(\lambda)\overline{f}}.
 \end{equation}
In particular, we have $\mathcal{R}_{\mathbb{H}}^\pm(\lambda)f \in C_b(\mathbb{R}).$
\end{lemma}

\begin{proof}
We use the decomposition \eqref{eq:RH-boundary}, namely the fact that
\begin{equation*}
 \mathcal R_{\mathbb H}^+(\lambda)f
=T^+(\lambda)f+ \mathcal R_{\mathbb H}(-\lambda) f \qquad\qquad(\lambda > 0,f\in L^2_s(\mathbb{R})), 
\end{equation*}
where $T^{+}(\lambda)$ has been defined by \eqref{def_pm_op}. For every
$\lambda > 0$ and $x\in \mathbb{R}$ we set 
    \[g_\lambda(x) = (\mathcal R_{\mathbb H}(-\lambda) f)(x) + g_{\lambda,1}(x)\]
    with
    \[g_{\lambda,1}(x) 
    =  (T^+(\lambda)f)(x) -  d_+(\lambda) \chi_+ (x) e^{i\lambda x} - d_-(\lambda) \chi_-(x)e^{-i\lambda x}.\] 
From the definitions of $T^{+}(\lambda)$ and $d_{+}$  
it follows that for every $\lambda>0$ and $x\in \mathbb{R}$ we have 

\begin{align*}
   & (T^+(\lambda)f)(x)- d_{+}(\lambda) e^{i\lambda x} \\
    & =  ie^{i\lambda x}\int_{-\infty}^{x} e^{- i\lambda r} f(r)\, {\rm d}r +
    ie^{-i\lambda x}\int_{x}^{\infty} e^{i\lambda r} f(r)\, {\rm d}r - i e^{i\lambda x} \int_{\mathbb{R}} e^{-i\lambda r} f(r) \,{\rm d}r \\
    & = -ie^{i\lambda x}\int_x^\infty e^{-i\lambda r}f(r)\,{\rm d}r
+ie^{-i\lambda x}\int_x^\infty e^{i\lambda r}f(r)\,{\rm d}r.
\end{align*}
Note that $g_{\lambda,1}(x) =(T^+(\lambda)f)(x)- d_{+}e^{i\lambda x}$ for $x > 2r_0.$  Thus, for $x > 2 r_0$ we have 
\begin{equation*}
    \left|g_{\lambda,1}(x) \right| \leqslant 2 \int_{x}^\infty |f(r)| \,{\rm d}r, \quad 
    \left|g_{\lambda,1}'(x) \right| \leqslant 2\lambda \int_{x}^\infty |f(r)| \,{\rm d}r. 
\end{equation*}
Applying the first inequality in \eqref{eq:hardy-tail} with $a=2r_0$
 to these two bounds gives
\[
 \|g_{\lambda,1}\|_{H^1(2r_0,\infty)}^2
 \leqslant 16(1+\lambda^2)
       \int_{2r_0}^\infty (r-2r_0)^2|f(r)|^2\,{\rm d}r
 \leqslant C_{r_0}(1+\lambda^2)\|f\|_{L_s^2(\mathbb R)}^2,
\]
where the last inequality uses $s\geqslant1$.
On the left half-line, the analogous formulas give
\[
|g_{\lambda,1}(x)|\leqslant2\int_{-\infty}^x|f(r)|\,{\rm d}r,\quad
|g_{\lambda,1}'(x)|\leqslant2\lambda\int_{-\infty}^x|f(r)|\,{\rm d}r.
\]
The second inequality in \eqref{eq:hardy-tail}, with $a=-2r_0$,
therefore gives $g_{\lambda,1}\in H^1(-\infty,-2r_0)$. Moreover, $T^+(\lambda) f$ belongs to $H^1(-2 r_0, 2r_0)$ since $f \in L^2_s(\mathbb{R}) \subset L^1(\mathbb{R})$. Consequently, $g_{\lambda,1} \in H^1(\mathbb{R})$. Finally, note that 
$\mathcal{R}_{\mathbb{H}}(-\lambda)f \in H^1(\mathbb{R})$, see \cref{thm:self-adjoint_plane}. Therefore, $g_\lambda \in H^1(\mathbb{R})$.

Finally, suppose that the same decomposition holds with coefficients
$\widetilde d_\pm$ and a remainder $\widetilde g_\lambda\in H^1(\mathbb R)$.
Subtracting the two decompositions shows that
\[
 (d_+-\widetilde d_+)\chi_+e^{i\lambda x}
 +(d_--\widetilde d_-)\chi_-e^{-i\lambda x}
 \in H^1(\mathbb R)\subset L^2(\mathbb R).
\]
On $(2r_0,\infty)$ this function equals
$(d_+-\widetilde d_+)e^{i\lambda x}$, whose modulus is constant.
Its square integrability forces $d_+=\widetilde d_+$.
The same argument on $(-\infty,-2r_0)$ gives
$d_-=\widetilde d_-$. Thus both coefficients, and hence the remainder
for the fixed cutoffs, are unique.
\end{proof}

Recall the spaces $\mathscr W_G$ defined in
\eqref{eq:robin_energy_space}, with norms given by
\eqref{eq:robin_energy_form}.

By \cref{lem:repRH-boundary}, $\mathcal R_{\mathbb H}^{\pm}(\lambda)f$
is bounded and continuous on $\mathbb R$. Its Poisson extension is
therefore bounded and continuous on $\overline{\mathbb H}$; see
Axler, Bourdon and Ramey \cite[Theorems 7.3,7.4,7.5]{axler2013harmonic}.
The next lemma lifts the outgoing decomposition \eqref{eq:outgoing}
to the half-plane, with vertical profile $e^{-\lambda y}$ and a
remainder in $\mathscr W_{\mathbb H}$.

\begin{lemma} \label{lem:DH-lift-boundary}
    Let $\lambda > 0,$ $1 \leqslant s <\frac32$ and $f \in L^2_s(\mathbb{R}).$  Let $\chi_\pm$ be  the cutoff functions introduced in \cref{lem:repRH-boundary} and let $D_{\mathbb{H}}$ be the Dirichlet map defined in \eqref{def_dir_h}. Then $\psi_{\lambda} := D_{\mathbb{H}}\mathcal R_{\mathbb H}^+(\lambda)f$ belongs to  $H^1(\mathbb{H} \cap B_R),$ for any $R > 0,$ and satisfies 
  \begin{equation}  \label{eq:1132}
    \left\{
    \begin{aligned}
     &\Delta \psi_{\lambda}(x,y) = 0 &&\left(\begin{bmatrix}x \\ y \end{bmatrix} \in \mathbb{H}\right)\\
     &-\partial_y \psi_{\lambda}(x,0)  - \lambda \psi_{\lambda}(x,0) = f(x) &&(x\in \mathbb{R}),   
    \end{aligned}
    \right.
    \end{equation}
    in the following sense: 
    \begin{equation} \label{eq:422}
\int_{\mathbb H}\nabla\psi_{\lambda}\cdot\overline{\nabla h} \;{\rm d}x\,{\rm d}y
=
\lambda\int_{\mathbb R} (\mathcal{R}_{\mathbb{H}}^{+}(\lambda) f)(x)\,\overline{h}(\cdot,0) \; {\rm d}x + \int_{\mathbb R} f\,\overline{h}(\cdot,0) \; {\rm d}x\qquad(h \in C_c^{\infty}(\overline{\mathbb{H}})).
\end{equation}
 Moreover,  for all $\begin{bmatrix} x\\ y\end{bmatrix} \in \mathbb{H}$ we have  
\begin{equation} \label{eq-343}
    \psi_\lambda(x,y) =  d_+(\lambda)\chi_+(x)e^{i\lambda x -\lambda y}+d_-(\lambda)\chi_-(x)e^{-i\lambda x-\lambda y} + \psi_{\lambda, \rm c}(x, y),
\end{equation}
where $\psi_{\lambda, \rm c} \in \mathscr{W}_{\mathbb{H}},$ and $d_{\pm}(\lambda)$ are uniquely determined by  $d_{\pm}(\lambda) = i\sqrt{2\pi}\widehat{f}(\pm\lambda).$
\end{lemma}

\begin{proof}
For notational simplicity, we omit the subscript $\lambda$ from
$\psi_\lambda$, $\psi_{\lambda,\rm c}$ and $d_\pm(\lambda)$.
Set
\begin{equation}\label{eq:1032}
 {\rm u}=\mathcal R_{\mathbb H}^+(\lambda)f,
 \qquad \widetilde f=f+(\lambda+1){\rm u}.
\end{equation}
By \cref{lem:repRH-boundary},
\[
 {\rm u}={\rm u}_0+{\rm u}_{\rm c},\qquad
 {\rm u}_0(x)=d_+\chi_+(x)e^{i\lambda x}
             +d_-\chi_-(x)e^{-i\lambda x},\qquad
 {\rm u}_{\rm c}\in H^1(\mathbb R),
\]
with $d_\pm$ given by \eqref{eq:dpm}. In particular, ${\rm u}\in C_b(\mathbb R)$,
so its Poisson extension $\psi=D_{\mathbb H}{\rm u}$ is bounded and
continuous on $\overline{\mathbb H}$, with boundary trace ${\rm u}$.

We first justify the local regularity and the weak boundary condition.
By \cref{prop:free-lap,lem:RH-boundary-extension},
$\widetilde f\in L^2_{-s}(\mathbb R)$ and
${\rm u}=\mathcal R_{\mathbb H}(-1)\widetilde f$.
The estimates \eqref{eq:1205} and \eqref{eq:1206} extend the maps
$D_{\mathbb H}\mathcal R_{\mathbb H}(-1)$ and
$\mathcal R_{\mathbb H}(-1)$ continuously from $L^2(\mathbb R)$ to
$L^2_{-s}(\mathbb R)$, with values in
$H^1(\mathbb H\cap B_R)$ and $L^2_{-s}(\mathbb R)$, respectively.
These extensions are consistent with the Poisson extension:
indeed, for fixed $\begin{bmatrix} x\\ y\end{bmatrix}\in\mathbb H$, the Poisson kernel
$r\mapsto y/[\pi((x-r)^2+y^2)]$ belongs to $L^2_s(\mathbb R)$
because $s<\frac32$. Thus approximation of $\widetilde f$ in
$L^2_{-s}(\mathbb R)$ identifies the local $H^1$ limit with
$D_{\mathbb H}{\rm u}$. Consequently,
$\psi\in H^1(\mathbb H\cap B_R)$ for every $R>0$.
For approximating data in $L^2(\mathbb R)$, the weak Robin identity
follows from \eqref{eq:DHRH-1}. Passing to the limit, using the local
$H^1$ convergence and the compact support of the test function, gives
\begin{equation}\label{eq:420}
 \int_{\mathbb H}\nabla\psi\cdot\overline{\nabla h}\,{\rm d}x\,{\rm d}y
 =\int_{\mathbb R}(\widetilde f-{\rm u})\,
                   \overline{h(x,0)}\,{\rm d}x
 \qquad\qquad(h\in C_c^\infty(\overline{\mathbb H})).
\end{equation}
Since $\widetilde f-{\rm u}=f+\lambda{\rm u}$, this yields \eqref{eq:422}.

We next prove \eqref{eq-343}. Define
\[
 \psi_0(x,y)={\rm u}_0(x)e^{-\lambda y},\qquad
 \phi=D_{\mathbb H}{\rm u}_0-\psi_0,\qquad
 \psi_{\rm c}=D_{\mathbb H}{\rm u}_{\rm c}+\phi.
\]
By linearity, $\psi=\psi_0+\psi_{\rm c}$.
The Fourier representation \eqref{def_dir_h} and Plancherel's theorem give
\[
 \|\nabla D_{\mathbb H}{\rm u}_{\rm c}\|_{L^2(\mathbb H)}^2
 =\int_{\mathbb R}|\xi|\,|\widehat{{\rm u}_{\rm c}}(\xi)|^2\,{\rm d}\xi
 <\infty.
\]
Its trace is ${\rm u}_{\rm c}\in L^2(\mathbb R)$, so
$D_{\mathbb H}{\rm u}_{\rm c}\in\mathscr W_{\mathbb H}$ by
\eqref{eq:robin_energy_space}. Thus it remains to prove that
$\nabla\phi\in L^2(\mathbb H)$, since $\phi$ has zero boundary trace.
Moreover, $\phi$ is bounded and belongs to $H^1(\mathbb H\cap B_R)$
for every $R>0$: boundedness follows from the Poisson formula, and
local regularity from
$\phi=\psi-D_{\mathbb H}{\rm u}_{\rm c}-\psi_0$.

A direct computation gives
\begin{equation}\label{eq:259}
 \left\{
 \begin{aligned}
  \Delta\phi(x,y)&=e^{-\lambda y}g(x)&&\left(\begin{bmatrix}x\\ y\end{bmatrix}\in\mathbb H\right),\\
  \phi(x,0)&=0&&(x\in\mathbb R),
 \end{aligned}
 \right.
\end{equation}
where
\[
 \begin{aligned}
 g(x)&=-\bigl({\rm u}_0''(x)+\lambda^2{\rm u}_0(x)\bigr)\\
 &=-d_+\bigl(\chi_+''(x)+2i\lambda\chi_+'(x)\bigr)e^{i\lambda x}
   -d_-\bigl(\chi_-''(x)-2i\lambda\chi_-'(x)\bigr)e^{-i\lambda x}
  \in C_c^\infty(\mathbb R).
 \end{aligned}
\]
Choose $\eta\in C_c^\infty(\mathbb R^2)$ with $0\leqslant\eta\leqslant1$,
$\eta=1$ on $B_1$ and $\eta=0$ outside $B_2$, and set
\[
\eta_R(x,y)=\eta\left(\frac{x}{R},\frac{y}{R}\right).
\]
In dimension two,
$\|\nabla\eta_R\|_{L^2(\mathbb H)}\leqslant C$, independently of $R$.
The function $\eta_R^2\,\phi$ has compact support and zero boundary
trace, and is therefore an admissible $H^1_0$ test function in
\eqref{eq:259}, by density. Taking real parts of the resulting
identity and applying Young's inequality yields
\[
 \int_{\mathbb H}\eta_R^2|\nabla\phi|^2\,{\rm d}x\,{\rm d}y
 \leqslant 2\int_{\mathbb H}e^{-\lambda y}|g(x)|\,|\phi(x,y)|
                       \,{\rm d}x\,{\rm d}y
       +4\int_{\mathbb H}|\nabla\eta_R|^2|\phi|^2\,{\rm d}x\,{\rm d}y.
\]
Consequently,
\[
 \int_{\mathbb H\cap B_R}|\nabla\phi|^2\,{\rm d}x\,{\rm d}y
 \leqslant \frac{2}{\lambda}\|g\|_{L^1(\mathbb R)}
                         \|\phi\|_{L^\infty(\mathbb H)}
           +C\|\phi\|_{L^\infty(\mathbb H)}^2,
\]
where $C$ is independent of $R$. Letting $R\to\infty$ proves
$\nabla\phi\in L^2(\mathbb H)$, hence
$\phi\in\mathscr W_{\mathbb H}$ and
$\psi_{\rm c}\in\mathscr W_{\mathbb H}$.

Finally, the trace on $\mathbb R$ of any remainder in
$\mathscr W_{\mathbb H}$ belongs to $L^2(\mathbb R)$.
Taking traces of two decompositions of the form \eqref{eq-343}
therefore reduces uniqueness of $d_\pm$ to exactly the half-line
$L^2$ argument at the end of the proof of \cref{lem:repRH-boundary}.
Thus the coefficients are unique and are given by \eqref{eq:dpm}.
\end{proof}

\section{The limiting absorption principle for $\Lambda_{\Omega}$}\label{sec:free}

In this section, we establish the limiting absorption principle for the operator $\Lambda_\Omega$, introduced in  \cref{prop:DtN_L2}, and thereby establish \cref{thm:mainlap}. As outlined in our proof strategy, the main idea is to view $\Lambda_\Omega$ as a compact perturbation of the free reference operator $\Lambda_{\mathbb{H}}$, introduced in \cref{sec:free_waves}. More precisely, to prove  \cref{thm:mainlap} we are going to apply \cref{thm:transfer_lap}. 

We continue to use the notation in the previous sections for the sets $\mathbb{H}, \Omega$, $\Gamma$, $\mathcal{E}, \mathcal{B}$
and for the operators $\Lambda_{\mathbb{H}}$, $\mathcal{R}_{\mathbb{H}},$ $D_\Omega$ and $\Lambda_\Omega$. 

Recall from \cref{thm:self-adjoint} that $\Lambda_\Omega$ is a non-negative self-adjoint operator on $L^2(\mathcal E)$ with domain $\mathcal D(\Lambda_\Omega)=H^1(\mathcal E)$. Consequently, its resolvent set contains $\mathbb C\setminus[0,\infty)$. For each
$z\in\mathbb C\setminus[0,\infty)$, the resolvent of
$\Lambda_\Omega$ is defined by
\begin{equation}\label{def-reolvent-main}
     R_{\Lambda_\Omega}(z)
    :=(\Lambda_\Omega-z)^{-1}
    \in\mathcal L(L^2(\mathcal E)).
\end{equation}
To extend $R_{\Lambda_\Omega}$ as an operator on
$L^2(\mathbb R)$, we first write 
\begin{equation*}
       L^2(\mathbb R)
    =L^2(\mathcal E)\oplus L^2(x_-,x_+). 
\end{equation*}
For each $z\in\mathbb C\setminus[0,\infty)$, we then extend
$\ R_{\Lambda_\Omega}(z)$ to $L^2(\mathbb R)$ by defining
\begin{equation}\label{def-resolvent-extension}
    \mathcal R_{\Lambda_\Omega}(z)
    :=  R_{\Lambda_\Omega}(z) \oplus 0
    \in\mathcal L(L^2(\mathbb R)).
\end{equation}
Thus, the extended resolvent coincides with the original
resolvent on $L^2(\mathcal E)$ and vanishes on
$L^2(x_-,x_+)$. More precisely, we have 
\begin{equation}
    \mathcal{R}_{\Lambda_{\Omega}}(z) = \mathfrak{E}R_{\Lambda_\Omega}(z) \mathfrak{E}^* \quad (z \in \mathbb{C} \setminus [0, \infty)),
\end{equation}
where $\mathfrak{E} : L^2(\mathcal{E}) \to L^2(\mathbb{R})$ is the extension by zero operator, i.e., 
\begin{equation} \label{def-zero-ext}
    (\mathfrak{E}f)(x)  = \begin{dcases}
        f(x) & x\in \mathcal{E}, \\ 0 & x \in (x_-, x_+).
    \end{dcases}
\end{equation}
Note that $\mathfrak{E}^*f  = f|_{\mathcal{E}},$ so that 
\begin{equation*}
    \left( \mathcal{R}_{\Lambda_{\Omega}}(z) f \right)(x) = \begin{dcases}
        \left( R_{\Lambda_{\Omega}}(z) f|_{\mathcal{E}} \right)(x) & x \in \mathcal{E}, \\ 0 & x \in (x_-, x_+).
    \end{dcases}
\end{equation*}

For notational simplicity, we henceforth denote the resolvent $R_{\Lambda_{\Omega}}$ by $R_{\Omega}$ and the  extended
resolvent $\mathcal R_{\Lambda_\Omega}(z)$ by $\mathcal R_\Omega(z)$.  

We start with a technical lemma. 

\begin{lemma} \label{lem-ext00}
        Let $\mathscr{W}_{\Omega}$ and $\mathscr{W}_{\mathbb{H}}$ be the spaces introduced in \eqref{eq:robin_energy_space}, with the norms defined in \eqref{eq:robin_energy_form}. 
        Then there exists a bounded operator $L\in  \mathcal{L}(\mathscr W_\Omega,\mathscr W_\mathbb{H})$  such that  for every $g\in \mathscr W_\Omega$ we have 
    \begin{equation} \label{est-ext00}
    Lg|_\Omega = g,\qquad
        \|Lg\|_{H^1(\mathcal{B})} + \|Lg(\cdot, 0)\|_{L^2(x_{-}, x_+)} \leqslant C \|g\|_{\mathscr{W}_{\Omega}}.
    \end{equation}
\end{lemma} 

\begin{proof}
Choose $R>0$ sufficiently large that
$\overline{\mathcal B}\subset B_{R/2}$, and put
$\Omega_R=\Omega\cap B_R$. By Assumption~\ref{ass:geometry},
$\Omega_R$ and $\mathcal B$ are bounded Lipschitz domains, so that we can
apply the Sobolev extension and trace theorems.

Let $g\in\mathscr W_\Omega$. By
\eqref{eq:energy_controls_local_H1}, its restriction to $\Omega_R$ belongs
to $H^1(\Omega_R)$ and satisfies
\[
 \|g\|_{H^1(\Omega_R)}\leqslant C_R\|g\|_{\mathscr W_\Omega}.
\]
The extension theorem for bounded Lipschitz domains provides a bounded
linear operator $\mathcal X_R:H^1(\Omega_R)\to H^1(\mathbb R^2)$ such that
$(\mathcal X_R h)|_{\Omega_R}=h$. Set
$\widetilde g=\mathcal X_R(g|_{\Omega_R})$ and define, almost everywhere
in $\mathbb H$,
\[
 g_{\textnormal e}(x,y)=
 \begin{dcases}
 g(x,y) &\left(\begin{bmatrix}x\\ y\end{bmatrix}\in\Omega\right),\\
 \widetilde g(x,y) &\left(\begin{bmatrix}x\\ y\end{bmatrix}\in\mathcal B\right).
 \end{dcases}
\]
On $\mathbb H\cap B_R$ this function agrees almost everywhere with
$\widetilde g$, because $\widetilde g=g$ on $\Omega_R$.
Outside $\overline{B_{R/2}}$ it coincides with $g$.
These two descriptions coincide on their overlap, so
$g_{\textnormal e}\in H^1_{\mathrm{loc}}(\mathbb H)$; in particular,
there is no jump across $\Gamma$ and its weak gradient is given by
the corresponding gradients on $\Omega$ and $\mathcal B$.
Moreover,
\[
 \|g_{\textnormal e}\|_{H^1(\mathcal B)}
 \leqslant\|\widetilde g\|_{H^1(\mathbb R^2)}
 \leqslant C_R\|g\|_{H^1(\Omega_R)}
 \leqslant C\|g\|_{\mathscr W_\Omega}.
\]
The trace theorem on $\mathcal B$, applied to its horizontal boundary
segment $(x_-,x_+)\times\{0\}$, therefore yields
\[
 \|g_{\textnormal e}(\cdot,0)\|_{L^2(x_-,x_+)}
 \leqslant C\|g_{\textnormal e}\|_{H^1(\mathcal B)}
 \leqslant C\|g\|_{\mathscr W_\Omega}.
\]
All these traces are understood as in
Remark~\ref{rem:homogeneous_boundary_traces}.
On $\mathcal E$, the trace of $g_{\textnormal e}$ equals that of $g$.
Consequently,
\[
 \begin{aligned}
 \|g_{\textnormal e}\|_{\mathscr W_{\mathbb H}}^2
 &=\|\nabla g\|_{L^2(\Omega)}^2
   +\|\nabla\widetilde g\|_{L^2(\mathcal B)}^2
   +\|g|_{\mathcal E}\|_{L^2(\mathcal E)}^2
   +\|g_{\textnormal e}(\cdot,0)\|_{L^2(x_-,x_+)}^2\\
 &\leqslant C\|g\|_{\mathscr W_\Omega}^2.
 \end{aligned}
\]
This proves membership in $\mathscr W_{\mathbb H}$ and
\eqref{est-ext00}. The radius $R$ was fixed in terms of the geometry,
so the constant is independent of $g$. Thus, it suffices to take $Lg=g_{\textnormal e}$.
\end{proof}

The following result is the main ingredient to prove the limiting absorption principle for $\Lambda_{\Omega}.$

\begin{proposition}\label{prop:compact_perturbation}
    Let $\frac12<s<\frac32$. The operator
    \begin{equation}\label{eq:compact_difference_at_minus_one}
        D_1
        := \mathcal{R}_{\Omega}(-1)  - \mathcal{R}_{\mathbb{H}}(-1) \in \mathcal{L}(L^2(\mathbb{R}))
    \end{equation}
 extends to a compact operator from $L^2_{-s}(\mathbb R)$ to $L^2_s(\mathbb R)$.
\end{proposition}

\begin{proof}
We first show that there exists a constant $C > 0$ with 
\begin{equation} \label{est-D-first}
    \|D_1 f\|_{L^2(\mathbb{R})} \leqslant C_s \|f\|_{L^2_{-s}(\mathbb{R})} \qquad\qquad (f \in L^2(\mathbb{R})).
\end{equation}
To this aim, recall the spaces $\mathscr W_G$ and the forms $a_G$ defined in \eqref{eq:robin_energy_space}--\eqref{eq:robin_energy_form}, for $G\in\{\Omega,\mathbb H\}$. We next take $f \in L^2(\mathbb{R})$ and we set 
    \begin{equation}\label{lipsea}
     {\rm v} = \mathcal{R}_{\Omega}(-1)f.
    \end{equation}
Note that ${\rm v}\in L^2(\mathbb{R})$ and its restriction to $\mathcal{E}$ lies in  $H^1(\mathcal{E})$. Moreover, according to \eqref{def-resolvent-extension}, we have that ${\rm v}$ vanishes on $(x_-,x_+)$. Therefore, above and in the remaining part of this proof we identify ${\rm v}$ and its restriction to $\mathcal{E}$. With this convention we can set
\begin{equation}\label{phi_prop_6.2}
 \phi = D_{\Omega} {\rm v},
\end{equation}
where $D_\Omega$ is the Dirichlet map introduced in \cref{prop:h1lifting}. It follows that   $\phi \in \mathscr{W}_\Omega$. Moreover, $\phi$ satisfies 
\begin{equation} \label{eq:var-phi}
    a_{\Omega}(\phi, g) = \int_{\mathcal{E}} f(x) \overline{g}(x,0) \, {\rm d}x  \qquad\qquad (g \in \mathscr{W}_{\Omega}).
\end{equation}

Indeed, from \cref{lem:operatorgreen}, for every $g \in \mathscr{W}_{\Omega}$ we have
\begin{equation}
    \int_\Omega \nabla \left(D_{\Omega} {\rm v}\right) \cdot \nabla \overline{g} \, {\rm d} x\,{\rm d}y = \int_{\mathcal{E}} (\Lambda_\Omega {\rm v})(x)\,\overline{g}(x,0) \, {\rm d}x = \int_{\mathcal{E}} (f(x) - {\rm v}(x))  \overline{g}(x,0) \ {\rm d} x.
\end{equation}
Rearranging the terms in the above identity, we obtain \eqref{eq:var-phi}. 

We next define  
\begin{equation} \label{DHRH-7.2}
    {\rm u} = \mathcal{R}_{\mathbb{H}}(-1)f, \qquad \psi = D_{\mathbb{H}} {\rm u}.
\end{equation}
Then by \cite[Proposition 3.4]{Ocqueteau2026}, we have ${\rm u} \in H^1(\mathbb{R})$, $\psi \in \mathscr{W}_{\mathbb{H}}$, and $\psi$ satisfies
\begin{equation} \label{eq:var-psi}
    a_{\mathbb{H}} (\psi, h) = \int_{\mathbb{R}} f(x) \overline h(x,0)\,{\rm d}x \qquad\qquad (h \in \mathscr{W}_{\mathbb{H}}). 
\end{equation}
Define $\vartheta = \phi  - \psi|_\Omega$, where $\phi$ has been defined in \eqref{phi_prop_6.2}.
Taking an arbitrary $g\in \mathscr W_\Omega$  in \eqref{eq:var-phi}, choosing its extension $h=Lg$ constructed in \cref{lem-ext00}, inserting it in \eqref{eq:var-psi}, and subtracting the two equalities, we obtain 
\begin{equation} \label{eq:518}
    a_{\Omega} (\vartheta, g) = \int_{\mathcal{B}} \nabla \psi \cdot \nabla \bar{h} \ {\rm d}x\, {\rm d}y + \int_{x_-}^{x_+} (\psi(x, 0) -  f(x)) {\bar h}(x, 0) \ {\rm d}x.
\end{equation}
By \cref{lem-local-estimate-free}, the trace theorem, and \eqref{est-ext00}, we have 
\begin{equation*}
      | a_{\Omega}(\vartheta, g)| \leqslant C \|f\|_{L^2_{-s}(\mathbb{R})} \|g\|_{\mathscr{W}_\Omega}.
\end{equation*}
 Taking $g = \vartheta$ yields
\begin{equation} \label{est-406}
    \| \vartheta\|_{\mathscr{W}_\Omega} \leqslant C \|f\|_{L^2_{-s}(\mathbb{R})}.
\end{equation}

Evaluating \eqref{eq:518} with test functions $g \in C_c^\infty(\Omega)$ yields $\Delta \vartheta = 0$ in the sense of distributions. Subsequently, applying the generalized Green's formula from \cref{lem:brown_traces_green} for arbitrary $g \in \mathscr{W}_\Omega$ recovers the boundary conditions, showing that $\vartheta$ satisfies the boundary value problem
\begin{equation} \label{def-D3}
\left\{
\begin{aligned}
        &\Delta \vartheta(x,y) = 0 && \left(\begin{bmatrix}x \\ y \end{bmatrix} \in \Omega\right), \\
        &\frac{\partial \vartheta}{\partial n}(x,y) = - \frac{\partial \psi}{\partial n}(x,y) && \left(\begin{bmatrix}x \\ y \end{bmatrix} \in \Gamma\right), \\
        &-\frac{\partial \vartheta}{\partial y}(x,0) + \vartheta(x,0) = 0 &&(x\in \mathcal{E}).
\end{aligned}
\right.
\end{equation}

To determine  $D_1 f$, we evaluate the resolvent difference pointwise. On the free surface $\mathcal{E}$, the definitions of the operators imply $(D_1 f)(x) = {\rm v}(x) - {\rm u}(x) = \phi(x,0) - \psi(x,0) = \vartheta(x,0)$. On the immersed interval $(x_-, x_+)$, the extended resolvent $\mathcal{R}_\Omega(-1)$ vanishes by the definition of the zero extension in \eqref{def-zero-ext}, leaving $(D_1 f)(x) = 0 - {\rm u}(x) = -\psi(x,0)$. Thus,
\begin{equation} \label{def-D2}
    (D_1 f)(x) = \begin{dcases}
        \vartheta(x, 0)  & (x \in \mathcal{E}),  \\
        - \psi(x, 0) & (x \in (x_-, x_+)).
    \end{dcases}
\end{equation}
The estimate \eqref{est-406} above and \cref{lem-local-estimate-free} yields \eqref{est-D-first}.

We next show that $D_1$ extends to a bounded linear operator from $L^2_{-s}(\mathbb{R})$ to $L^2_s(\mathbb{R}).$ The idea is to prove decay estimates for $(D_1f)(x)$ for $|x|$ large. To this aim,   we choose a cutoff function  $\chi\in C_c^\infty(\mathbb R^2)$ which is one on a neighbourhood of $\overline{\mathcal B}$, and fix $R_0>1$ so that $\mathrm{supp} \chi\subset B_{R_0}$. We define 
    \begin{equation*}
       \eta(x,y)  = 
       \begin{dcases}
           (1- \chi(x,y))\vartheta(x,y) & \left(\begin{bmatrix}x\\ y\end{bmatrix} \in \Omega\right) \\
           0 & \left(\begin{bmatrix}x\\ y\end{bmatrix} \in \mathcal{B}\right).
       \end{dcases}
    \end{equation*}
    Then $\eta \in \mathscr{W}_{\mathbb{H}}.$ Moreover, $\eta$ satisfies 
     \begin{equation} 
    \left\{
    \begin{aligned}
     &-\Delta \eta(x,y) = h(x,y) &&\left(\begin{bmatrix}x \\ y \end{bmatrix} \in \mathbb{H}\right)\\
     &-\partial_y \eta(x,0) + \eta(x,0) = 0 &&(x\in \mathbb{R}),  
    \end{aligned}
    \right.
    \end{equation}
    where $h = 2\nabla\chi \cdot \nabla \vartheta + (\Delta \chi)\vartheta.$  
 From \eqref{est-406}, we have 
    \begin{equation} \label{est-954}
        \|h\|_{L^2(\mathbb{H})} \leqslant C \|f\|_{L^2_{-s}(\mathbb{R})}. 
    \end{equation}
    Moreover, $\mathrm{supp}\, h \subset \mathbb{H} \cap B_{R_0}.$ By the  Fourier transform, $\eta$ admits the following integral representation 
    
\begin{equation}
    \eta (x, 0) =  \int_{\mathbb{H}} G(x-r, y) h (r, y) \ {\rm d}r\,{\rm d}y, 
\end{equation}
where the kernel $G$ is defined in \eqref{def-G}. Let  $M > 2R_0.$ Note that,
\begin{equation*}
    |x-r|\geqslant |x| - |r| \geqslant |x|  - R_0 \geqslant \frac{|x|}{2} \quad \left(|x| > M, (r, y) \in \mathrm{supp}\;h \right).
\end{equation*}
Then, \eqref{decay-G} and  \eqref{est-954} shows 
\begin{equation} 
    |\eta(x,0)| \leqslant C_{s, R_0} |x|^{-2} \|f\|_{L^2_{-s}(\mathbb{R})} \leqslant C_{s, R_0} \langle x \rangle^{-2}\|f\|_{L^2_{-s}(\mathbb{R})} \qquad (|x| \geqslant M),
\end{equation}
where $C_{s,R_0}$ is independent of $M.$ Therefore
\begin{equation} \label{eq:311}
    \int_{|x|\geqslant M} \langle x\rangle^{2s} |(D_1 f)(x)|^2 \ {\rm d}x \leqslant C_{s, R_0} \|f\|^2_{L^2_{-s}(\mathbb{R})} \int_{|x|\geqslant M} \langle x\rangle^{2s-4} {\rm d}x \leqslant C_{s, R_0} M^{2s-3} \|f\|^2_{L^2_{-s}(\mathbb{R})}.
\end{equation}
On the other hand, from the expression of $D_1f$ from \eqref{def-D2}, \cref{lem-local-estimate-free} and \eqref{est-406} it follows that, 
\begin{equation} \label{eq-255}
      \int_{|x|\leqslant M} \langle x\rangle^{2s} |(D_1 f)(x)|^2 \ {\rm d}x \leqslant \langle M \rangle^{2s} \left( \|\vartheta(\cdot, 0)\|^2_{L^2(\mathcal{E})} + \|\psi(\cdot, 0)\|^2_{L^2(x_-, x_+)}\right) \leqslant C_{s, M} \|f\|_{L^2_{-s}(\mathbb{R})}^2.
\end{equation}
Combining, the above two estimates we obtain 
\begin{equation*}
    \|D_1f \|_{L^2_s(\mathbb{R})} \leqslant C \|f\|_{L^2_{-s}(\mathbb{R})},
\end{equation*}
so that $D_1$ extends as a bounded linear operator from $L_{-s}^2(\mathbb{R})$ to $L_s^2(\mathbb{R})$.

Finally, we prove the fact that this extension is a compact operator from $L_{-s}^2(\mathbb{R})$ into $L_s^2(\mathbb{R})$. 

Let $ M > 2R_0,$ and we set 
\begin{equation*}
    ({D}_{1}^{M} f)(x) = \begin{dcases}
        (D_1 f)(x) & \mbox { if } |x| \leqslant M, \\
        0 & \mbox{ if } |x| > M.
    \end{dcases}
\end{equation*}
Then \cref{lem-local-estimate-free} and \eqref{est-406} show that 
\begin{equation*}
    \|{D}_1^M f\|_{X_M} \leqslant  C_M\left(
\|\vartheta\|_{H^1(\Omega\cap B_{2M})}
+\|\psi\|_{H^1(\mathbb H\cap B_{2M})} 
\right) \leqslant C_{M,s}\|f\|_{L^2_{-s}(\mathbb{R})},
\end{equation*}
where
\begin{equation*}
   X_M=
H^{\frac{1}{2}}(-M,x_-)\oplus
H^{\frac{1}{2}}(x_-,x_+)\oplus
H^{\frac{1}{2}}(x_+,M). 
\end{equation*}
Since  $X_M$ embeds compactly into $L^2(-M,M)$, and $L^2$  and $L^2_s$ norms are equivalent on this interval, we obtain  ${D}_1^M:L^2_{-s}(\mathbb R)\to L^2_s(\mathbb R)$ is compact. Finally, the estimate \eqref{eq:311} shows that 
\begin{equation*}
    \|D_1 -{D}_1^M\|_{\mathcal{L}(L^2_{-s}(\mathbb{R}), L^2_s(\mathbb{R}))} \leqslant C_s M^{2s -3} \xrightarrow[M\to \infty]{} 0,
 \end{equation*}
 since $s< \frac{3}{2}$. Thus $D_1$  extends to a compact operator from $L^2_{-s}(\mathbb R)$ to $L^2_s(\mathbb R)$. This completes the proof of the proposition.  
\end{proof}

\begin{remark}
    The restriction $s<\frac32$ cannot be removed from the statement of \cref{prop:compact_perturbation}. 
Indeed, take a nonzero, nonnegative $f\in C_c^\infty(\mathbb{R})$ such that $\mathrm{supp}\,f \subset (x_-, x_+).$ For such $f$, from \eqref{eq-formula-u}
it is easy to see that 
\begin{equation*}
    D_1(f) = - \int_{x_{-}}^{x_+} G(x- r, 0) f(r) {\rm d} r.
\end{equation*}
Thus
\begin{equation*}
x^2 D_1f(x)
=-\frac1\pi\int_{x_-}^{x_+} \int_0^\infty
z e^{-z}\,\frac{ x^2}{(x-r)^2+z^2}f(r)\,{\rm d}z\,{\rm d}r.    
\end{equation*}
By dominated convergence theorem, 
\begin{equation*}
\lim_{|x| \to \infty} x^2 D_1f(x) = - \frac{1}{\pi} 
\left(\int_0^\infty ze^{-z}\,{\rm d}z\right)
\left(\int_{x_-}^{x_+}f(r)\,{\rm d}r\right)
=-\frac1\pi \int_{x_-}^{x_+}f(r)\,{\rm d}r \neq 0.
\end{equation*}
In particular for sufficiently large $|x|,$ 
\begin{equation*}
    | D_1 f(x)| \geqslant \frac{1}{2\pi x^2} \int_{x_-}^{x_+}f(r)\,{\rm d}r.
\end{equation*}
Hence, 
\begin{equation*}
    \|{D}_1f\|_{L^2_s(\mathbb{R})}^2 \geqslant C \int_{R}^\infty x^{2s-4}\;{\rm d}x,
\end{equation*}
for some large $R > 0.$ Therefore, if $s \geqslant \frac{3}{2}$ then $D_1$ does not map $L^2_{-s}(\mathbb{R})$ into $L^2_s(\mathbb{R})$. 
\end{remark}

\begin{remark} \label{rem-local-estimate-R-Omega}
    Let $1 \leqslant s < \frac{3}{2},$ and let $R > 0$ such that $[x_-, x_+] \subset (-R, R).$ Then, by \eqref{lipsea},\eqref{phi_prop_6.2}, \eqref{DHRH-7.2}, \eqref{est-406} and \cref{lem-local-estimate-free}, there exists a constant $C_{R,s} > 0$ such that 
    \begin{equation}
        \norm{D_{\Omega} R_{\Omega}(-1)\mathfrak E^*f}_{H^1(\Omega\cap B_R)} \leqslant C_{R,s} \norm{f}_{L^2_{-s}(\mathbb{R})} \qquad (f \in L^2(\mathbb{R})).
    \end{equation}
\end{remark}

To prove \cref{thm:mainlap} we are going to apply \cref{thm:transfer_lap}.
To this aim, We introduce the several operators. 

Assume  $1 \leqslant s < \frac{3}{2}$ and recall the definition of $\mathcal{R}_{\mathbb{H}}^{\pm}(\lambda)$ from \cref{prop:free-lap}. Let $\iota:L^2_s(\mathbb{R})\hookrightarrow L^2_{-s}(\mathbb{R})$ denote  the natural inclusion map. For $\lambda > 0,$ define 
\begin{equation} \label{eq:def_K1}
    K_1^{\pm}(\lambda) = D_1 \left( \iota + (\lambda +1) \mathcal{R}_{\mathbb{H}}^{\pm}(\lambda) \right),
    \end{equation}
\begin{equation}\label{eq:def_F1}
    F_1^\pm(\lambda) = \mathbb{I} -(\lambda +1) K_1^{\pm}(\lambda), 
\end{equation}
where $\mathbb{I}$ denotes the identity operator on $L^2_s(\mathbb{R}).$ It follows from \cref{prop:free-lap} and \cref{prop:compact_perturbation} that 
\begin{equation*}
    K_1^{\pm}(\lambda), F_1^{\pm}(\lambda) \in \mathcal{L}(L^2_s(\mathbb{R}))
\end{equation*}
for any $\lambda > 0$ and $1 \leqslant  s< \frac{3}{2}.$ Moreover, according to \cref{lem:no_fixed_eigenvalues} and \cref{thm:transfer_lap}, $F_1^{\pm}(\lambda)$ is invertible. This allows us to define, for $\lambda > 0$
\begin{equation} \label{eq:ROmega-boundary-ext}
    \mathcal{R}_{\Omega}^{\pm}(\lambda) := \left[ \mathcal{R}_{\mathbb{H}}^{\pm}(\lambda) + \iota K_1^{\pm}(\lambda) \right] \left[ F_1^{\pm} (\lambda)\right]^{-1} \in \mathcal{L}(L^2_s(\mathbb{R}), L^2_{-s}(\mathbb{R})) \quad \left( 1 \leqslant s < \frac{3}{2} \right).
\end{equation}
Finally, for each $\lambda > 0,$ let $R_{\Omega}^{\pm}(\lambda) \in \mathcal{L}(L^2_s(\mathcal{E}), L^2_{-s}(\mathcal{E}))$ be defined by 
\begin{equation} \label{eq:ROmega-boundary}
    R_{\Omega}^{\pm}(\lambda) = \mathfrak{E}^{*} \mathcal{R}_{\Omega}^{\pm}(\lambda) \mathfrak{E},
\end{equation}
where $\mathfrak{E}$ is defined in \eqref{def-zero-ext}.
 
We are now in a position to prove \cref{thm:mainlap}.
\begin{proof}[Proof of \cref{thm:mainlap}]
To establish the limiting absorption principle for $\Lambda_{\Omega}$, we apply the abstract transfer result of \cref{thm:transfer_lap} to the extended operators acting on the full line. Let $1\leqslant s < \frac{3}{2}$. 

We verify the required hypotheses. The reference free operator $\Lambda_{\mathbb{H}}$ satisfies the stable LAP by \cref{prop:stable_lap_free}. The compatible weighted scales are provided by the zero-extension operator $\mathfrak{E}$ of \eqref{def-zero-ext} and its adjoint. The existence of the limits for the free resolvent is guaranteed by \cref{prop:free-lap}. Furthermore, \cref{prop:compact_perturbation} establishes that $D_1$ is a compact operator from $L^2_{-s}(\mathbb{R})$ to $L^2_s(\mathbb{R})$, which verifies the compact perturbation conditions of \cref{def:compact_perturbation} (as noted in \cref{remark:sufficient_resolvent_compactness}). Finally, $\Lambda_{\Omega}$ has no positive eigenvalues according to \cref{lem:no_fixed_eigenvalues}.

Since all conditions of \cref{thm:transfer_lap} are met, for any $\lambda > 0$, the limits
\begin{equation*}
    \lim_{\varepsilon \to 0^+} \mathcal{R}_{\Omega}(\lambda \pm i \varepsilon) = \mathcal{R}_{\Omega}^{\pm}(\lambda)
\end{equation*}
exist in $\mathcal{L}(L_s^2(\mathbb{R}),L_{-s}^2(\mathbb{R}))$. By applying the restriction formula \eqref{eq:ROmega-boundary}, we deduce that the limits
\begin{equation*}
    \lim_{\varepsilon \to 0^+} R_{\Omega}(\lambda \pm i \varepsilon) = R_{\Omega}^{\pm}(\lambda)
\end{equation*}
exist in $\mathcal{L}(L_s^2(\mathcal{E}),L_{-s}^2(\mathcal{E}))$, locally uniformly in $\lambda$. This concludes the proof.
\end{proof}

We now establish a representation of the harmonic potential associated
with $R_{\Omega}^{\pm}(\lambda)f$, followed by estimates used in the
proof of quantitative local energy decay.

We use the density result of \cref{lem:dense} to extend the variational
identity from smooth test functions to the space $\mathscr W_\Omega$
defined in \eqref{eq:robin_energy_space}.

\begin{lemma} \label{lem:ROmega-bdy-expression}
Under \cref{ass:geometry}, let $1\leqslant s<\frac32$, $\lambda>0$, and $f\in L^2_s(\mathcal{E})$.
Choose $r_0$ large enough that $[x_-,x_+]\subset(-r_0,r_0)$, and fix the cutoffs of \cref{lem:repRH-boundary}. Let us set 
\begin{equation}
    {\rm v} = R_\Omega^+(\lambda)f \in L^2_{-s}(\mathcal{E}), \quad h  = [F_1^+(\lambda)]^{-1} \mathfrak{E}f \in L^2_s(\mathbb{R}).
\end{equation}
Then there exists $\phi \in H^1(\Omega \cap B_R)$ for every $R > 0$, satisfying 
\begin{equation} \label{eq:ROmega-bdy-equation-1}
\left\{
\begin{aligned}
        &\Delta \phi(x,y)  = 0  &&\left(\begin{bmatrix}x\\ y\end{bmatrix}\in  \Omega\right), \\
        &\frac{\partial \phi}{\partial n}(x,y) = 0  &&\left(\begin{bmatrix}x\\ y\end{bmatrix}\in\Gamma\right), \\ 
        &-\partial_y\phi(x,0)-\lambda\phi(x,0)=f(x) &&(x\in \mathcal E),
        \\
        &\phi(x,0)  
        = {\rm v}(x) &&(x\in \mathcal E).
\end{aligned}
\right.
\end{equation}
Here the Dirichlet trace is the local Sobolev trace, and the Neumann and Robin conditions are understood weakly: for every $\eta\in C_c^\infty(\overline{\Omega})$,
\[
\int_\Omega\nabla\phi\cdot\overline{\nabla\eta}\,{\rm d}x\,{\rm d}y
-\lambda\int_{\mathcal E}{\rm v}\,\overline\eta\,{\rm d}x
=\int_{\mathcal E}f\,\overline\eta\,{\rm d}x.
\]
Furthermore, $\phi$ admits a decomposition into an explicit outgoing wave component $\phi_0$ and a localized near-field correction $\phi_{\rm c}$:
\begin{equation} \label{eq:1035}
    \phi(x,y) = \phi_0(x,y) + \phi_{\rm c}(x,y).
\end{equation}
Here, the outgoing component is given by
\begin{equation}\label{def-gg1}
    \phi_0(x,y) = e^{-\lambda y}\left(d_+(\lambda)\chi_+(x)e^{i\lambda x} + d_-(\lambda)\chi_-(x)e^{-i\lambda x}\right),
\end{equation}
with uniquely determined amplitudes $d_{\pm}(\lambda) = i\sqrt{2\pi}\,\widehat{h}(\pm\lambda)$. The remainder $\phi_{\rm c} \in \mathscr{W}_{\Omega}$ satisfies 
\begin{equation}\label{eq:phic}
\begin{aligned}
    \int_{\Omega} \nabla \phi_{\rm c} \cdot \overline{\nabla \eta} \, {\rm d}x\,{\rm d}y  - \lambda \int_{\mathcal{E}} \phi_{\rm c}(x, 0) {\bar \eta}(x, 0) \, {\rm d}x  \\
    =  \int_{\mathcal{E}} f(x) {\bar \eta}(x, 0) \, {\rm d}x  
    + \int_{\Omega} e^{-\lambda y} g_1(x) {\bar \eta}(x, y) \, {\rm d}x\,{\rm d}y \qquad ( \eta \in \mathscr{W}_{\Omega}),
\end{aligned}
\end{equation}
where $g(x)=d_+\chi_+(x)e^{i\lambda x}+d_-\chi_-(x)e^{-i\lambda x}$ and $g_1(x) = (d_+\chi_+e^{i\lambda x}+d_-\chi_-e^{-i\lambda x})'' + \lambda^2(d_+\chi_+e^{i\lambda x}+d_-\chi_-e^{-i\lambda x})$.
\end{lemma}

\begin{proof}
Let us set
\begin{equation}
    {\rm u} = \mathcal{R}_{\mathbb{H}}^+(\lambda)h, \quad \tilde{h} = h + (\lambda +1){\rm u}.
\end{equation} 
Then  ${\rm u}, \tilde{h} \in L^2_{-s}(\mathbb{R}).$ From the definition of $R_{\Omega}^+(\lambda)$ in \eqref{eq:ROmega-boundary}, it follows that 
\begin{equation} \label{eq:1049}
    {\rm v} = \mathfrak{E}^* {\rm u} + \mathfrak{E}^*D_1 \widetilde{h} = {\rm u}|_{\mathcal{E}} + D_1 \widetilde{h}|_{\mathcal{E}} \mbox{ and } h - (\lambda + 1) D_1 \tilde{h} = \mathfrak{E}f.
\end{equation}
We define 
\begin{equation} \label{eq:1048}
  \phi_{\mathbb{H}} := D_{\mathbb{H}}{\rm u}, \quad -\partial_y \phi_{\mathbb{H}}  - \lambda \phi_{\mathbb{H}} = h \mbox{ on } \mathbb{R},
\end{equation}
where $D_\mathbb{H}$ is the Dirichlet map defined in \eqref{def_dir_h}. From \cref{lem:DH-lift-boundary}, we have 
\begin{equation}
     \phi_{\mathbb{H}} =  D_{\mathbb{H}} \mathcal{R}_{\mathbb{H}}(-1) \tilde{h},   
\end{equation}
and
\begin{equation*}
  \phi_{\mathbb{H}} = d_+\chi_+(x)e^{i\lambda x -\lambda y}+d_-\chi_-(x)e^{-i\lambda x-\lambda y} +  \phi_{\rm{c}, \mathbb{H}} \quad \left(\begin{bmatrix} x\\ y\end{bmatrix} \in \mathbb{H}\right),
\end{equation*}
where $\phi_{\rm{c}, \mathbb{H}} \in \mathscr{W}_{\mathbb{H}}$.

Next, proceeding as in the proof of  \cref{prop:compact_perturbation} (see, for instance, \eqref{def-D3}--\eqref{def-D2}), there exists $\vartheta \in \mathscr{W}_{\Omega},$ satisfying 
\begin{equation} \label{eq:1039}
\left\{
\begin{aligned}
    &\Delta \vartheta(x,y)  = 0 && \left(\begin{bmatrix} x \\ y \end{bmatrix}\in \Omega\right), \\
    &\frac{\partial \vartheta}{\partial n}(x,y) = - \frac{\partial \phi_{\mathbb{H}}}{\partial n}(x,y)  &&\left(\begin{bmatrix} x \\ y \end{bmatrix}\in \Gamma\right), \\
    &-\partial_y\vartheta(x,0) + \vartheta(x,0) = 0 && (x \in \mathcal E), \\
    &\vartheta(x,0)  = (D_1 \widetilde{h})(x) &&(x\in  \mathcal{E}).
\end{aligned}
\right.
\end{equation}

For any $\begin{bmatrix}x \\ y\end{bmatrix} \in \Omega,$ we set 
\begin{equation}
    \phi = \phi_{\mathbb{H}}|_{\Omega} + \vartheta.
\end{equation}
By construction, $\Delta \phi = 0$ in $\Omega.$ From \eqref{eq:1049}, \eqref{eq:1039}, and \eqref{eq:1048}, we see that $\frac{\partial \phi}{\partial n} =0$ on $\Gamma$, $\phi|_{\mathcal{E}} = \phi_{\mathbb{H}}|_{\mathcal{E}} + \vartheta|_{\mathcal{E}} = {\rm u}|_{\mathcal{E}} + (D_1 \widetilde{h})|_{\mathcal{E}} = {\rm v}$, and for $x \in \mathcal{E}$:
\begin{equation*}
    -\partial_y \phi(x,0) - \lambda\phi(x,0) = h(x) - (1+\lambda) \vartheta(x, 0) =  h(x) - (1+\lambda)(D_1 \tilde h)(x) = f(x). 
\end{equation*}

This $\phi$ satisfies the decomposition \eqref{eq:1035} with $\phi_{\rm c} = \phi_{\rm{c}, \mathbb{H}}|_{\Omega} + \vartheta \in \mathscr{W}_{\Omega}.$ To prove the variational identity \eqref{eq:phic}, we note that $\phi_0$ defined in \eqref{def-gg1} satisfies
\begin{equation} \label{eq:phi0}
\left\{
\begin{aligned}
    &\Delta \phi_0(x,y) = e^{-\lambda y}g_1(x) &&\left(\begin{bmatrix} x \\ y \end{bmatrix}\in \Omega\right),  \\
    &\frac{\partial \phi_0}{\partial n}(x,y) = 0 && \left(\begin{bmatrix} x \\ y \end{bmatrix}\in \Gamma\right), \\ 
    &-\partial_y \phi_0(x,0) - \lambda \phi_0(x,0) = 0 &&(x\in \mathcal{E}).
  \end{aligned}  
\right.
\end{equation}  
Since $\phi = \phi_0 + \phi_{\rm c}$, we take $\eta \in C_c^\infty(\overline{\Omega})$. Using the weak boundary conditions in \eqref{eq:ROmega-bdy-equation-1} and integrating the smooth function $\phi_0$ by parts, we obtain
\begin{equation*}
\int_\Omega \nabla\phi\cdot\overline{\nabla \eta} \, {\rm d}x\,{\rm d}y
-\lambda\int_{\mathcal E} {\rm v}(x)\overline{\eta}(x, 0) \, {\rm d}x
=\int_{\mathcal E} f(x) \overline{\eta}(x,0) \, {\rm d}x, 
\end{equation*}
and 
\begin{equation*}
    \int_\Omega\nabla \phi_0 \cdot\overline{\nabla \eta} \, {\rm d}x \,{\rm d}y
-\lambda\int_{\mathcal E} \phi_0(x, 0) \overline{\eta}(x, 0) \, {\rm d}x
=-\int_\Omega e^{-\lambda y}g_1(x) \overline{\eta}(x,y) \, {\rm d}x\,{\rm d}y. 
\end{equation*}
Subtracting these two identities proves that \eqref{eq:phic} holds for all $\eta \in C_c^\infty(\overline{\Omega})$. Note that, 
\begin{equation*}
    \left| \int_{\Omega} \nabla \phi_{\rm c} \cdot \overline{\nabla \eta} \, {\rm d}x\,{\rm d}y  - \lambda \int_{\mathcal{E}} \phi_{\rm c}(x, 0) {\bar \eta}(x, 0) \, {\rm d}x  
    - \int_{\mathcal{E}} f(x) {\bar \eta}(x, 0) \, {\rm d}x \right| \leqslant C \left( \|\phi_{\rm c} \|_{\mathscr{W}_{\Omega}} + \|f \|_{L^2_s(\mathcal{E})} \right) \|\eta\|_{\mathscr{W}_\Omega}. 
\end{equation*}
For $x \in \mathcal{E},$ we have 
\begin{equation*}
    |\eta(x,y)| \leqslant |\eta(x, 0)| + \int_0^y|\partial_y \eta(x, r)| \,{\rm d}r \leqslant  |\eta(x, 0)| + \sqrt{y} \|\partial_y \eta(x, \cdot)\|_{L^2(0, \infty)}.
\end{equation*}
Consequently, as $g_1$ is supported in $\mathcal{E},$ we have 
\begin{multline} \label{est-753}
   \left| \int_{\Omega} e^{-\lambda y} g_1(x) {\bar \eta}(x, y) \, {\rm d}x\,{\rm d}y \right| \\
   \leqslant \int_{\mathcal{E}} |g_1(x)| \int_0^\infty e^{-\lambda y}|\eta(x, y)| \, {\rm d}y \,{\rm d}x 
   \leqslant \int_{\mathcal{E}} |g_1(x)| \left( \frac{1}{\lambda} |\eta(x, 0)| + \frac{\sqrt{\pi}}{2\lambda^{3/2}} \|\partial_y \eta(x, \cdot)\|_{L^2(0, \infty)}  \right) \, {\rm d}x\\
   \leqslant \|g_1\|_{L^2(\mathcal{E})} \left( \lambda^{-1} \|\eta(\cdot, 0)\|_{L^2(\mathcal{E})} + \frac{\sqrt{\pi}}{2}\lambda^{-3/2}\|\nabla \eta\|_{L^2(\Omega)} \right).
\end{multline}
From the above estimates and using the fact that $C_c^\infty(\overline{\Omega})$ is dense in $\mathscr{W}_{\Omega}$ (see \cref{lem:dense}), we see that \eqref{eq:phic} holds for all $\eta \in \mathscr{W}_{\Omega}$. Finally, taking traces on $\mathcal{E}$, the uniqueness of $d_\pm(\lambda)$ follows from the same half-line $L^2$ argument as in the proof of \cref{lem:repRH-boundary}.
This completes the proof of the lemma. 
\end{proof}

We now use the decomposition of \cref{lem:ROmega-bdy-expression} to estimate the weighted norm of resolvent trace at every positive frequency.

\begin{proposition} \label{prop-lambda4}
      Under the assumptions and with the notation in  \cref{lem:ROmega-bdy-expression}, 
      let ${\rm v} = R_{\Omega}^{+}(\lambda)f.$ Then there exists $C>0$, depending only on $s$ and on $\Omega$ such that
      \begin{equation} \label{est-ROmega-bdy}
          \| {\rm v} \|_{L^2_{-s}(\mathcal{E})} \leqslant C (1+\lambda^{-4}) \|f\|_{L^2_s(\mathcal{E})}
          \qquad\qquad(\lambda>0,\ f\in L_s^2(\mathcal E)).
      \end{equation}
\end{proposition}
\begin{proof}
    We use the notation of \cref{lem:ROmega-bdy-expression}. Throughout the proof, the constant $C$ depends only on $s, \Omega$ and $\chi_{\pm},$ and is independent of $\lambda, f$ and $\phi.$ The constant may change from line to line.   

    We are going to prove the following two estimates:
    \begin{equation} \label{est-v-1}
        \|{\rm v} \|_{L^2_{-s}(\mathcal{E})} \leqslant C \left( (1+ \lambda^{-1}) \| f \|_{L^2_s(\mathcal{E})} + (1+ \lambda^{-2}) \sqrt{|d_+|^2+|d_-|^2}  \right),
    \end{equation}
    and 
    \begin{equation} \label{est-d}
        |d_+|^2+|d_-|^2 \leqslant 2  \|{\rm v} \|_{L^2_{-s}(\mathcal{E})} \| f \|_{L^2_s(\mathcal{E})}.
    \end{equation}
    We first show that \eqref{est-ROmega-bdy} follows from the above two estimates. To this aim, combining  \eqref{est-v-1}, \eqref{est-d} and Young's inequality, we obtain 
    \begin{multline*}
          \|{\rm v} \|_{L^2_{-s}(\mathcal{E})} \leqslant C \left( (1+ \lambda^{-1}) \| f \|_{L^2_s(\mathcal{E})} + (1+ \lambda^{-2}) \sqrt{2  \|{\rm v} \|_{L^2_{-s}(\mathcal{E})} \| f \|_{L^2_s(\mathcal{E})} }\right) \\
          \leqslant C (1+ \lambda^{-1}) \| f \|_{L^2_s(\mathcal{E})} + \frac{1}{2}  \|{\rm v} \|_{L^2_{-s}(\mathcal{E})} + C (1+ \lambda^{-2})^2  \| f \|_{L^2_s(\mathcal{E})}. 
    \end{multline*}
    Therefore,
    \begin{equation*}
        \|{\rm v} \|_{L^2_{-s}(\mathcal{E})} \leqslant C \Big( (1+ \lambda^{-1}) + (1+ \lambda^{-2})^2 \Big)\| f \|_{L^2_s(\mathcal{E})} \leqslant C (1+\lambda^{-4}) \| f \|_{L^2_s(\mathcal{E})}.
    \end{equation*}
    We next prove \eqref{est-v-1}. Recall that $\phi$ satisfies \eqref{eq:ROmega-bdy-equation-1} --\eqref{eq:1035}. Similar to \cref{lem:no_fixed_eigenvalues}, we define, for $x \in \mathcal{E}$
\begin{equation} \label{def-F0-new}
    F(x):=\frac1\lambda\left(
\phi_{\rm c}(x, 0)+\int_0^\infty e^{-\lambda y}\partial_y\phi_{\rm c}(x,y)\,{\rm d}y
\right).
\end{equation}
 Proceeding as in the proof of \cref{lem:no_fixed_eigenvalues} and testing \eqref{eq:phic} with $\eta(x,y) = \vartheta_{\rm e}(x)e^{-\lambda y}$ for $\vartheta \in C_c^\infty(\mathcal{E})$, we obtain $F \in H^1(\mathcal{E})$ and 
  \begin{equation} \label{def-F-1135}
     F(x)=\int_0^\infty e^{-\lambda y}\phi_{\rm c}(x,y)\,{\rm d}y \quad (x \in \mathcal{E}), \quad   F'' + \lambda^2 F = -\left( f + \frac{g_1}{2\lambda} \right),
  \end{equation}
  in the sense of distributions. Since $F \in H^1(\mathcal{E}),$ we obtain 
  \begin{equation}
      F(x)  = \begin{dcases}
             -\frac{1}{\lambda} \int_x^\infty \sin (\lambda(r -x))  \left( f(r) + \frac{g_1(r)}{2\lambda} \right) \; {\rm d}r & x \in (x_+, \infty) \\
              -\frac{1}{\lambda} \int_{-\infty}^x \sin (\lambda(x - r))  \left( f(r) + \frac{g_1(r)}{2\lambda} \right) \; {\rm d}r & x \in (- \infty, x_{-}).
      \end{dcases}
  \end{equation}

To justify this representation, we first note that the  integrals are absolutely convergent, since $ f+\frac{g_1}{2\lambda}\in L^2_1(\mathcal E)\hookrightarrow L^1(\mathcal E).$  Applying \eqref{eq:hardy-tail} to both the sine kernel and its derivative shows that these integrals define an $H^1$ solution of the same one-dimensional equation as $F$. Their difference therefore satisfies the corresponding homogeneous equation $F''+\lambda^2 F =0,$ which admits only the trivial solution in $L^2$ on each half-line. This establishes the desired representation.  

We now compute $L^2$ norm of $F$. Applying \eqref{eq:hardy-tail} on the two half-lines, with endpoints $x_+$ and $x_-$, gives 
  \begin{equation} \label{eq:est-402}
      \|F\|_{L^2(\mathcal{E})} \leqslant C \lambda^{-1} \Big\| \left( f + \frac{g_1}{2\lambda} \right) \Big\|_{L^2_1(\mathcal{E})} \leqslant C  \left( \lambda^{-1} \|f\|_{L^2_s(\mathcal{E})} + \lambda^{-2} \|g_1\|_{L^2_1(\mathcal{E})} \right).
  \end{equation}
  From the definition of $g$ and $g_1$ in \eqref{def-gg1}, we get 
  \begin{equation}
      g_1(x) =  d_+(\chi_+''+2i\lambda\chi_+')e^{i\lambda x} +d_-(\chi_-''-2i\lambda\chi_-')e^{-i\lambda x}.
  \end{equation}
  Thus $g_1$ is compactly supported, and 
  \begin{equation} \label{est-g1}
   \|g\|_{L^2_{-s}(\mathcal{E})} \leqslant C  \sqrt{|d_+|^2+|d_-|^2}, \quad   \|g_1\|_{L^2_1(\mathcal{E})} \leqslant C (1+\lambda) \sqrt{|d_+|^2+|d_-|^2}.
  \end{equation}
  Consequently, from \eqref{eq:est-402}, we deduce that 
  \begin{equation}  \label{est-F-1038}
       \|F\|_{L^2(\mathcal{E})} \leqslant C \lambda^{-1} \left( \|f\|_{L^2_s(\mathcal{E})} + (1+\lambda^{-1})\sqrt{|d_+|^2+|d_-|^2} \right).
  \end{equation}
From \eqref{def-F0-new}, Cauchy--Schwarz inequality  and
$|a+b|^2\leqslant3|a|^2+\frac32|b|^2$ give
\begin{multline*}
    \|\phi_{\rm c}(\cdot, 0)\|_{L^2(\mathcal{E})}^2 \leqslant 3\lambda^2\|F\|_{L^2(\mathcal{E})}^2 + \frac{3}{2} \left\| \int_0^\infty e^{-\lambda y}\partial_y\phi_{\rm c}(\cdot,y)\,{\rm d}y  \right\|^2_{L^2(\mathcal{E})} \leqslant 3\lambda^2\|F\|_{L^2(\mathcal{E})}^2 + \frac{3}{4\lambda} \|\nabla\phi_{\rm c}\|^2_{L^2(\Omega)}.  
\end{multline*}
Thus 
\begin{equation} \label{est-trace-phic}
    \lambda  \|\phi_{\rm c}(\cdot, 0)\|_{L^2(\mathcal{E})}^2 \leqslant 3\lambda^3\|F\|_{L^2(\mathcal{E})}^2 + \frac{3}{4} \|\nabla\phi_{\rm c}\|^2_{L^2(\Omega)}.
\end{equation}
Recall that $\phi_{\rm c}$ satisfies \eqref{eq:phic} and belongs to the space $\mathscr W_\Omega$ defined in \eqref{eq:robin_energy_space}. Thus, by taking $\eta = \phi_{\rm c}$ in \eqref{eq:phic}, we get 
\begin{equation} \label{eq:1121}
    \norm{\nabla \phi_{\rm c}}^2_{L^2(\Omega)} - \lambda \norm{\phi_{\rm c}(\cdot, 0)}^2_{L^2(\mathcal{E})} = \int_{\mathcal{E}} f \; {\bar \phi_{\rm c}}(x, 0) \; {\rm d}x  + \int_{\Omega} e^{-\lambda y} g_1(x) {\bar \phi_{c}}(x, y) \; {\rm d}x{\rm d}y. 
\end{equation}
Taking $\eta  = \phi_{\rm c}$ in \eqref{est-753}, and using \eqref{est-g1} we deduce 
\begin{align*}
    & \left|\int_{\mathcal{E}} f \; {\bar \phi_{\rm c}}(x, 0) \; {\rm d}x  + \int_{\Omega} e^{-\lambda y} g_1(x) {\bar \phi_{c}}(x, y) \; {\rm d}x\,{\rm d}y \right| \\
    & \leqslant \norm{f}_{L^2_s(\mathcal{E})}\norm{\phi_{\rm c}(\cdot, 0)}_{L^2(\mathcal{E})} + C(1+\lambda) \sqrt{|d_+|^2+|d_-|^2} \left( \lambda^{-1} \|\phi_{\rm c}(\cdot, 0)\|_{L^2(\mathcal{E})} + \lambda^{-3/2}\|\nabla \phi_{\rm c}\|_{L^2(\Omega)} \right).
\end{align*}
Consequently, by Young's inequality, for any $\varepsilon > 0$ there exists $C_\varepsilon > 0$ such that
\begin{multline*}
       \norm{\nabla \phi_{\rm c}}^2_{L^2(\Omega)} - \lambda \norm{\phi_{\rm c}(\cdot, 0)}^2_{L^2(\mathcal{E})}  \\
    \leqslant \varepsilon \norm{\nabla \phi_{\rm c}}^2_{L^2(\Omega)} + \varepsilon\lambda \norm{\phi_{\rm c}(\cdot, 0)}^2_{L^2(\mathcal{E})} + C_{\varepsilon} \lambda^{-1}\norm{f}_{L^2_s(\mathcal{E})}^2 + C_{\varepsilon} (1+\lambda)^2 \lambda^{-3}(|d_+|^2+|d_-|^2). 
\end{multline*}
Using $(1+\lambda)^2\lambda^{-3} \leqslant 2(\lambda^{-1} + \lambda^{-3}),$ for fixed small $\varepsilon > 0$ the above inequality can be rewritten as
\begin{multline*}
    (1 - \varepsilon) \norm{\nabla \phi_{\rm c}}^2_{L^2(\Omega)}\leqslant (1 + \varepsilon)\lambda \norm{\phi_{\rm c}(\cdot, 0)}^2_{L^2(\mathcal{E})} + C \lambda^{-1}\norm{f}_{L^2_s(\mathcal{E})}^2 + C (\lambda^{-1} +  \lambda^{-3})(|d_+|^2+|d_-|^2). 
\end{multline*}
The above inequality together with \eqref{est-trace-phic} yields 
\begin{equation*}
    \left(\frac14 - \frac{7\varepsilon}{4} \right) \norm{\nabla \phi_{\rm c}}^2_{L^2(\Omega)} \leqslant 3(1+ \varepsilon)\lambda^3\|F\|_{L^2(\mathcal{E})}^2 +  C \lambda^{-1}\norm{f}_{L^2_s(\mathcal{E})}^2 + C (\lambda^{-1} +  \lambda^{-3})(|d_+|^2+|d_-|^2). 
\end{equation*}
By choosing $\varepsilon$ small, the above inequality can be written as 
\begin{equation*}
    \norm{\nabla \phi_{\rm c}}^2_{L^2(\Omega)} \leqslant C\lambda^3\|F\|_{L^2(\mathcal{E})}^2 +  C \lambda^{-1}\norm{f}_{L^2_s(\mathcal{E})}^2 + C (\lambda^{-1} +  \lambda^{-3})(|d_+|^2+|d_-|^2). 
\end{equation*}
Combining the above estimate together with \eqref{est-F-1038} and \eqref{est-trace-phic}, we deduce 
\begin{align*}
     &  \|\phi_{\rm c}(\cdot, 0)\|_{L^2(\mathcal{E})}^2 \leqslant 3\lambda^2\|F\|_{L^2(\mathcal{E})}^2 + \frac{3}{4} \lambda^{-1} \|\nabla\phi_{\rm c}\|^2_{L^2(\Omega)} \\
     & \leqslant C \lambda^2\|F\|_{L^2(\mathcal{E})}^2 +  C \lambda^{-2} \norm{f}_{L^2_s(\mathcal{E})}^2 +  C (\lambda^{-2} +  \lambda^{-4})(|d_+|^2+|d_-|^2) \\
     & \leqslant C(1 + \lambda^{-2}) \norm{f}_{L^2_s(\mathcal{E})}^2 + C( 1+ \lambda^{-2} + \lambda^{-4} ) (|d_+|^2+|d_-|^2). 
\end{align*}
    Note that $\sqrt{1+\lambda^{-2}} \leqslant 1 + \lambda^{-1}$ and $\sqrt{1+\lambda^{-2} + \lambda^{-4}} \leqslant 1 + \lambda^{-2}.$ Consequently, from the above inequality, we infer that 
    \begin{equation} \label{est-phic-final}
         \|\phi_{\rm c}(\cdot, 0)\|_{L^2_{-s}(\mathcal{E})}  \leqslant \|\phi_{\rm c}(\cdot, 0)\|_{L^2(\mathcal{E})} \leqslant   C \left( (1+ \lambda^{-1}) \| f \|_{L^2_s(\mathcal{E})} + (1+ \lambda^{-2}) \sqrt{|d_+|^2+|d_-|^2}  \right).
    \end{equation}
We note that, from \eqref{eq:ROmega-bdy-equation-1} and \eqref{eq:1035}, ${\rm v}(x) = g(x) + \phi_{\rm c}(x, 0).$ Therefore, from \eqref{est-phic-final} and \eqref{est-g1}, we deduce \eqref{est-v-1}. 
\medskip
\noindent We now give a proof of \eqref{est-d}. Taking imaginary part in \eqref{eq:1121}, we get 
\begin{equation}\label{eq:imaginary_energy}
\mathrm{Im}\left(
\int_{\mathcal E}f\overline {\phi_c}(x, 0) \; {\rm d}x 
+\int_\Omega e^{-\lambda y} g_1(x) \overline{\phi_{\rm c}}(x, y) \; {\rm d}x{\rm d}y
\right)=0.
\end{equation}
Recall that
\begin{equation*}
   F''+\lambda^2F=- \left( f+\frac{g_1}{2\lambda} \right) ,\qquad 
g_1 = g''+\lambda^2g. 
\end{equation*}
Since $F \in H^1(\mathcal{E}),$ $f, g_1 \in L^2(\mathcal{E})$ we have $F \in H^2(\mathcal{E}).$ In particular, $F$ and $F'$ tend to zero at the infinite ends of the
two half-lines. Since $g$ and $g'$ are bounded and vanish near the
finite endpoints, integration by parts on truncated half-lines has
vanishing boundary terms in the limit and gives 
\begin{equation} \label{eq:200}
    - \int_{\mathcal{E}} \left( f(x) + \frac{g_1(x)}{2\lambda} \right) \overline{g}(x) \; {\rm d}x  = \int_{\mathcal{E}} (F''(x) + \lambda^2 F(x)) \overline{g}(x) \; {\rm d}x = \int_{\mathcal{E}} F \overline{g_1}(x) \; {\rm d}x.  
\end{equation}
Using \eqref{def-F-1135} and \eqref{def-gg1} we obtain 
\begin{equation*}
    \int_{\mathcal{E}}  \overline{F}(x) {g_1}(x) \; {\rm d}x = \int_{\mathcal{E}} \int_0^\infty e^{-\lambda y} \overline{\phi_{\rm c}}(x, y) g_1(x) {\rm d}y {\rm d}x = \int_{\Omega} e^{-\lambda y}g_1(x)  \overline{\phi_{\rm c}}(x, y) {\rm d}x {\rm d}y, 
\end{equation*}
and 
\begin{equation*}
    \int_{\Omega} e^{-\lambda y}g_1(x) {\overline{\phi_0}}(x, y)\;{\rm d}x {\rm d}y = \left( \int_{0}^\infty  e^{-2\lambda y} \; {\rm d}y \right)\left( \int_{\mathcal{E}} g_1(x) \bar{g}(x) \; {\rm d}x  \right) = \frac{1}{2\lambda} \int_{\mathcal{E}} g_1(x) \bar{g}(x) \; {\rm d}x. 
\end{equation*}
Consequently, \eqref{eq:200} reads 
\begin{equation}
\int_{\mathcal E}f(x)\overline{g}(x) \; {\rm d}x +   \int_{\Omega} e^{-\lambda y}g_1(x) {\overline{\phi_0}}(x, y)\;{\rm d}x\, {\rm d}y 
=
-\overline{\int_{\Omega} e^{-\lambda y}g_1(x)  \overline{\phi_{\rm c}}(x, y) {\rm d}x\, {\rm d}y}.
\end{equation}
Taking imaginary parts and using \eqref{eq:imaginary_energy}, together with ${\rm v} = g(x) +\phi_{\rm c}(x, 0)$ gives 
\begin{equation}
    \mathrm{Im}  \int_{\Omega} e^{-\lambda y}g_1(x) {\overline{\phi_0}}(x, y)\;{\rm d}x\, {\rm d}y = - \mathrm{Im} \int_{\mathcal{E}} f(x) \bar{\rm v}(x) \; {\rm d}x.
\end{equation}
Now using the definitions of $\phi_0, g$ and $g_1$ in \eqref{def-gg1}, and noting that $\chi_+$ and $\chi_-$ have disjoint supports in $\mathcal{E}$ with $\mathrm{Im}\bigl((\chi_\pm'' \pm 2i\lambda\chi_\pm')\chi_\pm\bigr) = \pm 2\lambda \chi_\pm'\chi_\pm$, we obtain
\begin{multline*}
     \mathrm{Im}  \int_{\Omega} e^{-\lambda y}g_1(x) {\overline{\phi_0}}(x, y)\;{\rm d}x \,{\rm d}y = \frac{1}{2\lambda} \mathrm{Im} \int_{\mathcal{E}} g_1(x) \overline{g}(x) \; {\rm d}x \\
     = |d_{+}|^2 \int_{\mathbb{R}} \chi_{+}'(x) \chi_{+}(x) \; {\rm d}x - |d_{-}|^2 \int_{\mathbb{R}} \chi_{-}'(x) \chi_{-}(x) \; {\rm d}x  = \frac{1}{2} (|d_{+}|^2 + |d_{-}|^2).
\end{multline*}
Therefore, 
\begin{equation*}
    |d_{+}|^2 + |d_{-}|^2 = - 2\mathrm{Im} \int_{\mathcal{E}} f(x) \bar{\rm v}(x) \; {\rm d}x \leqslant 2 \norm{\rm v}_{L^2_{-s}(\mathcal{E})} \|f\|_{L^2_s(\mathcal{E})}.
\end{equation*}
This completes the proof of the proposition. 
\end{proof} 

Since $R_\Omega^-(\lambda)f = \overline{R_\Omega^+(\lambda)\overline{f}}$ and the identity $\mathcal{R}_\Omega(z) = \mathfrak{E}R_\Omega(z)\mathfrak{E}^*$ extends to the limit as $z \to \lambda \pm i0,$ we obtain $\mathcal{R}_\Omega^\pm(\lambda) = \mathfrak{E}R_\Omega^\pm(\lambda)\mathfrak{E}^*$. 
Consequently, \cref{prop-lambda4} yields the following estimate:

\begin{corollary} \label{cor-lambda4}
There exists a constant $C > 0,$ depending only on $s$ and $\Omega,$ and independent of $\lambda > 0$ such that 
\begin{equation}
   \norm{ \mathcal{R}_{\Omega}^{\pm}(\lambda)}_{\mathcal{L}(L^2_s(\mathbb{R}),L^2_{-s}(\mathbb{R}))} \leqslant C (1 + \lambda^{-4})  \qquad\qquad \left(1 \leqslant s < \frac{3}{2}, \lambda > 0\right).
\end{equation}
\end{corollary}

The main result of this section is:

\begin{proposition} \label{prop-Holder-ROmega}
Let $1 \leqslant s<\frac32,$  $0<\beta< s-\frac12$ and $0 < a < 1 < b < \infty.$ Let $R_{\Omega}^\pm(\lambda) \in \mathcal{L}(L^2_s(\mathcal{E}), L^2_{-s}(\mathcal{E}))$ be the operator defined in  \eqref{eq:ROmega-boundary}.
There is a constant $C_{s,\beta}>0$, independent of $a$ and $b$, such that
\begin{equation}\label{eq:sup-ROmega}
 \sup_{\lambda\in[a,b]}
 \|R_{\Omega}^{\pm}(\lambda)\|
 _{\mathcal L(L_s^2(\mathcal{E}),L_{-s}^2(\mathcal{E}))}
 \leqslant C_{s,\beta} a^{-4},
\end{equation}
and, for every $\lambda,\mu\in[a,b]$,
\begin{equation}\label{eq:holder-ROmega}
\| R_{\Omega}^{\pm}(\lambda)
       - R_{\Omega}^{\pm}(\mu)\|
 _{\mathcal L(L_s^2(\mathcal{E}),L_{-s}^2(\mathcal{E}))} \leqslant C_{s,\beta} a^{-(8+\beta)} b^4
 |\lambda-\mu|^\beta.
\end{equation}

\end{proposition}
\begin{proof}
Throughout  the proof, the constant $C$ depends only on $s$ and $ \beta$, and is independent of $a$ and $b$.  The constant may change from line to line.  

We first note that \eqref{eq:sup-ROmega} follows from \cref{prop-lambda4}. To prove \eqref{eq:holder-ROmega}, we first prove a similar estimate for the operator   $\mathcal{R}_{\Omega}^{\pm}(\lambda) \in \mathcal{L}(L^2_s(\mathbb{R}), L^2_{-s}(\mathbb{R}))$ defined in \eqref{eq:ROmega-boundary-ext}. 
We recall that, for $\lambda > 0,$ $\mathcal{R}_{\Omega}^{\pm}(\lambda)$ is given  by 
\begin{equation} \label{eq:1051}
   \mathcal{R}_{\Omega}^{\pm}(\lambda) :=  \left[ \mathcal{R}_{\mathbb{H}}^{\pm}(\lambda) + \iota K_1^{\pm}(\lambda) \right] \left[ F_1^{\pm} (\lambda)\right]^{-1},
\end{equation}
where the operators $\mathcal{R}_{\mathbb H}^{\pm}(\lambda),$ $K_1^{\pm}(\lambda)$ and $F_1^{\pm}(\lambda)$ are defined in \cref{prop:free-lap}, \eqref{eq:def_K1} and \eqref{eq:def_F1}, respectively. Moreover, $\iota:L^2_s(\mathbb{R})\hookrightarrow L^2_{-s}(\mathbb{R})$ denotes  the natural inclusion map.

 Since $0 < a < 1,$ from \cref{prop-holder-RH},  we obtain  
    \begin{equation} \label{eq-948}
        \norm{\mathcal{R}_{\mathbb{H}}^{\pm} (\lambda)}_{\mathcal L(L_s^2(\mathbb{R}),L_{-s}^2(\mathbb{R}))} \leqslant C(1+ {\rm log}(1 +a^{-1})) \leqslant C a^{-\beta} \qquad\qquad (\lambda \in [a,b]),
    \end{equation}
    and 
    \begin{equation} \label{eq-948-bis}
\| \mathcal{R}_{\mathbb{H}}^{\pm}(\lambda)
       - \mathcal{R}_{\mathbb{H}}^{\pm}(\mu)\|
 _{\mathcal L(L_s^2(\mathbb R),L_{-s}^2(\mathbb R))} \leqslant C a^{-\beta}
 |\lambda-\mu|^\beta \qquad (\lambda, \mu \in [a,b]).
\end{equation}
From the definitions of $F_1^{\pm}(\lambda)$ and $D_1$ in \eqref{eq:def_F1} and \eqref{eq:compact_difference_at_minus_one}, we see that, by some simple algebric calculation, $F_1^{\pm}(\lambda)$ and its inverse  can  written as 
\begin{align}
& F_1^\pm(\lambda) = \mathbb{I} -(\lambda+1)D_1 \iota -(\lambda+1)^2D_1 \mathcal{R}^{\pm}_{\mathbb{H}}(\lambda) \quad (\lambda > 0) \label{def-F1new}
\\
 & [F_1^\pm(\lambda)]^{-1}
= \mathbb{I} +(\lambda+1)D_1 \iota 
+(\lambda+1)^2D_1 \mathcal{R}_\Omega^\pm(\lambda) \quad (\lambda > 0).
\end{align}
Therefore, by \cref{prop:compact_perturbation} and \cref{cor-lambda4}, we deduce 
\begin{equation}\label{eq-sup-F1inv}
\sup_{\lambda\in[a,b]}\|[F_1^\pm(\lambda)]^{-1}\|_{
\mathcal L(L_s^2(\mathbb R))}
\leqslant C\sup_{\lambda\in[a,b]}
\bigl(1+(\lambda+1)+(\lambda+1)^2(1+\lambda^{-4})\bigr)
\leqslant Ca^{-4}b^2.
\end{equation}
We now estimate the H\"older seminorm defined in
\eqref{eq:holder_norms} of $F_1^{\pm}(\lambda)$ with respect to $\lambda.$  From \eqref{def-F1new}, for any $\lambda, \mu \in [a, b],$ we get 
\begin{multline*}
F_1^{\pm}(\lambda)-F_1^{\pm}(\mu)
=-(\lambda-\mu)D_1 \iota 
-(\lambda+1)^2D_1 \left(\mathcal{R}_{\mathbb{H}}^{\pm}(\lambda)-\mathcal{R}_{\mathbb{H}}^{\pm}(\mu) \right) \\
-(\lambda-\mu)(\lambda+\mu+2)D_1 \mathcal{R}_{\mathbb{H}}^{\pm}(\mu) \qquad\qquad(\lambda,\mu\in [a,b]).
\end{multline*}
Consequently, by \eqref{eq-948}, \eqref{eq-948-bis} and \cref{prop:compact_perturbation},  for any $\lambda, \mu \in [a, b],$ we have 
\begin{multline*}
   \norm{F_1^\pm(\lambda) - F_1^\pm(\mu)}_{\mathcal L(L_s^2(\mathbb{R}),L_{s}^2(\mathbb{R}))} \leqslant 
 C \Big( |\lambda - \mu| + a^{-\beta} (\lambda + 1)^2  |\lambda - \mu|^\beta 
 \\ + a^{-\beta}  |\lambda - \mu| |\lambda + \mu +2| \Big).
\end{multline*}
Since $a < 1 < b,$ and $\beta < 1,$ we deduce  
\begin{multline} \label{eq-holder-F1}
  \norm{F_1^\pm(\lambda) - F_1^\pm(\mu)}_{\mathcal L(L_s^2(\mathbb{R}),L_{s}^2(\mathbb{R}))} \leqslant 
 C |\lambda - \mu|^{\beta} \left( |\lambda - \mu|^{1-\beta}  + 4a^{-\beta}b^2 + 4a^{-\beta}b|\lambda - \mu|^{1-\beta}\right) \\
 \leqslant C |\lambda - \mu|^{\beta} \left( b^{1-\beta}+ a^{-\beta}b^2+ a^{-\beta} b^{2 -\beta} \right)
 \leqslant C a^{-\beta}b^2 |\lambda - \mu|^\beta \qquad (\lambda, \mu \in [a, b]).
\end{multline}
We set 
\begin{equation}
    P^{\pm}(\lambda) := \mathcal{R}_{\mathbb{H}}^{\pm}(\lambda) + \iota K_1^{\pm}(\lambda) \in \mathcal L(L_s^2(\mathbb{R}),L_{-s}^2(\mathbb{R})) \qquad\qquad (\lambda > 0).
\end{equation}
Indeed, \eqref{eq:def_K1} gives
\[
 P^\pm(\lambda)=\mathcal R_{\mathbb H}^\pm(\lambda)
       +\iota D_1\iota
       +(\lambda+1)\iota D_1\mathcal R_{\mathbb H}^\pm(\lambda).
\]
Consequently, its difference at $\lambda$ and $\mu$ is
\[
 \begin{aligned}
 P^\pm(\lambda)-P^\pm(\mu)
 ={}&[\mathbb I+(\lambda+1)\iota D_1]
       [\mathcal R_{\mathbb H}^\pm(\lambda)
        -\mathcal R_{\mathbb H}^\pm(\mu)]\\
 &+(\lambda-\mu)\iota D_1\mathcal R_{\mathbb H}^\pm(\mu).
 \end{aligned}
\]
Using \eqref{eq-948}--\eqref{eq-948-bis} and
$|\lambda-\mu|\leqslant b^{1-\beta}|\lambda-\mu|^\beta$,
we obtain the following bounds (using also $a^{-\beta}\leqslant a^{-1}$):

\begin{equation*}
    \sup_{\lambda \in [a, b]} \norm{P^{\pm}(\lambda)}_{\mathcal L(L_s^2(\mathbb{R}),L_{-s}^2(\mathbb{R}))} \leqslant Cba^{-1}, 
\end{equation*}
and 
\begin{equation} \label{eq-holde-Slambda}
  \norm{ P^{\pm}(\lambda) - P^{\pm}(\mu) }_{\mathcal L(L_s^2(\mathbb{R}),L_{-s}^2(\mathbb{R}))} \leqslant C a^{-\beta}b  |\lambda - \mu|^\beta \qquad\qquad (\lambda, \mu \in [a, b]).
\end{equation}
Note that from \eqref{eq:1051}, we have $\mathcal{R}_{\Omega}^{\pm}(\lambda) = P^{\pm}(\lambda)\left[ F_1^{\pm} (\lambda)\right]^{-1}.$ This implies 
\begin{equation*}
    \mathcal{R}_{\Omega}^{\pm}(\lambda) - \mathcal{R}_{\Omega}^{\pm}(\mu) = \left( P^{\pm}(\lambda) - P^{\pm}(\mu) \right) \left[ F_1^{\pm} (\lambda)\right]^{-1} -  \mathcal{R}_{\Omega}^{\pm}(\mu) \left( F_1^\pm(\lambda) - F_1^\pm(\mu) \right) \left[ F_1^{\pm} (\lambda)\right]^{-1}.
\end{equation*}
Consequently, by \cref{cor-lambda4}, \eqref{eq-holde-Slambda}, \eqref{eq-holder-F1} and \eqref{eq-sup-F1inv}, we deduce that 
 for any $\lambda, \mu \in [a, b],$ 
 \begin{multline}
    \norm{ \mathcal{R}_{\Omega}^{\pm}(\lambda)
       - \mathcal{R}_{\Omega}^{\pm}(\mu)}
 _{\mathcal L(L_s^2(\mathbb R),L_{-s}^2(\mathbb R))} \\
 \leqslant C |\lambda - \mu|^{\beta} \left( a^{-(4+\beta)}b^3  + a^{-(8+\beta)} b^4 \right) \leqslant  C  a^{-(8+\beta)} b^4  |\lambda - \mu|^{\beta}.
 \end{multline}
 Finally, \eqref{eq:holder-ROmega} follows from the above estimate and \eqref{eq:ROmega-boundary}.

\end{proof}

\section{Proof of \cref{local_energy_decay}}\label{sec:local_energy_decay}

Within this section we denote by $H$ the Hilbert space $L^2(\mathcal{E})$.
Moreover, let $A:\mathcal{D}(A)\to H$ be the operator defined by 
\begin{equation}\label{prop_a_1}
 \mathcal{D}(A)=H^1(\mathcal E),
\end{equation}
\begin{equation}\label{prop_a_2}
  A {\rm v}=g \Lambda_\Omega {\rm v} \qquad\qquad({\rm v}\in \mathcal{D}(A)),  
\end{equation}
where $\Lambda_\Omega$
is the Dirichlet-to-Neumann operator introduced in \Cref{prop:DtN_L2}. Thus, according to
\Cref{thm:self-adjoint}, we have  that $A$
is self-adjoint, injective and non-negative on $L^2(\mathcal{E})$. Moreover, still from \Cref{thm:self-adjoint},
we know that
\begin{equation}\label{eq:domains_A}
 \mathcal{D}(A^{\frac{1}{2}})=H^{\frac{1}{2}}(\mathcal E),\qquad
 \|A^{\frac{1}{2}} {\rm v}\|_{L^2(\mathcal{E})}^2
   =g\|\nabla (D_\Omega {\rm v})\|_{L^2(\Omega)}^2 \qquad\qquad({\rm v}\in H^1(\mathcal E)).
\end{equation}

\begin{lemma}
\label{lem:local_energy_compactness}
Let $\mathscr W_\Omega$ be the space introduced in \eqref{eq:robin_energy_space} and let $D_\Omega\in \mathcal{L}\left(H^1(\mathcal{E}),\mathscr W_\Omega\right)$ be the Dirichlet map introduced in \Cref{prop:h1lifting}.
Then $D_\Omega$ admits a unique extension, denoted in the same way, to an operator in $\mathcal{L}(H^{\frac12}(\mathcal E),\mathscr W_\Omega)$. Moreover, for each $R>0$, the restriction map $C_{1,R} f = f|_{\mathcal{E}_R}$ defines an operator $C_{1,R}\in \mathcal{L}\left(H,L^2(\mathcal{E}_R)\right)$ of norm at most one. Furthermore, the linear map defined on the dense subspace $\operatorname{Ran}(A^{\frac12})$ by 
\begin{equation}\label{eq:physical_local_energy_maps}
    A^{\frac12}{\rm u} \longmapsto \sqrt g\,\nabla (D_\Omega {\rm u})|_{\Omega_R} \qquad ({\rm u}\in \mathcal{D}(A^{\frac12})),
\end{equation}
is a contraction, and thus extends uniquely to an operator $C_{0,R}\in \mathcal{L}\left(H, L^2(\Omega_R)^2\right)$. Finally,
\begin{equation}\label{eq:physical_spectral_compactness}
    C_{j,R}E_A([a,b])\hbox{ is compact}
       \qquad(0<a<b<\infty,\ j=0,1),
\end{equation}
where $E_A$ is the spectral measure associated with $A$ (see \cref{sec:lap_background}).
\end{lemma}

\begin{proof}
By \eqref{eq:domains_A} and the fact that the trace of
$D_\Omega {\rm u}$ on $\mathcal E$ is ${\rm u}$, we have
\[
 \|D_\Omega {\rm u}\|_{\mathscr W_\Omega}^2
 =\|{\rm u}\|_H^2+g^{-1}\|A^{\frac{1}{2}}{\rm u}\|_H^2
 \qquad({\rm u}\in\mathcal{D}(A)).
\]
The above estimate, combined with the density of $\mathcal{D}(A)=H^1(\mathcal E)$ in
$\mathcal{D}(A^{\frac{1}{2}})=H^{\frac{1}{2}}(\mathcal E)$, implies that indeed
$D_\Omega$ admits a unique extension to an operator in
 $\mathcal{L}(H^{\frac{1}{2}}(\mathcal E),\mathscr W_\Omega)$. Moreover, for every
${\rm u}\in H^{\frac{1}{2}}(\mathcal E)$ we have
\begin{equation}\label{eq:form_energy_lifting_identity}
 g\|\nabla D_\Omega {\rm u}\|_{L^2(\Omega)^2}^2
       =\|A^{\frac{1}{2}}{\rm u}\|_H^2,
 \qquad
 \|D_\Omega {\rm u}\|_{\mathscr W_\Omega}^2
       =\|{\rm u}\|_H^2+g^{-1}\|A^{\frac{1}{2}}{\rm u}\|_H^2.
\end{equation}
Using \eqref{eq:energy_controls_local_H1} it follows that for every sufficiently large $M$ we have
\begin{equation}\label{eq:form_lifting_local_bound}
 \|D_\Omega {\rm u}\|_{H^1(\Omega\cap B_M)}
 \leqslant C_M\bigl(\|{\rm u}\|_H+
                              \|A^{\frac{1}{2}}{\rm u}\|_H\bigr) \qquad\qquad({\rm u}\in H^{\frac{1}{2}}(\mathcal E)).
\end{equation}

The operator $C_{1,R}$ is the restriction to $\mathcal E_R$ so it clearly belongs to $\mathcal{L}\left(H,L^2(\mathcal E_R)\right)$ and its operator norm is at most $1$.

Concerning $C_{0,R}$, we first remark that, due to the injectivity of $A^{\frac{1}{2}}$, this operator
is well defined on 
$\operatorname{Ran} A^{\frac{1}{2}}\subset H$. Moreover, \eqref{eq:physical_local_energy_maps}
and \eqref{eq:form_energy_lifting_identity} give
\begin{equation}\label{eq:estimate_C0R}
 \|C_{0,R}(A^{\frac{1}{2}}{\rm u})\|_{L^2(\Omega_R)^2}
                 \leqslant\|A^{\frac{1}{2}}{\rm u}\|_H \qquad\qquad ({\rm u}\in \mathcal{D}(A^{\frac{1}{2}})).
\end{equation}

Since $A$ is self-adjoint and injective, its square root $A^{\frac{1}{2}}$ is also injective, ensuring that $\operatorname{Ran}A^{\frac{1}{2}}$ is dense in $H$. Combining this fact and
\eqref{eq:estimate_C0R}, it follows that indeed
$C_{0,R}$ extends to a linear contraction from $H$ to $L^2(\Omega_R)^2$.

We still have to prove \eqref{eq:physical_spectral_compactness}. We begin with the surface observation, i.e., we show that \eqref{eq:physical_spectral_compactness} holds for $j=1$.
Fix $0<a<b<\infty$ and write $P=E_A([a,b])$.

Let $(f_n)$ be a sequence weakly convergent to zero in $H$. The map
$P$ belongs to $\mathcal{L}(H, \mathcal{D}(A^{\frac{1}{2}}))$ since
\[
    \|Pf\|_H\leqslant\|f\|_H,
    \qquad\|A^{\frac{1}{2}}Pf\|_H\leqslant\sqrt b\,\|f\|_H.
\]
It follows that $Pf_n\to 0$ weakly in $\mathcal{D}(A^{\frac{1}{2}})$.
Choose a sufficiently large $M>R$ so that
$\Omega\cap B_M$ is a bounded Lipschitz domain.
Estimate~\eqref{eq:form_lifting_local_bound} implies next that
\[
    D_\Omega(Pf_n)\to 0
            \quad\hbox{weakly in }H^1(\Omega\cap B_M).
\]
The trace map from $H^1$ of a bounded Lipschitz domain to
$L^2$ of its boundary is compact. Restricting the trace to
$\mathcal E_R$ therefore gives
$\|C_{1,R}Pf_n\|_{L^2(\mathcal E_R)}\to0$.
This proves compactness of $C_{1,R}P$.

It remains to prove compactness of $C_{0,R}P$.
For the same sequence $(f_n)$ weakly convergent to zero, set
\[
    {\rm u}_n= \int_a^b \lambda^{-\frac{1}{2}} \, {\rm d}E_A(\lambda) f_n,\qquad U_n=D_\Omega {\rm u}_n.
\]
The spectral integral defining ${\rm u}_n$ yields a bounded operator from $H$ to $\mathcal{D}(A)$ with respect to the graph norm, since
\[
    \|{\rm u}_n\|_H\leqslant a^{-\frac{1}{2}}\|f_n\|_H,
    \qquad\|A {\rm u}_n\|_H\leqslant\sqrt b\,\|f_n\|_H.
\]
Thus ${\rm u}_n \to 0$ weakly in $\mathcal{D}(A)$, the sequence
$\Lambda_\Omega {\rm u}_n$ is bounded in $H$, and
$\nabla U_n$ is bounded in $L^2(\Omega)$ by
\eqref{eq:form_energy_lifting_identity}. On every sufficiently
large fixed truncation, \eqref{eq:form_lifting_local_bound}
also gives weak convergence of $U_n$ to zero in $H^1$.
By the Rellich--Kondrachov theorem and compactness of the trace map, we extract a subsequence strongly convergent in local $L^2$ and local boundary traces. By uniqueness of the weak limit, the entire sequence converges strongly to zero. Thus,
\begin{equation}\label{eq:compact_lifting_and_trace_limits}
    U_n\to 0\ \hbox{in }L^2(\Omega\cap B_M),
    \qquad
    {\rm u}_n\to 0\ \hbox{in }L^2(\mathcal E\cap(-M,M)).
\end{equation}

Choose a real $\eta\in C_c^\infty(B_M)$ equal to $1$ on $B_R$. The function $\eta^2U_n$ has an $L^2$ boundary trace, and
its gradient is $\eta^2\nabla U_n+2\eta U_n\nabla\eta$,
which belongs to $L^2(\Omega)^2$ by the local $H^1$ bound and
the compact support of $\eta$. Thus $\eta^2U_n\in\mathscr W_\Omega$.
Since ${\rm u}_n\in\mathcal{D}(A)=H^1(\mathcal E)$, we may apply
\eqref{eq:operatorgreen} with ${\rm u}={\rm u}_n$ and test function
$\eta^2U_n$. Taking real parts and expanding its gradient yields
\begin{equation}\label{eq:localized_gradient_identity}
\begin{split}
    \int_\Omega\eta^2|\nabla U_n|^2
       ={}&\operatorname{Re}\int_{\mathcal E}
             \eta^2(\Lambda_\Omega {\rm u}_n)\overline{{\rm u}_n}\\
          &-2\operatorname{Re}\int_\Omega
             \eta\,\overline{U_n}\,\nabla U_n\cdot\nabla\eta.
\end{split}
\end{equation}
Note that
\[
    \left|\operatorname{Re}\int_{\mathcal E}
             \eta^2(\Lambda_\Omega {\rm u}_n)\overline{{\rm u}_n}\,{\rm d}x\right|\leqslant \|\eta\|_{L^\infty(B_M)}^2
       \|\Lambda_\Omega {\rm u}_n\|_H
       \|{\rm u}_n\|_{L^2(\mathcal E\cap(-M,M))},
\]
and
\[
  \left| 2\operatorname{Re}\int_\Omega
             \eta\,\overline{U_n}\,\nabla U_n\cdot\nabla\eta\,{\rm d}x\,{\rm d}y \right| \leqslant  2\|\eta\nabla\eta\|_{L^\infty(B_M)}
       \|U_n\|_{L^2(\Omega\cap B_M)}
       \|\nabla U_n\|_{L^2(\Omega)}.
\]
Both bounds tend to zero by
\eqref{eq:compact_lifting_and_trace_limits}. Since $\eta=1$ on
$B_R$, we obtain
\[
    \|C_{0,R}Pf_n\|_{L^2(\Omega_R)^2}^2
        =g\int_{\Omega_R}|\nabla U_n|^2\to 0.
\]
This proves compactness of $C_{0,R}P$, which ends the proof.
\end{proof}
We can now apply the abstract decay result (\cref{prop:abstract_local_energy_decay}) to the water-wave equation.

\begin{proof}[Proof of \cref{local_energy_decay}]
We first establish the absolute continuity of the spectrum. \cref{lem:no_fixed_eigenvalues} excludes positive eigenvalues, and the zero eigenvalue was excluded in \cref{thm:self-adjoint}. Since \cref{thm:mainlap} holds for any $s \in [1, 3/2)$, we fix such an $s$ to obtain the locally uniform LAP for $\Lambda_\Omega$ on $(0,\infty)$. Since $A=g\Lambda_\Omega$ with $g>0$, the same locally uniform LAP holds for $A$ on $(0,\infty)$. The weighted space $L^2_s(\mathcal E)$ is dense in $H=L^2(\mathcal E)$, and ${\rm Ker}\, A=\{0\}$. The Stone formula and the density argument in Remark~\ref{rem:lap_to_abstract_decay} therefore show that $A$ has purely absolutely continuous spectrum (see Amrein, Boutet de Monvel and Georgescu \cite[Proposition~7.1.3]{amrein1996c0}).

Our next aim is to apply this spectral property to the energy variables. The form-domain identification in \cref{thm:self-adjoint} yields $\mathcal{D}(A^{\frac{1}{2}})=H^{\frac{1}{2}}(\mathcal E)$. Thus, the initial data in the statement satisfy the hypotheses required in \cref{prop:abstract_local_energy_decay}. The solution supplied by \cref{thm:well_posedness_energy} satisfies $v_{tt}+Av=0$. The abstract decay result (\cref{prop:abstract_local_energy_decay}) relies on the Riemann-Lebesgue lemma applied to the absolutely continuous spectral measure to ensure the weak decay of the half-wave groups. Consequently, \eqref{eq:abstract_weak_energy_decay} yields
\[
    (A^{\frac{1}{2}}v)(t,\cdot) \to 0 \quad \text{and} \quad v_t(t,\cdot) \to 0 \quad \text{weakly in } H \text{ as } |t|\to\infty.
\]
Since $A^{\frac{1}{2}}=\sqrt g\,\Lambda_\Omega^{\frac{1}{2}}$, this proves the weak convergence assertion.

Finally, to deduce the decay of the physical local energy, we fix $R>0$. \cref{lem:local_energy_compactness} provides the bounded operators $C_{0,R}$ and $C_{1,R}$, and verifies their compactness after every spectral truncation $E_A([a,b])$ with $0<a<b<\infty$. Applying \cref{prop:abstract_local_energy_decay} with these two operators, it follows that
\[
    \|(C_{0,R}A^{\frac{1}{2}}v)(t,\cdot)\|_{L^2(\Omega_R)^2} + \|(C_{1,R}v_t)(t,\cdot)\|_{L^2(\mathcal E_R)} \longrightarrow 0 \qquad (|t| \to \infty).
\]
From the definitions in \eqref{eq:physical_local_energy_maps}, we see that
\begin{equation} \label{def-ER-new}
        E_R(t)=\frac{1}{2}\|(C_{0,R}A^{\frac{1}{2}}v)(t,\cdot)\|_{L^2(\Omega_R)^2}^2 +\frac{1}{2}\|(C_{1,R}v_t)(t,\cdot)\|_{L^2(\mathcal E_R)}^2,
\end{equation}
which implies that $E_R(t)\to0$, proving the local energy decay assertion. The approximation in the abstract proposition is uniform in time with respect to the conserved energy norm. Therefore, this conclusion holds for the stated weak initial data without requiring additional spatial weights or higher regularity assumptions.
\end{proof}

\section{Quantitative local energy decay}\label{sec:quantitative_decay}
We prove \cref{thm:full_polynomial_decay,cor:polynomial_decay_consequences}
by separating the spectrum into low, intermediate, and high frequencies.
The low-frequency part has a small total energy by
\cref{lem:low_spectral_mass}; the high-frequency part has a small total
energy due to the additional Sobolev regularity of the data.
For intermediate frequencies, the H\"older continuity of the resolvent
in \cref{prop-Holder-ROmega} gives the decay through an elementary
Fourier estimate. Conservation of energy controls the two tails
uniformly in time. We finally let the lower threshold tend to zero
and the upper threshold tend to infinity at rates chosen to balance
the three estimates.

Throughout this section, we write $A=g\Lambda_\Omega$ for the evolution operator appearing
in \eqref{prop_a_1}--\eqref{prop_a_2}. All spectral cutoffs below are
applied to $\Lambda_\Omega$, so that $\lambda$ is its spectral parameter
and $r=\sqrt{g\lambda}$ is the temporal frequency. Equivalently,
$\chi(\Lambda_\Omega)=\widetilde\chi(A)$ with
$\widetilde\chi(\mu)=\chi(\mu/g)$. Hence, the spectral-calculus results
of \cref{appendix:lap} apply with this rescaling. Throughout, 
constants may depend on the fixed $g>0$, the geometry,  and the chosen
cutoff functions, but are independent of the spectral thresholds and time.

\subsection{Decay on an intermediate spectral interval}
We begin with the Fourier estimate used to exploit fractional
regularity in the temporal frequency. The H\"older norm and the
compact-support convention are those of
\eqref{eq:holder_norms}--\eqref{eq:compact_holder_space}.

  \begin{lemma} \label{lem:holder_fourier_decay}
Let $Y$ be a Banach space, $0<\beta<1$, and let
$I\subset\mathbb R$ be a bounded interval. Suppose that
$B\in C_c^{0,\beta}(\mathbb R;Y)$ with
$\operatorname{supp}B\subseteq I$.
Then, for every $t\in\mathbb R$,
\[
\left\|\int_{\mathbb R}e^{itr}B(r)\,{\rm d}r\right\|_Y
\leqslant C_{\beta}(1+ |I|)\,
\|B\|_{C^{0,\beta}(\mathbb{R};Y)}
(1+|t|)^{-\beta},
\]
where $C_{\beta}>0$ depends only on $\beta.$ 
\end{lemma}

\begin{proof}
For $0<|h|\leqslant1$, Hölder continuity and the support condition give
\begin{equation*}
    \int_{\mathbb R}\|B(r+h)-B(r)\|_Y\,{\rm d}r
\leqslant (|I|+1)\norm{B}_{C^{0,\beta}(\mathbb R;Y)}|h|^\beta.
\end{equation*}
A change of variables gives 
\begin{equation}
   ( e^{-ith}-1)\int_{\mathbb R}e^{itr}B(r)\,{\rm d}r = \int_{\mathbb R}e^{itr}\bigl(B(r+h)-B(r)\bigr)\,{\rm d}r. 
\end{equation}
For $|t|\geqslant\pi$, choosing $h=\pi/|t|$ we obtain 
\begin{equation*}
    2\norm{\int_{\mathbb R}e^{itr}B(r)\,{\rm d}r}_Y \leqslant (|I|+1)\pi^\beta
\norm{B}_{C^{0,\beta}(\mathbb R;Y)}|t|^{-\beta}.
\end{equation*}
For $|t|<\pi$, we simply use
\begin{equation*}
    \norm{\int_{\mathbb R}e^{itr}B(r)\,{\rm d}r}_Y  \leqslant \|B\|_{L^1(\mathbb R;Y)}
\leqslant  |I|\|B\|_{C^{0,\beta}(\mathbb R;Y)}.
\end{equation*}
Combining these estimates proves the result.
\end{proof}

Let $\rho_1,\rho_2\in C^\infty(\mathbb R)$ be two fixed functions satisfying
\begin{gather} \label{def-rho12}
    0 \leqslant \rho_1(x), \rho_2(x) \leqslant 1 \mbox{ for all }  x\in \mathbb{R},  \\
   \rho_1(x) =  \begin{dcases}
        0  &  x \in (-\infty,1] \\ 1  &  x \in [2,\infty)
    \end{dcases} \quad \mbox{ and } \quad 
   \rho_2(x) =  \begin{dcases}
        1 &  x \in (-\infty,1] \\ 
        0  & x \in [2,\infty),
    \end{dcases} 
\end{gather}
For $0<\varepsilon<1<N<\infty$, define
\begin{equation} \label{def-835}
\chi_{\varepsilon,N}(\lambda)
:=\rho_1\left(\frac{\lambda}{\varepsilon}\right)
\rho_2\left(\frac{\lambda}{N}\right).
\end{equation}

Recall that the operator $A$ and its domain $\mathcal D(A)$ are defined in \eqref{prop_a_1}--\eqref{prop_a_2}. Let $f_0,f_1\in L_s^2(\mathcal E)$, where $1\leqslant s<\frac32$. Consider the Cauchy problem
\begin{equation}\label{wave-spectral}
    \left\{
    \begin{aligned}
                &w_{tt}(t)+A w(t)=0\qquad(t\in\mathbb{R}),\\
        &w(0)= w_0 := \chi_{\varepsilon,N}(\Lambda_\Omega)f_0 ,\quad w_t(0)=w_1 := \chi_{\varepsilon,N}(\Lambda_\Omega)f_1.
    \end{aligned}
    \right.
\end{equation}
By \cref{est-spectral-cutoff}, we have $w_0,w_1\in\mathcal D(A)$. Consequently, \cref{thm:well_posedness_energy} yields
\begin{equation}
 w\in C(\mathbb R;H^1(\mathcal E))
\cap C^2(\mathbb R;L^2(\mathcal E)).   
\end{equation}
For $R>0$ such that $[x_-,x_+]\subset(-R,R)$, we define the local energy associated with the solution $w$ of \eqref{wave-spectral} by
\begin{equation} \label{def-E-eps-N}
E_R^{\varepsilon,N}(t)
:=\frac{g}{2}\int_{\Omega_R}
|\nabla(D_\Omega w)(t,\cdot)|^2
+\frac12\int_{\mathcal E_R}
|w_t(t,\cdot)|^2.
\end{equation}

We now establish the following local energy decay estimate for $w.$

\begin{proposition} \label{prop:local-spectral-estimate}
Let $1\leqslant s<\frac32$ and $0<\beta<s-\frac12$. For every $R>0$ such that $[x_-,x_+]\subset(-R,R)$, there exists a constant $C_{R,s,\beta}>0$, independent of $\varepsilon$ and $N$, such that, for all $0<\varepsilon<1<N<\infty,$ $f_0, f_1 \in L^2_s(\mathcal{E})$ and $t\in\mathbb R$,
\begin{equation} \label{est:ER-eps-N}
E_R^{\varepsilon,N}(t)
\leqslant C_{R,s,\beta}
(1+|t|)^{-2\beta}
 \varepsilon^{-(16+2\beta)} N^{12+\beta} 
\left(
\norm{f_0}_{L_s^2(\mathcal E)}
+\norm{f_1}_{L_s^2(\mathcal E)}
\right)^2.
\end{equation}
\end{proposition}

\begin{proof}
Throughout  the proof, the constant $C$ depends only on $R, s$ and $ \beta$, and is independent of $\varepsilon$ and $N$.  The constant may change from line to line.

Recall from \eqref{eq:ROmega-boundary} and the proof of
\cref{thm:mainlap} that
\[
 R_\Omega(\lambda\pm i0)
 :=\lim_{\eta\downarrow0}R_\Omega(\lambda\pm i\eta)
 =R_\Omega^\pm(\lambda).
\]
For $\lambda>0$, define the spectral density
\begin{equation}\label{eq:spectral_density_local_decay}
 S(\lambda)=\frac{1}{2\pi i}
       \bigl(R_\Omega(\lambda+i0)-R_\Omega(\lambda-i0)\bigr).
\end{equation}
According to \cref{thm:mainlap}, $S(\lambda) \in \mathcal{L}(L^2_s(\mathcal{E}), L^2_{-s}(\mathcal{E})).$ Note that $\mathrm{supp} \chi_{\varepsilon, N} \subset [\varepsilon, 2N].$  By \cref{rem-sol-rep-cutoff} and Stone's formula, we deduce that the solution $w$ of \eqref{wave-spectral} can be written as 
\begin{equation} \label{exp-w-345}
 w(t) =\int_{\varepsilon}^{2N} \chi_{\varepsilon, N}(\lambda)S(\lambda) \left(  \cos(t\sqrt{g\lambda}) f_0 + \frac{\sin(t\sqrt{g\lambda})}{\sqrt{g\lambda}} f_1 \right) \; {\rm d}\lambda. 
 \qquad\qquad(t\in\mathbb R).
\end{equation}
The above relation holds in $L^2_{-s}(\mathcal{E}),$ see, for instance, \cref{corollary_formula_w_resolvents}. 

By the definition of $C_{0,R}$ in \eqref{eq:physical_local_energy_maps} and of the restriction operator $C_{1,R}$ in \cref{lem:local_energy_compactness}, we have
\begin{equation}
    \sqrt{g} \nabla(D_\Omega w(t))|_{\Omega_R} = C_{0,R} A^{\frac{1}{2}} w(t), \qquad w_t(t)|_{\mathcal{E}_R} = C_{1,R} w_t(t)  \quad (t \in \mathbb{R}).
\end{equation}
For any $f \in L^2(\mathcal{E}),$ we define   
\begin{equation}
    \mathcal{T}_R f = C_{0,R} A^{\frac{1}{2}} R_{\Omega}(-1)f. 
\end{equation}
Since $A=g\Lambda_\Omega$, the spectral calculus gives
$\|A^{\frac12}R_\Omega(-1)\|_{\mathcal L(L^2(\mathcal E))}
\leqslant\sup_{\lambda\geqslant0}\sqrt{g\lambda}/(1+\lambda)
=\sqrt g/2$. Together with the boundedness of $C_{0,R}$ from
\cref{lem:local_energy_compactness}, this yields
$\mathcal{T}_R\in\mathcal{L}(L^2(\mathcal E),L^2(\Omega_R)^2)$. Moreover, from \eqref{eq:physical_local_energy_maps} and \cref{rem-local-estimate-R-Omega}, there exists a constant $C_{R,s} > 0$ such that 
\begin{equation} \label{est-TR}
    \norm{\mathcal{T}_R f}_{ L^2(\Omega_R)^2} \leqslant C_{R,s} \|f\|_{L^2_{-s}(\mathcal{E})} \quad (f \in L^2(\mathcal E)).
\end{equation}
Thus $\mathcal{T}_R$ extends to a bounded operator from $L^2_{-s}(\mathcal{E})$ to $L^2(\Omega_R)^2.$ Since $\mathcal{E}_R$ is bounded, we have
\begin{equation} \label{est-C1R}
    \norm{C_{1, R} f}_{L^2(\mathcal{E}_R)} \leqslant (1 + R^2)^{s/2} \|f\|_{L^2_{-s}(\mathcal{E})}  \qquad (f \in L^2(\mathcal{E})).
\end{equation}
Thus $C_{1,R}$ extends to a bounded operator from $L^2_{-s}(\mathcal{E})$ to $L^2(\mathcal{E}_R).$ 

With the above notation, we see that 
\begin{equation} \label{eq-exp-356}
        \sqrt{g} \nabla(D_\Omega w(t))|_{\Omega_R}  = \mathcal{T}_R (\Lambda_{\Omega } + 1) w(t), \qquad  w_t(t)|_{\mathcal{E}_R} = C_{1,R} w_t(t) \qquad (t \in \mathbb{R}). 
\end{equation}
The spectral calculus and Stone's formula, applied separately to
these compactly supported spectral multipliers, give 
\begin{align*}
 (\Lambda_\Omega+1)w(t)&=\int_{\varepsilon}^{2N}(\lambda+1)\chi_{\varepsilon, N}(\lambda)S(\lambda)
   \left(\cos(t\sqrt{g\lambda})f_0+
            \frac{\sin(t\sqrt{g\lambda})}{\sqrt{g\lambda}}f_1\right)
                                                        \,{\rm d}\lambda,\\
 w_t(t)&=\int_{\varepsilon}^{2N} \chi_{\varepsilon, N}(\lambda)S(\lambda)
   \left(-\sqrt{g\lambda} \sin(t\sqrt{g\lambda} )f_0+
                \cos(t\sqrt{g\lambda})f_1\right) \,{\rm d}\lambda.
\end{align*}
These identities hold in $L^2_{-s}(\mathcal E)$. Indeed, the continuity of the operator-valued spectral densities on the compact spectral interval ensures that the integrals are well-defined as Bochner integrals in this space. Applying the bounded maps
$\mathcal T_R$ and $C_{1,R}$ and changing variables to
$r=\sqrt{g\lambda}$, with ${\rm d}\lambda=(2r/g){\rm d}r$, gives 
\begin{equation} \label{eq:local-energy-924}
    \begin{bmatrix}
         \sqrt{g} \nabla(D_\Omega w(t))|_{\Omega_R} \\ w_t(t)|_{\mathcal{E}_R}
    \end{bmatrix} = \int_{\sqrt{g\varepsilon}}^{\sqrt{2gN}} \chi_{\varepsilon, N}(r^2/g)
                   [\cos(tr)B_1(r)+\sin(tr)B_2(r)]\,{\rm d}r,
\end{equation}
where 
\begin{equation}
    B_1(r) =\frac{2}{g}
\begin{bmatrix}r(\frac{r^2}{g}+1)\mathcal{T}_RS(\frac{r^2}{g})f_0\\
            r C_{1,R}S(\frac{r^2}{g})f_1\end{bmatrix} \quad \mbox{ and } \quad 
          B_2(r) =\frac{2}{g}
       \begin{bmatrix}\dfrac{r^2+g}{g}\mathcal{T}_R S(\frac{r^2}{g})f_1\\  -r^2C_{1,R}S(\frac{r^2}{g})f_0\end{bmatrix}.   
\end{equation}
Write $Y_R=L^2(\Omega_R)^2\times L^2(\mathcal E_R)$ and
$I_{\varepsilon,N}=[\sqrt{g\varepsilon},\sqrt{2gN}]$.
For $j\in\{1,2\}$, define
\begin{equation}
 L_j(r)=
 \begin{cases}
 \chi_{\varepsilon,N}(r^2/g)B_j(r),&r\in I_{\varepsilon,N},\\
 0,&r\notin I_{\varepsilon,N}.
 \end{cases}
\end{equation}
Thus $L_j$ is supported in the positive interval $I_{\varepsilon,N}$.
Then, from \eqref{eq:local-energy-924} and \eqref{def-E-eps-N}, it follows that 
\begin{equation} \label{est-EE-eps-N-bis}
    E^{\varepsilon, N}_R(t) \leqslant C \norm{ \int_{\mathbb{R}}\cos(tr) L_1(r)\,{\rm d}r}^2_{L^2(\Omega_R)^2 \times L^2(\mathcal{E}_R)} + C \norm{ \int_{\mathbb{R}}\sin(tr) L_2(r)\,{\rm d}r}^2_{L^2(\Omega_R)^2 \times L^2(\mathcal{E}_R)}.
\end{equation}
Now, to prove \eqref{est:ER-eps-N}, we use \cref{lem:holder_fourier_decay}. To this end, we now estimate the $C^{0,\beta}$ norm of $L_j(r),$ for $j=1,2.$ It is easy to see that $\mathrm{supp} L_j \subset I_{\varepsilon,N}.$ Let us set 
\begin{equation}
    p_1(r) =  r\left(1+ \frac{r^2}{g}\right), \quad p_2(r) = r, \quad p_3(r) = 1+\frac{r^2}{g}, \quad p_4(r) = -r^2.
\end{equation}
On $I_{\varepsilon,N}$, the polynomials $p_k$ satisfy
$\sup|p_k|\leqslant CN^{3/2},$ and $\sup|p_k'|\leqslant CN$.
For $|r_1-r_2|\leqslant1,$ the mean value theorem gives the required H\"older estimate using the derivative bound, while  for $|r_1-r_2|>1$ the same estimate follows from the supremum bound. This proves
\begin{equation}\label{est-pk}
 \|p_k\|_{C^{0,\beta}(I_{\varepsilon,N})}
 \leqslant CN^{3/2}\qquad\qquad(k=1,2,3,4).
\end{equation}
Similarly, the composed cutoff $r\mapsto\chi_{\varepsilon,N}(r^2/g)$ is bounded by one, and its derivative is bounded by
$C\varepsilon^{-1}\sqrt N$. Applying the same argument, using the derivative bound when $|r_1 - r_2|\leqslant 1$ and the supremum bound otherwise, yields
\begin{equation}\label{est-chi-eps-N}
 \begin{split}
 &\sup_{r\in I_{\varepsilon,N}}
       |\chi_{\varepsilon,N}(r^2/g)|\leqslant1,\\
 &|\chi_{\varepsilon,N}(r_1^2/g)
       -\chi_{\varepsilon,N}(r_2^2/g)|
 \leqslant C\varepsilon^{-1}\sqrt N\,|r_1-r_2|^\beta
 \quad(r_1,r_2\in I_{\varepsilon,N}).
 \end{split}
\end{equation}
For any  $r_1, r_2 \in I_{\varepsilon,N},$ we note that 
\begin{equation*}
    \left| \frac{r_1^2}{g} - \frac{r^2_2}{g}\right| \leqslant C\sqrt{N} |r_1 - r_2|.
\end{equation*}
Consequently, from \cref{prop-Holder-ROmega}, we deduce that 
\begin{equation} \label{est-S-sup}
\sup_{r \in I_{\varepsilon,N}}    \norm{S\left( \frac{r^2}{g}\right)}_{\mathcal{L}(L^2_s(\mathcal{E}),L^2_{-s}(\mathcal{E}))}\leqslant C\varepsilon^{-4}, 
\end{equation}
and 
\begin{equation} \label{est-S-holder}
    \norm{S\left( \frac{r_1^2}{g}\right) -  S\left( \frac{r_2^2}{g}\right)}_{\mathcal{L}(L^2_s(\mathcal{E}),L^2_{-s}(\mathcal{E}))} \leqslant C  \varepsilon^{-(8+\beta)}N^{4 + \beta/2}  |r_1 - r_2|^\beta \qquad (r_1, r_2 \in I_{\varepsilon, N}). 
\end{equation}
We write 
\begin{equation}
    M_s = \left(
\norm{f_0}_{L_s^2(\mathcal E)}
+\norm{f_1}_{L_s^2(\mathcal E)}
\right).
\end{equation}
Combining the estimates \eqref{est-TR}, \eqref{est-C1R}, \eqref{est-S-sup}, \eqref{est-S-holder} and  \eqref{est-pk}, we infer that
\begin{equation}
    \sup_{r \in I_{\varepsilon, N}} \norm{B_j(r)}_{L^2(\Omega_R)^2 \times L^2(\mathcal{E}_R)} \leqslant C \varepsilon^{-4}N^{3/2} M_s, \quad (j =1,2),
\end{equation}
and 
\begin{multline}
  \norm{B_j(r_1) - B_j(r_2)}_{L^2(\Omega_R)^2 \times L^2(\mathcal{E}_R)} 
  \leqslant C \left( \varepsilon^{-(8+\beta)}N^{(11+\beta)/2}  + \varepsilon^{-4}N^{3/2} \right) M_s|r_1-r_2|^\beta \\
 \leqslant  C \varepsilon^{-(8+\beta)}N^{(11+\beta)/2} M_s |r_1-r_2|^\beta \quad (r_1, r_2 \in I_{\varepsilon, N}, j=1,2).
\end{multline}
Consequently, the above estimate together with \eqref{est-chi-eps-N} yields
\begin{equation}
    \sup_{r \in I_{\varepsilon, N}} \norm{L_j(r)}_{L^2(\Omega_R)^2 \times L^2(\mathcal{E}_R)} \leqslant C \varepsilon^{-4}N^{3/2} M_s,
\end{equation}
and 
\begin{multline}
    \norm{L_j(r_1) - L_j(r_2)}_{L^2(\Omega_R)^2 \times L^2(\mathcal{E}_R)} \\
    \leqslant CM_s|r_1 -r_2|^{\beta} \left( \varepsilon^{-(8+\beta)}N^{(11+\beta)/2}+ \varepsilon^{-4}N^{3/2} \varepsilon^{-1}N^{1/2} \right) \\
    \leqslant C M_s   \varepsilon^{-(8+\beta)}N^{(11+\beta)/2}  |r_1 -r_2|^{\beta} \qquad (r_1, r_2 \in I_{\varepsilon, N}, j=1,2).
\end{multline}
Since the cutoff vanishes at both endpoints of $I_{\varepsilon,N}$, so does $L_j$. Its extension by zero outside this interval, therefore, preserves the H\"older estimate. Indeed, if one point lies within the interval and the other outside, the difference can be estimated by comparing the interior value with its value at the intervening endpoint, where $L_j$ vanishes. Consequently,
\[
 \|L_j\|_{C^{0,\beta}(\mathbb R;Y_R)}
 \leqslant C \varepsilon^{-(8+\beta)}N^{(11+\beta)/2} M_s.
\]
We use this bound in \eqref{est-EE-eps-N-bis} and use \cref{lem:holder_fourier_decay} to obtain 
\begin{equation*}
     E^{\varepsilon, N}_R(t)
     \leqslant C (1+|I_{\varepsilon, N}|)^2 (1+|t|)^{-2\beta}  \varepsilon^{-(16+2\beta)}N^{11+\beta}   M_s^2 
     \leqslant C  \varepsilon^{-(16+2\beta)} N^{12+\beta}  (1+|t|)^{-2\beta}M_s^2.
\end{equation*}
This completes the proof of the proposition. 
\end{proof}

\subsection{Removal of the spectral cutoffs}

To prove \cref{thm:full_polynomial_decay}, we now combine the decay estimate on the intermediate spectral interval with uniform bounds on the energies in the low- and high-frequency tails.

\begin{proof}[Proof of \cref{thm:full_polynomial_decay}]
We use the data norm $\mathcal N_s$ defined in \eqref{eq:data_norm_decay}.
    Let $0 < \varepsilon < \frac{1}{2}$ and  $N > 1$. Let us recall the definition of $\rho_1$ and $\rho_2$ from \eqref{def-rho12}, and definition of $\chi_{\varepsilon, N}$ from \eqref{def-835}.  For $\lambda \geqslant 0,$ we set 
    \begin{equation}
       \chi^{\rm l}(\lambda) = 1 - \rho_1 \left(\frac{\lambda}{\varepsilon}\right), \; \chi^{\rm m}(\lambda) = \rho_1\left(\frac{\lambda}{\varepsilon}\right)
\rho_2\left(\frac{\lambda}{N}\right),\;  \chi^{\rm h}(\lambda) = \rho_1 \left(\frac{\lambda}{\varepsilon}\right) \left[1 - \rho_2\left(\frac{\lambda}{N}\right)\right ].
    \end{equation}
    Note that, 
    \begin{equation}
        \chi^{\rm l}(\lambda) +  \chi^{\rm m}(\lambda) + \chi^{\rm h}(\lambda) = 1  \qquad (\lambda \geqslant 0).
        \end{equation}
        For $j = 0,1,$ we set 
        \begin{equation}
            v_j^{\rm l} = \chi^{\rm l}(\Lambda_\Omega)v_j, \;  v_j^{\rm m} = \chi^{\rm m}(\Lambda_\Omega)v_j, \;  v_j^{\rm h} = \chi^{\rm h}(\Lambda_\Omega)v_j.
        \end{equation}
        For ${\rm k} \in \{\rm l, m, h\},$ let $v^{\rm k}\in C(\mathbb R;H^{\frac12}(\mathcal E))\cap C^1(\mathbb R;L^2(\mathcal E)),$ be the solution of 
        \begin{equation} \label{eq-wave-split}
        \begin{dcases}
             v_{tt}^{\rm k} + Av^{\rm k} = 0  \quad (t \in \mathbb{R}) \\
             v^{\rm k}(0) = v^{\rm k}_0, \quad v_t^{\rm k}(0) = v^{\rm k}_1.
        \end{dcases}
        \end{equation}
Then by linearity 
\begin{equation}
    v(t) = v^{\rm l}(t) + v^{\rm m}(t) + v^{\rm h}(t) \quad (t \in \mathbb{R}).
\end{equation}
For ${\rm k} \in \{\rm l, m, h\},$  the corresponding energy $E^{\rm k}(t)$ is given by  
\begin{equation} \label{energy-split} 
    E^{\rm k}(t) = \frac12\|v^{\rm k}_t(t,\cdot)\|_{L^2(\mathcal{E})}^2+
           \frac 12\|(A^{\frac{1}{2}}v^{\rm k})(t,\cdot)\|_{L^2(\mathcal{E})}^2 = \frac12\|v^{\rm k}_1\|_{L^2(\mathcal{E})}^2+
           \frac 12\norm{ A^{\frac{1}{2}}v^{\rm k}_0}_{L^2(\mathcal{E})}^2.
\end{equation}
For ${\rm k} \in \{\rm l, m, h\},$ we define the local energy associated with the solution $v^{\rm k}$ of \eqref{eq-wave-split} by
\begin{equation}
E_R^{\rm k}(t)
:=\frac{g}{2}\int_{\Omega_R}
|\nabla(D_\Omega v^{\rm k})(t,\cdot)|^2
+\frac12\int_{\mathcal E_R}
|v^{\rm k}_t(t,\cdot)|^2.
\end{equation}
We also recall that, the local energy associated with the solution $v$ of \eqref{v_eq}, is denoted by $E_R(t)$ and it is introduced in \cref{local_energy_decay}. Then we have 
\begin{equation} \label{eq-local-decay-all}
    E_R(t) \leqslant 3 \left(E^{\rm l}_R(t) + E^{\rm m}_R(t) + E^{\rm h}_R(t) \right) \leqslant 3 \left(E^{\rm l}(t) + E^{\rm m}_R(t) + E^{\rm h}(t) \right). 
\end{equation}
The low-frequency multiplier is supported in $[0,2\varepsilon]$
and bounded by one. Recall that
$P_{2\varepsilon}=\mathbf1_{[0,2\varepsilon]}(\Lambda_\Omega)$
from \cref{lem:low_spectral_mass}. Since $A=g\Lambda_\Omega$,
\cref{est-spectral-cutoff,lem:low_spectral_mass} give
\begin{equation}\label{eq-decay-low-frequency}
 \begin{split}
 E^{\rm l}(t)
 &\leqslant\frac12\|P_{2\varepsilon}v_1\|_{L^2(\mathcal{E})}^2
                  +g\varepsilon\|v_0\|_{L^2(\mathcal{E})}^2\\
 &\leqslant C\varepsilon\left(
   \log(2/\varepsilon)\|v_1\|_{L^1(\mathcal E)}^2
       +\|v_1\|_{L^2(\mathcal{E})}^2+\|v_0\|_{L^2(\mathcal{E})}^2\right)\\
 &\leqslant C\varepsilon\log(2/\varepsilon)
                       \mathcal N_s(v_0,v_1)^2.
 \end{split}
\end{equation}
Here $2\varepsilon<1$, and the last inequality uses
\[
 \|f\|_{L^1(\mathcal E)}
 \leqslant\left(\int_{\mathcal E}(1+x^2)^{-s}\,{\rm d}x\right)^{1/2}
                      \|f\|_{L_s^2(\mathcal E)},
\]
which holds since $s>1/2$.
To estimate $E^{\rm m}_R(t)$ we note that $\chi^{\rm m} = \chi_{\varepsilon, N}.$  Therefore, by \cref{prop:local-spectral-estimate}, we have 
\begin{equation} \label{eq-decay-middle-frequency}
    E^{\rm m}_R(t) \leqslant C
(1+|t|)^{-2\beta}
 \varepsilon^{-(16+2\beta)} N^{12+\beta}
\left(
\norm{v_0}_{L_s^2(\mathcal E)}
+\norm{v_1}_{L_s^2(\mathcal E)}
\right)^2.
\end{equation}
From the definition of $\chi^{\rm h},$ we see that $\chi^{\rm h}(\lambda) =0$ for $\lambda\leqslant N,$ and $\|\chi^{\rm h}\|_\infty\leqslant1$. Its support corresponds to
$A$-frequencies at least $gN$. Consequently, by
\cref{est-spectral-cutoff}, with the fixed factors of $g$ absorbed in $C$,
\begin{multline} \label{eq-decay-high-frequency}
    E^{\rm h}(t) = \frac{1}{2} \norm{\chi^{\rm h}(\Lambda_\Omega) v_1}^2_{L^2(\mathcal{E})} + \frac{1}{2} \norm{A^{\frac 12} \chi^{\rm h}(\Lambda_\Omega)v_0}_{L^2(\mathcal{E})}^2 
     \leqslant \frac{C}{N} \left( \norm{A^{\frac 12} v_1}_{L^2(\mathcal{E})}^2 + \norm{Av_0}_{L^2(\mathcal{E})}^2 \right).
\end{multline}
The domain identifications in \cref{thm:self-adjoint} imply
$\|Av_0\|_{L^2(\mathcal{E})}\leqslant C\|v_0\|_{H^1}$ and
$\|A^{1/2}v_1\|_{L^2(\mathcal{E})}\leqslant C\|v_1\|_{H^{1/2}}$.
Combining the estimates \eqref{eq-local-decay-all} -- \eqref{eq-decay-high-frequency}, we obtain 
\begin{equation}
    E_R(t)
    \leqslant C \mathcal N_s^2(v_0,v_1) \left( \varepsilon \mathrm{log}\left(\frac{2}{\varepsilon}\right) + (1+|t|)^{-2\beta}
 \varepsilon^{-(16+2\beta)} N^{12+\beta}  + N^{-1} \right).
\end{equation}

We now choose the spectral cutoffs to balance the low-, intermediate-, and high-frequency contributions. To determine the optimal algebraic decay rate, we temporarily disregard the logarithmic factor and set
\begin{equation*}
 \varepsilon=\frac14(1+|t|)^{-a},
\qquad
N=2(1+|t|)^b,   
\end{equation*}
where $a,b>0$. The corresponding algebraic decay exponents are
\begin{equation*}
  a,\qquad
2\beta-(16+2\beta)a-(12+\beta)b,
\qquad b.  
\end{equation*}
Hence, the decay rate is determined by
\begin{equation*}
  \min\left\{ a,b,2\beta-(16+2\beta)a-(12+\beta)b \right\}
\end{equation*}
To obtain decay of order $(1+|t|)^{-\gamma}$, we require $a,b\geqslant\gamma$ and
\begin{equation*}
  \begin{aligned}
\gamma
&\leqslant 2\beta-(16+2\beta)a-(12+\beta)b\\
&\leqslant 2\beta-(28+3\beta)\gamma.
\end{aligned}  
\end{equation*}
It follows that
\begin{equation*}
    \gamma\leqslant\frac{2\beta}{29+3\beta}.
\end{equation*}
The largest admissible exponent is therefore
$\gamma=\frac{2\beta}{29+3\beta},$
which is achieved by choosing $a=b=\gamma$. With this choice, all three contributions have the same algebraic decay rate. This choice gives the best polynomial decay exponent $\gamma$ obtainable from the preceding estimate using power-law cutoffs, while retaining the logarithmic factor in the final bound. Accordingly, we choose 
\begin{equation}
    \varepsilon = \frac{1}{4}(1+|t|)^{-\gamma} \mbox{ and } N = 2(1+ |t|)^\gamma. 
\end{equation}
These choices satisfy $0<\varepsilon<1/2$ and $N>1$, also at
$t=0$. Then, we have 
\begin{equation*}
\varepsilon\mathrm{log}\left(\frac2\varepsilon\right)
=\frac14 (1+|t|)^{-\gamma}
\left(\mathrm{log} 8+\gamma\mathrm{log} (1+|t|) \right)
\leqslant C(1+|t|)^{-\gamma}\mathrm{log}(2+|t|),
\end{equation*}
and 
\begin{equation*}
     (1+|t|)^{-2\beta} \varepsilon^{-(16+2\beta)}N^{12 + \beta} = 2^{44 +5\beta} (1+|t|)^{-2\beta+(28 +3\beta)\gamma} \leqslant C (1+|t|)^{-\gamma}.
\end{equation*}
Consequently, 
\begin{equation*}
    \varepsilon \mathrm{log}\left(\frac{2}{\varepsilon}\right) + (1+|t|)^{-2\beta} \varepsilon^{-(16+2\beta)} N^{12+\beta}  + N^{-1} \leqslant C (1+|t|)^{-\gamma} \ln(2+|t|).
\end{equation*}
Therefore, 
\begin{equation*}
     E_R(t) \leqslant C \mathcal N_s^2(v_0,v_1) (1+|t|)^{-\gamma} \ln(2+|t|).
\end{equation*}
This completes the proof of the theorem. 
\end{proof}

\begin{proof}[Proof of \cref{cor:polynomial_decay_consequences}]
Since $v_0$ and $v_1$ are compactly supported, they belong to $L^2_s(\mathcal{E})$ for any $s\geqslant0.$ We may therefore take $s$ arbitrarily close to $3/2$ and $\beta$ close to $1$. More precisely, 
let $0 < \delta < \frac{1}{16}.$ We choose $0 < \beta < 1$ such that 
\begin{equation*}
    \delta < \gamma = \frac{2\beta}{29+3\beta}.
\end{equation*}
Next choose $s$ satisfying $\max\{1,\beta+1/2\}<s<3/2$.
Since the supports of $v_0$ and $v_1$ are contained in the interval $[-M,M]$, we have, for $j =0,1$, 
\begin{equation*}
    \|v_j\|_{L_s^2(\mathcal E)}
 \leqslant(1+M^2)^{s/2}\|v_j\|_{L^2(\mathcal E)}.
\end{equation*}
Thus, the assumptions of \cref{thm:full_polynomial_decay}  are satisfied. Moreover, our choice of $\beta$ ensures that
\begin{equation*}
  \sup_{t\in\mathbb R}(1+|t|)^{\delta-\gamma}
                         \log(2+|t|)<\infty,   
\end{equation*}
Consequently, \cref{thm:full_polynomial_decay} yields \eqref{eq:bounded_support_decay}, with a constant depending only on $R,M,\delta$, and the fixed parameters of the problem.
\end{proof}

\subsection{Improved decay under higher spectral regularity}
\label{subsec:higher_regularity_decay}
The high-frequency estimate can be strengthened when the data belong
to higher-power domains of $\Lambda_\Omega$. This gives the following
extension, in particular, for smooth compactly supported data satisfying
the stated domain conditions.

\begin{proposition}\label{prop:higher_regularity_decay}
Under Assumption~\ref{ass:geometry}, let $m\geqslant1$,
$1\leqslant s<3/2$ and $0<\beta<s-1/2$. Suppose that
\[
 v_0\in L_s^2(\mathcal E)\cap
       \mathcal D(\Lambda_\Omega^{(m+1)/2}),\qquad
 v_1\in L_s^2(\mathcal E)\cap
       \mathcal D(\Lambda_\Omega^{m/2}),
\]
and set
\begin{equation}\label{eq:higher_data_norm}
 \mathcal N_{s,m}(v_0,v_1)
 =\|v_0\|_{L_s^2(\mathcal E)}+\|v_1\|_{L_s^2(\mathcal E)}
  +\|\Lambda_\Omega^{(m+1)/2}v_0\|_{L^2(\mathcal{E})}
  +\|\Lambda_\Omega^{m/2}v_1\|_{L^2(\mathcal{E})}.
\end{equation}
For every $R>0$ with $[x_-,x_+]\subset(-R,R)$, the solution of
\eqref{v_eq} satisfies
\begin{equation}\label{eq:higher_regularity_decay}
 E_R(t)\leqslant C_{R,s,\beta,m}\mathcal N_{s,m}(v_0,v_1)^2
 (1+|t|)^{-\gamma_m}\log(2+|t|),
\end{equation}
where 
\begin{equation*}
\gamma_m = 
    \frac{2\beta m}{(17+2\beta)m+12+\beta}.
\end{equation*}
If, in addition, both data vanish for $|x|>M$, then, for every
$0<\delta<2m/(19m+13)$,
\begin{multline}\label{eq:higher_bounded_support_decay}
 E_R(t)\leqslant C_{R,M,m,\delta}(1+|t|)^{-\delta}\bigl(
\|v_0\|_{L^2(\mathcal{E})}^2+\|\Lambda_\Omega^{(m+1)/2}v_0\|_{L^2(\mathcal{E})}^2\\
+\|v_1\|_{L^2(\mathcal{E})}^2+\|\Lambda_\Omega^{m/2}v_1\|_{L^2(\mathcal{E})}^2\bigr).
\end{multline}
\end{proposition}

\begin{proof}
For a non-negative self-adjoint operator, the spectral theorem gives
$\mathcal D(\Lambda_\Omega^b)\subset\mathcal D(\Lambda_\Omega^a)$
for $0\leqslant a\leqslant b$, with continuous inclusion for the graph
norms, since $\lambda^{2a}\leqslant1+\lambda^{2b}$.
Thus the data have the regularity required for well-posedness.
Use the decomposition from the proof of
\cref{thm:full_polynomial_decay}, with cutoffs applied to
$\Lambda_\Omega$. The low-frequency estimate
\eqref{eq-decay-low-frequency} and the intermediate-frequency estimate
\eqref{eq-decay-middle-frequency} remain valid with
$\mathcal N_{s,m}$ in place of $\mathcal N_s$.

For the high-frequency part, write $E_{\Lambda_\Omega}$ for the
spectral measure of $\Lambda_\Omega$, as in Appendix~\ref{appendix:lap}.
The support of $\chi^{\rm h}$ lies in $[N,\infty)$ and
$|\chi^{\rm h}|\leqslant1$. Consequently,
\[
 \begin{aligned}
 \|\chi^{\rm h}(\Lambda_\Omega)v_1\|_{L^2(\mathcal{E})}^2
 &\leqslant N^{-m}\int_{[N,\infty)}\lambda^m
                   \,{\rm d}\|E_{\Lambda_\Omega}(\lambda)v_1\|_{L^2(\mathcal{E})}^2
 \leqslant N^{-m}\|\Lambda_\Omega^{m/2}v_1\|_{L^2(\mathcal{E})}^2,\\
 \|\Lambda_\Omega^{1/2}\chi^{\rm h}(\Lambda_\Omega)v_0\|_{L^2(\mathcal{E})}^2
 &\leqslant N^{-m}\int_{[N,\infty)}\lambda^{m+1}
                   \,{\rm d}\|E_{\Lambda_\Omega}(\lambda)v_0\|_{L^2(\mathcal{E})}^2
 \leqslant N^{-m}\|\Lambda_\Omega^{(m+1)/2}v_0\|_{L^2(\mathcal{E})}^2.
 \end{aligned}
\]
Conservation of energy and $A=g\Lambda_\Omega$ therefore improve
\eqref{eq-decay-high-frequency} to a bound of order $N^{-m}$.
Combining the three parts, using \eqref{eq-local-decay-all}, yields
\begin{equation*}
  E_R(t)\leqslant C\mathcal N_{s,m}(v_0,v_1)^2
 \left( \varepsilon \mathrm{log}\left(\frac{2}{\varepsilon}\right) + (1+|t|)^{-2\beta} \varepsilon^{-(16+2\beta)} N^{12+\beta} + N^{-m} \right).   
\end{equation*}
Take
\begin{equation*}
   \varepsilon=\tfrac14(1+|t|)^{-\gamma_m},\qquad
 N=2(1+|t|)^{\gamma_m/m}.  
\end{equation*}
The definition of $\gamma_m$ gives
\[
 -2\beta+(16+2\beta)\gamma_m
       +\frac{12+\beta}{m}\gamma_m=-\gamma_m.
\]
Thus the intermediate- and high-frequency terms are bounded by
$C(1+|t|)^{-\gamma_m}$, while the low-frequency term is bounded by
$C(1+|t|)^{-\gamma_m}\log(2+|t|)$.
This proves \eqref{eq:higher_regularity_decay}. Similarly, \eqref{eq:higher_bounded_support_decay} follows by the same argument as in the proof of \cref{cor:polynomial_decay_consequences}.
\end{proof}

\begin{remark}\label{rem:higher_domains_decay}
The higher-power assumptions are additional to smoothness and spatial
localization. In particular, \eqref{eq:higher_bounded_support_decay}
applies to $v_0,v_1\in C_c^\infty(\mathcal E)$, provided that they belong to the indicated higher-power domains. We do not assert that all such test
functions satisfy these domain conditions.

For $0\leqslant\theta\leqslant1$, the interpolation argument in the
proof of \cref{thm:self-adjoint}, applied to
$\mathbb I+\Lambda_\Omega$, gives
\[
 \mathcal D(\Lambda_\Omega^\theta)=H^\theta(\mathcal E).
\]
Here the spectral theorem identifies the shifted and unshifted domains,
with equivalent graph norms. The multiplication rule for spectral
integrals then gives the exact characterization
\begin{equation}\label{eq:higher_domain_characterization}
 \mathcal D(\Lambda_\Omega^{1+\theta})
 =\{{\rm u}\in H^1(\mathcal E):
                  \Lambda_\Omega{\rm u}\in H^\theta(\mathcal E)\}.
\end{equation}
Thus, for $1<m\leqslant2$, the domain assumptions in
\cref{prop:higher_regularity_decay} are equivalently
\[
 v_0\in H^1(\mathcal E),\qquad
 \Lambda_\Omega v_0\in H^{(m-1)/2}(\mathcal E),\qquad
 v_1\in H^{m/2}(\mathcal E).
\]
A description using only higher Sobolev regularity of $v_0$ would
require additional contact-point analysis; no such identification is
used here.

For $m=2$, the estimate holds for every decay exponent $\delta < \frac{4}{51}$. If the initial data satisfy the required domain conditions for all $m$, then any exponent $\delta<\frac{2}{19}$ is admissible, with constants depending on the corresponding higher-order norms. We emphasize that this limiting exponent reflects the restrictions of the present method, which improves only the high-frequency contribution, and is not claimed to be optimal for the underlying evolution.
\end{remark}

\section*{Acknowledgments}
The authors thank Christophe Hazard, Jean-Fran\c cois Bony and David Lannes for helpful discussions and suggestions.

\section*{AI disclosure statement}

The main mathematical ideas in this work,  together with the essential part its write-up were generated without the assistance of AI.  During the final stages of this project, GPT-6 Astra and GPT-5.6 Sol
 were  used to prove a small number of technical lemmas, namely, \cref{lem:quarter_plane_comparison} and \cref{lem:low_spectral_mass}. Moreover, GPT-6 Astra has been used for research of references on Sobolev spaces in unbounded domain and abstract versions of the Limiting Absorption Principle.


\appendix
\section{Construction of a Conformal Map}\label{appendix:conformal_mapping}

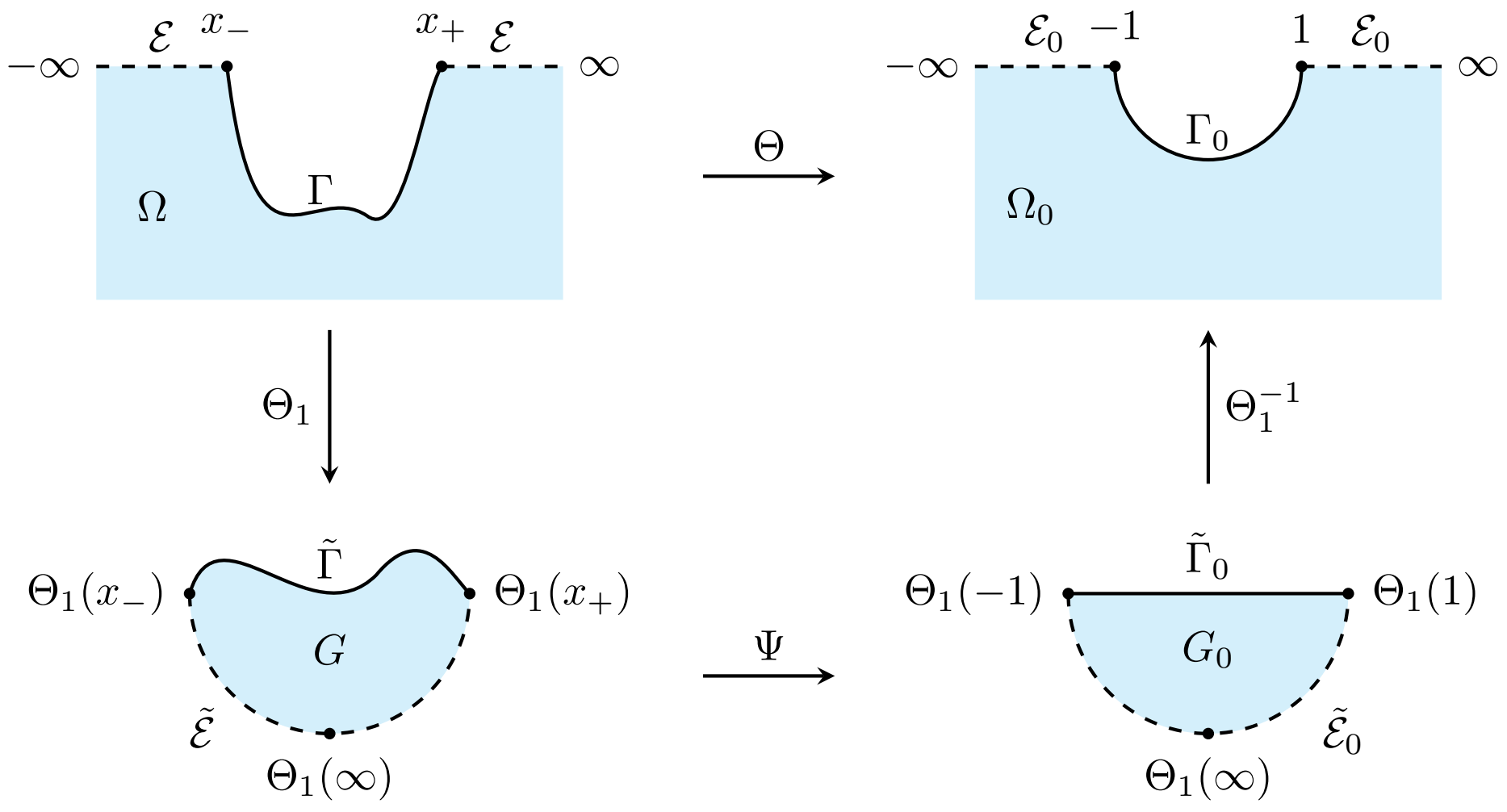
\begin{figure}[H]
        \centering
        \begin{tikzpicture}[scale=0.9, >=stealth]
            
            \begin{scope}[shift={(0,4.8)}, scale=0.85]
                \fill[cyan!15] (-2.5, -2.5) -- (-2.5, 0) -- (-1.1, 0)
                    .. controls (-0.8, -2.5) and (-0.2, -1.2) .. (0.4, -1.6)
                    .. controls (0.8, -1.9) and (1.0, -0.3) .. (1.2, 0)
                    -- (2.5, 0) -- (2.5, -2.5) -- cycle;
                
                \draw[black, thick, dashed] (-2.5, 0) -- (-1.1, 0);
                \draw[black, thick, dashed] (1.2, 0) -- (2.5, 0);
                \node[above] at (-1.8, 0) {$\mathcal{E}$};
                \node[above] at (1.85, 0) {$\mathcal{E}$};
                
                \draw[black, thick] (-1.1, 0) 
                    .. controls (-0.8, -2.5) and (-0.2, -1.2) .. (0.4, -1.6)
                    .. controls (0.8, -1.9) and (1.0, -0.3) .. (1.2, 0);
                \node[below] at (-0.1, -1.0) {$\Gamma$};

                \node at (-1.9, -1.5) {$\Omega$};

                \filldraw (-1.1, 0) circle (1.5pt) node[above=2pt] {$x_-$};
                \filldraw (1.2, 0) circle (1.5pt) node[above=2pt] {$x_+$};
                \node[left] at (-2.5, 0) {$-\infty$};
                \node[right] at (2.5, 0) {$\infty$};
            \end{scope}

            \begin{scope}[shift={(8,4.8)}, scale=0.85]
                \fill[cyan!15] (-2.5, -2.5) -- (-2.5, 0) -- (-1, 0)
                    arc (180:360:1)
                    -- (2.5, 0) -- (2.5, -2.5) -- cycle;
                
                \draw[black, thick, dashed] (-2.5, 0) -- (-1, 0);
                \draw[black, thick, dashed] (1, 0) -- (2.5, 0);
                \node[above] at (-1.75, 0) {$\mathcal{E}_0$};
                \node[above] at (1.75, 0) {$\mathcal{E}_0$};
                
                \draw[black, thick] (-1, 0) arc (180:360:1);
                \node[below] at (0, -0.35) {$\Gamma_0$};

                \node at (-1.9, -1.5) {$\Omega_0$};

                \filldraw (-1, 0) circle (1.5pt) node[above=2pt] {$-1$};
                \filldraw (1, 0) circle (1.5pt) node[above=2pt] {$1$};
                \node[left] at (-2.5, 0) {$-\infty$};
                \node[right] at (2.5, 0) {$\infty$};
            \end{scope}

            \begin{scope}[shift={(0,0)}, scale=0.85]
                \fill[cyan!15] (-1.5, 0) arc (180:360:1.5)
                    .. controls (1.3, 0.2) and (1.0, 0.8) .. (0.5, 0.2)
                    .. controls (-0.2, -0.5) and (-1.2, 1.0) .. (-1.5, 0);
                
                \draw[black, thick, dashed] (-1.5, 0) arc (180:360:1.5);
                \node[below left] at (-1.06, -1.06) {$\tilde{\mathcal{E}}$};
                
                \draw[black, thick] (-1.5, 0) 
                    .. controls (-1.2, 1.0) and (-0.2, -0.5) .. (0.5, 0.2)
                    .. controls (1.0, 0.8) and (1.3, 0.2) .. (1.5, 0);
                \node[above] at (0, 0) {$\tilde{\Gamma}$};
                
                \node at (0, -0.6) {$G$};

                \filldraw (-1.5, 0) circle (1.5pt) node[left=2pt] {$\Theta_1(x_-)$};
                \filldraw (1.5, 0) circle (1.5pt) node[right=2pt] {$\Theta_1(x_+)$};
                \filldraw (0, -1.5) circle (1.5pt) node[below=2pt] {$\Theta_1(\infty)$};
            \end{scope}

            \begin{scope}[shift={(8,0)}, scale=0.85]
                \fill[cyan!15] (-1.5, 0) arc (180:360:1.5) -- (-1.5, 0);
                
                \draw[black, thick, dashed] (-1.5, 0) arc (180:360:1.5);
                \node[below right] at (1.06, -1.06) {$\tilde{\mathcal{E}}_0$};
                
                \draw[black, thick] (-1.5, 0) -- (1.5, 0);
                \node[above] at (0, 0) {$\tilde{\Gamma}_0$};
                
                \node at (0, -0.6) {$G_0$};

                \filldraw (-1.5, 0) circle (1.5pt) node[left=2pt] {$\Theta_1(-1)$};
                \filldraw (1.5, 0) circle (1.5pt) node[right=2pt] {$\Theta_1(1)$};
                \filldraw (0, -1.5) circle (1.5pt) node[below=2pt] {$\Theta_1(\infty)$};
            \end{scope}

            \draw[->, thick] (3.4, 3.8) -- (4.6, 3.8) node[midway, above] {$\Theta$};
            
            \draw[->, thick] (0, 2.4) -- (0, 1.0) node[midway, left] {$\Theta_1$};
            
            \draw[->, thick] (3.4, -0.75) -- (4.6, -0.75) node[midway, above] {$\Psi$};
            
            \draw[->, thick] (8, 1.0) -- (8, 2.4) node[midway, right] {$\Theta_1^{-1}$};

        \end{tikzpicture}
        \caption{Commutative diagram for the conformal mapping $\Theta$, highlighting the boundary correspondence that uniquely determines $\Psi$.}
        \label{fig:conformal_mappings}
\end{figure}

In this appendix, we give a detailed construction of the conformal map used in \cref{sec:spaces}.
Identify the plane with $\mathbb C$ and let
\[
 \Omega_0=\{z:\operatorname{Im}z>0,\ |z|>1\},\quad
 \mathcal E_0=(-\infty,-1)\cup(1,\infty),\quad
 \Gamma_0=\{|z|=1,\ \operatorname{Im}z\geq0\}.
\]
\begin{lemma}\label{lemma:conformal_mapping}
There exists a unique conformal map $\Theta:\Omega\to\Omega_0$
which extends to a homeomorphism from $\overline{\Omega}$ onto $\overline{\Omega}_0$
and satisfies
\[
 \Theta(x_-)=-1,\qquad\Theta(x_+)=1,\qquad\Theta(\infty)=\infty.
\]
It maps $\Gamma$ onto $\Gamma_0$ and each free-surface component onto
the corresponding component of $\mathcal E_0$. 
\end{lemma}
\begin{proof}[Proof of Lemma \ref{lemma:conformal_mapping}]
    We will construct the map $\Theta$ in intermediate steps, as illustrated in Figure \ref{fig:conformal_mappings}. Let $\mathbb{C}_{+}$ denote the upper half plane and let $\Theta_1 : \mathbb{C}_+ \to \mathbb{D}$ be the rational transformation defined by 
    \begin{equation}\label{Theta_1}
        \Theta_1(z) = \frac{z - i}{z + i}.
    \end{equation}
    By the symmetry principle for Möbius transformations (see Ahlfors \cite[Chapter 3, Section 3.3]{ahlfors1979complex}), this mapping sends the real axis to the unit circle and the interior point $i \in \mathbb{C}_+$ to the origin, thereby mapping the upper half-plane conformally onto the unit disk. Since the fluid domain $\Omega$ and the reference domain $\Omega_0$ are simply connected regions strictly contained within $\mathbb{C}_+$, their images under $\Theta_1$, denoted as $G = \Theta_1(\Omega)$ and $G_0 = \Theta_1(\Omega_0)$, are strictly bounded simply connected domains in $\mathbb{C}$. The boundaries $\partial G$ and $\partial G_0$ are closed Jordan curves formed by the union of circular arcs on $\partial\mathbb{D}$ (the images of the free surfaces $\mathcal{E}$ and $\mathcal{E}_0$, including the point at infinity which maps to $1$) and interior arcs (the images of the wetted boundaries $\Gamma$ and $\Gamma_0$).

    By the Riemann mapping theorem, there exist conformal mappings between $G$ and $G_0$. Because $G$ and $G_0$ are bounded domains enclosed by Jordan curves, Carathéodory's extension theorem, as stated in Pommerenke \cite[Theorem 2.6]{pommerenke1992boundary}, guarantees that any such mapping extends to a global homeomorphism between their closures $\overline{G}$ and $\overline{G_0}$. To fix the map uniquely, choose Riemann maps from $G$ and $G_0$
to the disk, extended homeomorphically to the boundary. Disk automorphisms
act simply transitively on ordered triples of distinct boundary points
with the same cyclic order. The two endpoint--infinity triples here
have the same cyclic order. There is therefore exactly one conformal
map between $G$ and $G_0$ matching these three boundary points.
Let $x_{-}, x_{+} \in \partial\Omega$ denote the two intersection points between $\mathcal{E}$ and $\Gamma$, and note that $-1, 1 \in \partial\Omega_0$ are the corresponding intersection points between $\mathcal{E}_0$ and $\Gamma_0$. We define $\Psi : G \to G_0$ as the unique conformal mapping whose continuous extension satisfies the boundary correspondence $\Psi(\Theta_1(x_{-})) = \Theta_1(-1)$, $\Psi(\Theta_1(x_{+})) = \Theta_1(1)$, and $\Psi(\Theta_1(\infty)) = \Theta_1(\infty)$.

Furthermore, since the boundary arcs $\Gamma$ and $\mathcal{E}$ are smooth away from the contact points $x_\pm$, and their images under $\Theta_1$ are analytic arcs (portions of the unit circle and smooth interior curves), the conformal mapping $\Psi$, and consequently the global map $\Theta$, extend smoothly to the boundary away from the corners (see, e.g., \cite[Theorem 3.5]{pommerenke1992boundary}).

    Finally, we define the desired conformal mapping $\Theta:\Omega\to\Omega_0$ as the composition 
    \begin{equation}
        \Theta = \Theta_1^{-1} \circ \Psi \circ\Theta_1.
    \end{equation}
    This construction directly guarantees $\Theta(x_{-}) = -1$, $\Theta(x_{+}) = 1$, and $\Theta(\infty) = \infty$, thereby mapping the wetted surface $\Gamma$ uniquely onto $\Gamma_0$ and the free surface $\mathcal{E}$ onto $\mathcal{E}_0$.
\end{proof}

\begin{remark}
The normalization at infinity is essential for the uniqueness statement.
Preserving the two boundary pieces alone does not fix a conformal map.
Moreover, a homeomorphic boundary extension does not imply uniform
bounds on its derivative at corners. None of the domain identifications
or compactness arguments in the main text relies on such bounds.
For finite-energy functions, conformal invariance of the Dirichlet
integral follows by the chain rule and change of variables on compact
subdomains, followed by exhaustion.
\end{remark}
\section{The abstract wave equation and decay in the presence of a LAP}\label{appendix:lap}

In this appendix, we introduce this principle and the basic spectral notions required to formulate it in an abstract setting. We also prove an abstract local energy decay result, in which absolute continuity of the spectrum and compactness of localised observations are stated as separate hypotheses.

Throughout this section, $H$ denotes a Hilbert space. Let $A$ be a self-adjoint operator densely defined on $H$. By the spectral theorem (see Rudin \cite[Theorem~13.33]{rudin1991functional}), there exists a unique spectral family (or resolution of the identity) $\{E(\lambda)\}_{\lambda \in \mathbb{R}}$ of orthogonal projections on $H$ associated with $A$. For any $v, u \in H$, the map $\lambda \mapsto \langle E(\lambda)v, u \rangle_H$ is a function of bounded variation, thereby defining a complex Lebesgue--Stieltjes measure on $\mathbb{R}$. This measure rigorously defines the spectral representation of $A$ via the inner product:
\begin{equation}
    \langle Av, u \rangle_H = \int_{\mathbb{R}} \lambda \, {\rm d}\langle E(\lambda)v, u \rangle_H \qquad (v \in \mathcal{D}(A), u \in H).
\end{equation}
For convenience, we will also use the following vector and operator notations:
\begin{equation}
    Av = \int_{\mathbb{R}} \lambda \, {\rm d}E(\lambda)v \qquad (v \in \mathcal{D}(A)), \qquad A = \int_{\mathbb{R}} \lambda \, {\rm d}E(\lambda).
\end{equation}

This spectral representation extends naturally to Borel measurable functions $f: \mathbb{R} \to \mathbb{C}$ via the functional calculus. The operator $f(A)$ and its corresponding domain are defined as:
\begin{gather}
    f(A) = \int_{\mathbb{R}} f(\lambda) \, {\rm d}E(\lambda), \label{eq:functional_calculus} \\
    \mathcal{D}(f(A)) = \left\{ u \in H : \int_{\mathbb{R}} |f(\lambda)|^2 \, {\rm d}\|E(\lambda)u\|_H^2 < \infty \right\}. \label{def_domain_fA}
\end{gather}

A notable example is when $A$ is additionally positive (i.e., $\langle Au,u\rangle_H \geqslant 0$ for all $u\in \mathcal{D}(A)$). In this case, its spectral family is supported on $[0,\infty)$ and its square root is well defined. Evaluating $f(\lambda)=\lambda^{\frac12}$ yields
\begin{equation}
    A^{\frac12} = \int_{0}^\infty \lambda^{\frac12} \, {\rm d}E(\lambda)
\end{equation}
Another notable example is the resolvent $R_A(z)=(A-z \mathbb{I})^{-1}$ of $A$, defined for $z$ in its resolvent set $\rho(A)$. Taking $f(\lambda)=(\lambda-z)^{-1}$, we get
\begin{equation}
    R_A(z)=\int_{\mathbb{R}}\frac{1}{\lambda-z}{\rm d}E(\lambda).
\end{equation}

The following spectral cutoff estimates are standard consequences
of the Borel functional calculus; see \cite[Chapter~13]{rudin1991functional}.
We include their short proof to record the domain assertions and the
precise dependence on the cutoff thresholds.

\begin{lemma} \label{est-spectral-cutoff}
    Let $A$ be a nonnegative self-adjoint operator on a Hilbert space $H$.  Suppose that $\chi:\mathbb R\to\mathbb C$ is a bounded Borel measurable function. Then the following assertions hold.
    \begin{enumerate}
        \item For every $s\geqslant 0$ and $f\in \mathcal{D}(A^s),$ we have $\chi(A)f \in \mathcal{D}(A^s),$ and 
    \begin{equation}
        \norm{A^s\chi(A)f}_H \leqslant \left(\sup_{\lambda\in\mathbb R}|\chi(\lambda)|\right)
 \norm{A^sf}_H.
 \end{equation}
 \item If,  for some $a > 0$ 
 \begin{equation*}
   \chi(\lambda)=0,\qquad \lambda<a. 
\end{equation*}
Then for every $s \geqslant 0,$ and $f \in \mathcal{D}(A^{s+1/2}),$ we have 
\begin{equation}
    \norm{A^s\chi(A)f}_{H} \leqslant \frac{1}{\sqrt{a}} \left(\sup_{\lambda\in\mathbb R}|\chi(\lambda)|\right)  \norm{A^{s+1/2} f}_{H}
\end{equation}
 \item If, for some $b > 0,$ 
\begin{equation*}
   \chi(\lambda)=0,\qquad \lambda>b. 
\end{equation*} 
Then, for every $s\geqslant 0$ and $f\in H$, we have $\chi(A)f\in\mathcal{D}(A^s)$ and 
\begin{equation} 
  \norm{A^s\chi(A)f}_H
\leqslant  b^s\left(\sup_{\lambda\in\mathbb R}|\chi(\lambda)|\right)\norm{f}_H.  
\end{equation}
    \end{enumerate}

\end{lemma}

\begin{proof}
Put $M=\sup_{\lambda\in\mathbb R}|\chi(\lambda)|$.
We use the spectral family $E$ and the functional calculus introduced
in \eqref{eq:functional_calculus}--\eqref{def_domain_fA}.
The multiplication rule for spectral integrals implies that, for
$f\in H$ and $s\geqslant0$,
\[
 \chi(A)f\in\mathcal D(A^s)
 \quad\Longleftrightarrow\quad
 \int_{[0,\infty)}\lambda^{2s}|\chi(\lambda)|^2
                  \,{\rm d}\|E(\lambda)f\|_H^2<\infty.
\]
Whenever this integral is finite,
\begin{equation}\label{eq:spectral_multiplier_norm}
 \|A^s\chi(A)f\|_H^2
 =\int_{[0,\infty)}\lambda^{2s}|\chi(\lambda)|^2
                  \,{\rm d}\|E(\lambda)f\|_H^2.
\end{equation}
The interval includes zero; when $s=0$, the factor $\lambda^{2s}$
is understood to be identically one, consistently with $A^0=\mathbb I$.

For the first assertion, let $f\in\mathcal D(A^s)$. Then
\[
 \int_{[0,\infty)}\lambda^{2s}|\chi(\lambda)|^2
                  \,{\rm d}\|E(\lambda)f\|_H^2
 \leqslant M^2\int_{[0,\infty)}\lambda^{2s}
                  \,{\rm d}\|E(\lambda)f\|_H^2
 =M^2\|A^sf\|_H^2.
\]
This proves both preservation of $\mathcal D(A^s)$ and the stated
bound after taking square roots.

Under the hypothesis of the second assertion, the integrand in
\eqref{eq:spectral_multiplier_norm} vanishes on $[0,a)$.
For $\lambda\geqslant a$, we have
$\lambda^{2s}\leqslant a^{-1}\lambda^{2s+1}$. Thus, for
$f\in\mathcal D(A^{s+1/2})$,
\[
 \int_{[0,\infty)}\lambda^{2s}|\chi(\lambda)|^2
                  \,{\rm d}\|E(\lambda)f\|_H^2
 \leqslant\frac{M^2}{a}
       \int_{[0,\infty)}\lambda^{2s+1}
                  \,{\rm d}\|E(\lambda)f\|_H^2
 =\frac{M^2}{a}\|A^{s+1/2}f\|_H^2.
\]
Again the integral is finite, and taking square roots gives the
second estimate.

Finally, under the hypothesis of the third assertion, the integral
is supported in $[0,b]$, where $\lambda^{2s}\leqslant b^{2s}$.
Since the total mass of ${\rm d}\|E(\lambda)f\|_H^2$ is $\|f\|_H^2$,
\[
 \int_{[0,\infty)}\lambda^{2s}|\chi(\lambda)|^2
                  \,{\rm d}\|E(\lambda)f\|_H^2
 \leqslant b^{2s}M^2\int_{[0,b]}{\rm d}\|E(\lambda)f\|_H^2
 \leqslant b^{2s}M^2\|f\|_H^2.
\]
This proves that $\chi(A)f\in\mathcal D(A^s)$ for every $f\in H$,
as well as the third estimate.
\end{proof}

\subsection{The abstract wave equation and spectral representation}\label{sec:lap_background}
Assume now that $A$ is a self-adjoint positive operator. Let $w_0\in \mathcal{D}(A^{\frac12})$ and $w_1\in H$. We consider the initial value problem associated with the abstract wave equation
\begin{equation}\label{abstract_wave_pb}
    \left\{
    \begin{aligned}
        &\ddot w(t)+A w(t)=0\qquad(t\in\mathbb{R}),\\
        &w(0)=w_0,\quad \dot w(0)=w_1.
    \end{aligned}
    \right.
\end{equation}
The connection between the abstract wave equation and its spectral representation is a well-known consequence of the functional calculus. We state and prove the following proposition for the sake of completeness, with no claim of originality.

\begin{proposition}\label{prop_abstract_eq_spectral}
    The problem \eqref{abstract_wave_pb} admits a unique mild solution satisfying
    \begin{equation}
      w \in C(\mathbb{R},\mathcal{D}(A^{\frac12}))\cap C^1(\mathbb{R},H).  
    \end{equation}
    Furthermore, if $w_0\in \mathcal{D}(A)$ and $w_1\in \mathcal{D}(A^{\frac12})$, the solution is classical, satisfying
    \begin{equation}
     w \in C(\mathbb{R}, \mathcal{D}(A))\cap C^1(\mathbb{R},\mathcal{D}(A^{\frac12}))\cap C^2(\mathbb{R},H).     
    \end{equation}
    
    Moreover, the solution $w(t)$ admits the following spectral representation:
    \begin{equation}\label{spectral_rep_abstract_sol}
        w(t)=\int_0^{\infty} \cos(t\sqrt{\lambda} )\,{\rm d}E(\lambda)w_0+\int_0^{\infty}\frac{\sin(t\sqrt{\lambda})}{\sqrt{\lambda}}\,{\rm d}E(\lambda)w_1 \qquad(t\in\mathbb{R}).
    \end{equation}
    For either class of data, the energy is conserved:
    \begin{equation}\label{eq:abstract_energy_conservation}
        \|A^{\frac{1}{2}}w(t)\|_H^2+\|\dot w(t)\|_H^2 =\|A^{\frac{1}{2}}w_0\|_H^2+\|w_1\|_H^2.
    \end{equation}
\end{proposition}

\begin{proof}
    The well-posedness result follows by recasting \eqref{abstract_wave_pb} as a first-order system. We rewrite the equation as $\ddot{w}(t) + (\mathbb{I}+A) w(t) - w(t) = 0$, where $\mathbb{I}$ is the identity operator on $H$. Setting $z(t) = [w(t), \dot{w}(t)]^\top$, this is equivalent to the initial value problem $\dot{z}(t) = (\mathcal{A} + \mathcal{P})z(t)$ on the energy space $X = \mathcal{D}(A^{\frac12}) \times H$, where
    \begin{equation}
        \mathcal{A} = \begin{bmatrix} 0 & \mathbb{I} \\ -(\mathbb{I}+A) & 0 \end{bmatrix}, \qquad \mathcal{P} = \begin{bmatrix} 0 & 0 \\ \mathbb{I} & 0 \end{bmatrix}.
    \end{equation}
    We equip $X$ with the equivalent norm $\|[u, v]^\top\|_X^2 = \|(\mathbb{I}+A)^{\frac12} u\|_H^2 + \|v\|_H^2$. Since $\mathbb{I}+A$ is strictly positive, \cite[Proposition~3.7.6]{Obs_Book} ensures that the operator $\mathcal{A}$ with domain $\mathcal{D}(\mathcal{A}) = \mathcal{D}(A) \times \mathcal{D}(A^{\frac12})$ is skew-adjoint, and thus it generates a unitary group on $X$. The operator $\mathcal{P}$ is clearly a bounded linear perturbation on $X$. Thus, the sum $\mathcal{A} + \mathcal{P}$ generates a strongly continuous group on $X$, yielding the stated well-posedness.

    For the spectral representation, write $c_t(\lambda)=\cos(t\sqrt{\lambda})$ and let $s_t(\lambda)$ be $\sin(t\sqrt{\lambda})/\sqrt{\lambda}$ for $\lambda>0$, with $s_t(0)=t$. On any bounded time interval $|t|\leqslant T$, we have the uniform bounds
    \[
     |c_t(\lambda)|\leqslant 1,\qquad |s_t(\lambda)|\leqslant T,\qquad \sqrt{\lambda}|s_t(\lambda)|\leqslant 1.
    \]
    Consequently, $c_t(A)$ preserves $\mathcal{D}(A^{\frac12})$, while $s_t(A)$ maps $H$ into $\mathcal{D}(A^{\frac12})$. These estimates and the dominated convergence theorem applied to the spectral measure in \eqref{def_domain_fA} prove continuity in the graph norm of $W(t)=c_t(A)w_0+s_t(A)w_1$.
    
    Differentiation of the spectral integrals is justified by dominated convergence. The difference quotients of $c_t(\lambda)$ with respect to $t$ are bounded by $\sqrt{\lambda}$, which is square-integrable against the spectral measure of $w_0 \in \mathcal{D}(A^{\frac12})$; the difference quotients of $s_t(\lambda)$ are bounded by 1. This yields the strong derivative in $H$:
    \[
     \dot W(t)=-\sin(t\sqrt{A})A^{\frac12}w_0+\cos(t\sqrt{A})w_1.
    \]
    For the stronger initial data $w_0 \in \mathcal{D}(A)$ and $w_1 \in \mathcal{D}(A^{\frac12})$, the inequalities $\lambda|s_t(\lambda)|\leqslant \sqrt{\lambda}$ and $|\partial_t c_t(\lambda)|\leqslant \sqrt{\lambda}$ yield continuity of $W$ in $\mathcal{D}(A)$ and of $\dot W$ in $\mathcal{D}(A^{\frac12})$, respectively. A second differentiation then gives $\ddot W=-AW$ in $H$, with the prescribed initial data. Uniqueness identifies $W$ with the classical solution, and a density argument extends this to the mild solution.

    Finally, the energy variables satisfy
    \[
     \begin{bmatrix}A^{\frac12}w(t)\\ \dot w(t)\end{bmatrix}
     =\begin{bmatrix}\cos(t\sqrt{A})&\sin(t\sqrt{A})\\
     -\sin(t\sqrt{A})&\cos(t\sqrt{A})\end{bmatrix}
     \begin{bmatrix}A^{\frac12}w_0\\ w_1\end{bmatrix}.
    \]
    This matrix is unitary on $H\times H$ by the spectral calculus, which proves the energy conservation \eqref{eq:abstract_energy_conservation}.
\end{proof}

At the spectral point zero the solution contains
$E_A(\{0\})w_0+tE_A(\{0\})w_1$. Thus well-posedness does not
require $\ker A=\{0\}$. Conversely, energy conservation alone
need not give a uniform bound on $\|w(t)\|_H$ when zero belongs
to the spectrum. The inhomogeneous space used above avoids any
coercivity assumption on $A$.

\begin{remark}\label{rem-sol-rep-cutoff}
Let $A$ be a non-negative self-adjoint operator on $H$, let
$w_0\in\mathcal D(A^{\frac12})$ and $w_1\in H$, and let $w$ be the
solution of \eqref{abstract_wave_pb}. Suppose that
$\chi:\mathbb R\to\mathbb C$ is bounded and Borel measurable, and set
\[
 v(t)=\chi(A)w(t)\qquad\qquad(t\in\mathbb R).
\]
Then
\[
 v\in C(\mathbb R,\mathcal D(A^{\frac12}))\cap C^1(\mathbb R,H),
\]
and $v$ is the unique mild solution of
\[
 \left\{
 \begin{aligned}
  &\ddot v(t)+Av(t)=0\qquad(t\in\mathbb R),\\
  &v(0)=\chi(A)w_0,\qquad \dot v(0)=\chi(A)w_1.
 \end{aligned}
 \right.
\]
Indeed, \cref{est-spectral-cutoff} shows that $\chi(A)$ is bounded
both on $H$ and on $\mathcal D(A^{\frac12})$ equipped with its graph
norm. Thus $v$ has the asserted regularity and
$\dot v(t)=\chi(A)\dot w(t)$.
Moreover, the multiplication rule for the functional calculus,
applied to \eqref{spectral_rep_abstract_sol}, gives
\begin{equation}
 v(t) =\int_{[0,\infty)}\cos(t\sqrt\lambda)\chi(\lambda)
                         \,{\rm d}E(\lambda)w_0 + \int_{[0,\infty)}
       \frac{\sin(t\sqrt\lambda)}{\sqrt\lambda}\chi(\lambda)
                         \,{\rm d}E(\lambda)w_1
 \qquad\qquad(t\in\mathbb R),
\end{equation}
where $\sin(t\sqrt\lambda)/\sqrt\lambda$ is assigned the value $t$
at $\lambda=0$. This is precisely the solution formula in
\cref{prop_abstract_eq_spectral} for the initial data
$\chi(A)w_0$ and $\chi(A)w_1$, which proves the assertion.
If $\chi$ vanishes for all sufficiently large positive $\lambda$,
\cref{est-spectral-cutoff} also places these filtered data in
$\mathcal D(A^s)$ for every $s\geqslant0$; in particular, $v$ is then
a classical solution.
\end{remark}

\subsection{The Limiting Absorption Principle}\label{subsec:abstract_local_energy_decay}

In applications to wave propagation and scattering, understanding the limiting behaviour of the resolvent $R_A(z)=(A-z)^{-1}$ near the continuous spectrum is crucial. The Limiting Absorption Principle (LAP) provides a rigorous framework to make these limits well-defined in suitable weighted spaces; see Agmon \cite{agmon1975spectral} and \cite{HazardLoret2007}.

\begin{definition}\label{def_lap}
Let $A$ be a self-adjoint operator on a Hilbert space $H$. 
Let $V_{\downarrow}$ be a continuously and densely embedded subspace of $H$, and
let $V_{\uparrow}$ be the dual of $V_{\downarrow}$ with respect to the pivot space $H$.
We say that $A$ satisfies the \emph{Limiting Absorption Principle (LAP)} with respect to the scale of spaces 
\begin{equation}
V_{\downarrow} \hookrightarrow H \hookrightarrow V_{\uparrow},
\end{equation} 
if for each $\lambda$ in the continuous spectrum of $A$, the limits
    \begin{equation}
        R_A^\pm(\lambda) := \lim_{\varepsilon\to 0^+} R_A(\lambda \pm i \varepsilon)
    \end{equation}
    exist in $\mathcal{L}(V_{\downarrow},V_{\uparrow})$.
\end{definition}

On an open interval $I$ where the LAP limits exist locally uniformly in $\mathcal L(V_\downarrow,V_\uparrow)$, Stone's formula gives an absolutely continuous spectral measure. For the operator $A$, its density, understood through duality with $V_\downarrow$, is
\begin{equation}\label{stone_formula_density}
\frac{{\rm d}E(\lambda)}{{\rm d}\lambda} = \frac{1}{2\pi i} \left[ R_A^+(\lambda) - R_A^-(\lambda) \right]
\end{equation}
(see \cite[Proposition~7.1.3]{amrein1996c0}).

\begin{remark}
\label{rem:lap_to_abstract_decay}
For $f,g\in V_\downarrow$ and any Borel set $B$, the scalar spectral measure is given by $\mu_{f,g}(B)=\langle E_A(B)f,g\rangle_H$. On each compact interval in $(0,\infty)$, Stone's formula gives its density:
\[
    \frac{{\rm d}\mu_{f,g}}{{\rm d}\lambda}(\lambda)
       =\frac{1}{2\pi i}
          \langle[R_A^+(\lambda)-R_A^-(\lambda)]f,g\rangle_{V_{\uparrow},V_{\downarrow}}.
\]
The locally uniform resolvent convergence justifies passage to the limit under the integral against a compactly supported continuous function of $\lambda$; in particular this scalar measure is absolutely continuous on $(0,\infty)$. If $N\subset(0,\infty)$ is a Borel set of Lebesgue measure zero, then
\[
    \|E_A(N)f\|_H^2=\langle E_A(N)f,f\rangle_H=0
          \qquad(f\in V_\downarrow).
\]
Density of $V_\downarrow$ and boundedness of $E_A(N)$ then give $E_A(N)=0$ on $H$. This proves that the spectrum of $A$ in $(0,\infty)$ is purely absolutely continuous.
\end{remark}

Combining this spectral measure representation with the well-posedness result of Proposition~\ref{prop_abstract_eq_spectral} yields the following explicit formula for the wave evolution.

\begin{corollary}\label{corollary_formula_w_resolvents}
Let $A$ be a non-negative self-adjoint operator with purely absolutely continuous spectrum, satisfying the locally uniform LAP on $(0,\infty)$. For initial data $w_0 \in V_\downarrow \cap \mathcal D(A^{\frac12})$ and $w_1 \in V_\downarrow$, the solution $w(t)$ of the abstract wave problem \eqref{abstract_wave_pb} satisfies 
\begin{equation}\label{formula_w_resolvents}
\begin{aligned}
w(t) &= \frac{1}{2\pi i} \int_0^{\infty} \cos(t\sqrt{\lambda}) \, [R_A^+(\lambda) - R_A^-(\lambda)] w_0 \, {\rm d}\lambda \\
&\quad + \frac{1}{2\pi i} \int_0^{\infty} \frac{\sin(t\sqrt{\lambda})}{\sqrt{\lambda}} \, [R_A^+(\lambda) - R_A^-(\lambda)] w_1 \, {\rm d}\lambda,
\end{aligned}
\end{equation}
where the integrals are understood in the weak* sense in $V_\uparrow$.
\end{corollary}

\begin{remark}\label{rem:weak_stone_formulation}
Writing the representation formula of \cref{corollary_formula_w_resolvents} in duality form against a test function $\tilde w \in V_\downarrow$ yields a duality pairing formulation on both sides:
\begin{equation}\label{eq:weak_stone_formulation}
\begin{aligned}
\langle w(t), \tilde w \rangle_{V_\uparrow, V_\downarrow} &= \frac{1}{2\pi i} \int_0^\infty \cos(t\sqrt{\lambda}) \, \langle [R_A^+(\lambda) - R_A^-(\lambda)] w_0, \tilde w \rangle_{V_\uparrow, V_\downarrow} \, {\rm d}\lambda \\
&\quad+ \frac{1}{2\pi i} \int_0^{\infty} \frac{\sin(t\sqrt{\lambda})}{\sqrt{\lambda}} \, \langle[R_A^+(\lambda) - R_A^-(\lambda)] w_1, \tilde w\rangle_{V_{\uparrow},V_{\downarrow}} \, {\rm d}\lambda.
\end{aligned}
\end{equation}
By the well-posedness result of \cref{prop_abstract_eq_spectral}, the solution satisfies $w(t) \in H$, which allows the left-hand side to be alternatively written as the standard inner product in $H$.
\end{remark}

We next describe a perturbation approach, proposed in \cite{HazardLoret2007}, allowing the transfer of the LAP property from a reference self-adjoint operator $\widetilde A$ to a compact (in an appropriate sense) perturbation, such as our operator of interest $A$. To this end, an important assumption introduced in \cite{HazardLoret2007} is that $\widetilde A$ satisfies the \textit{stable} LAP condition, in the sense of the definition below. Physically, this property means that if a wave does not radiate energy, its total energy must be finite.

\begin{definition}\label{def:stable_lap}
    Let $\widetilde A$ be a self-adjoint operator on $\widetilde H$ satisfying the LAP with respect to the scale $\widetilde V_\downarrow \hookrightarrow \widetilde H \hookrightarrow \widetilde V_\uparrow$. We say that $\widetilde A$ satisfies the stable LAP if, for every $\lambda$ in the continuous spectrum of $\widetilde A$ and every $u \in \widetilde V_\downarrow$, the condition
    \begin{equation}
        \operatorname{Im} \langle \widetilde{\mathcal{R}}^\pm(\lambda)u, u \rangle_{\widetilde V_\uparrow, \widetilde V_\downarrow} = 0
    \end{equation}
    implies that $\widetilde{\mathcal{R}}^\pm(\lambda)u \in \widetilde H$.
\end{definition}

Let $A$ and $\widetilde A$ be non-negative self-adjoint operators on the Hilbert spaces
$H$ and $\widetilde H$, respectively. We assume that these Hilbert spaces are 
closed subspaces of another Hilbert space $\mathcal H$, so that
we have the orthogonal decompositions
\begin{equation}\label{eq:common_hilbert_space}
    \mathcal H=H\oplus H_{\mathrm{ext}}
              =\widetilde H\oplus\widetilde H_{\mathrm{ext}}.
\end{equation}
For each $z\in\mathbb C\setminus[0,\infty)$ we extend the resolvents of $A$ and $\widetilde A$ to operators
\[ \mathcal R(z),\ \widetilde{\mathcal R}(z)\in \mathcal{L}(\mathcal{H}), \]
by setting
\begin{equation}\label{eq:extended_original_resolvents}
\begin{aligned}
    \mathcal R(z)&=(A-z)^{-1}\oplus0,\\
    \widetilde{\mathcal R}(z)&=(\widetilde A-z)^{-1}\oplus0.
\end{aligned}
\end{equation}

We next state several assumptions which play an important role in the remaining part of this work.

(H1) $\widetilde A$ is a self-adjoint operator on $\widetilde H$ satisfying the stable LAP with respect to the scale $\widetilde V_\downarrow \hookrightarrow \widetilde H \hookrightarrow \widetilde V_\uparrow$.

(H2) There exist two other scales of Hilbert spaces with
\begin{equation}\label{eq:common_weighted_scale}
    \mathcal V_\downarrow\hookrightarrow\mathcal H
    \hookrightarrow\mathcal V_\uparrow,
    \qquad \mathcal V_\uparrow \subset (\mathcal V_\downarrow)',
\end{equation}
\begin{equation}\label{eq:common_weighted_scale_drept}
     V_\downarrow\hookrightarrow H
    \hookrightarrow V_\uparrow,
    \qquad V_\uparrow \subset (V_\downarrow)',
\end{equation}
with continuous dense embeddings and with 
\begin{equation}\label{eq:compatible_weighted_decompositions}
\begin{aligned}
    \mathcal V_\downarrow
    &=V_\downarrow\oplus H_{\mathrm{ext}}
      =\widetilde V_\downarrow\oplus\widetilde H_{\mathrm{ext}},\\
    \mathcal V_\uparrow
    &=V_\uparrow\oplus H_{\mathrm{ext}}
      =\widetilde V_\uparrow\oplus\widetilde H_{\mathrm{ext}},
\end{aligned}
\end{equation}
with equivalent norms. Thus the complementary Hilbert spaces are
included in both weighted spaces. This is the setting of
\cite[Section~2.5.3]{HazardLoret2007}.

(H3) Let $I\subset(0,\infty)$ be an open interval on which the $\widetilde A$
operator satisfies the LAP. We assume that the limits
\begin{equation}\label{eq:free_extended_boundary_values}
    \widetilde{\mathcal R}^\pm(\lambda)
    =\lim_{\varepsilon\to0^+}
      \widetilde{\mathcal R}(\lambda\pm i\varepsilon)
\end{equation}
exist in $\mathcal L(\mathcal V_\downarrow,\mathcal V_\uparrow)$,
locally uniformly for $\lambda\in I$.

For the following definition and transfer theorem, we set
\begin{equation}\label{eq:fixed_resolvent_difference}
    D_1=\mathcal R(-1)
                  -\widetilde{\mathcal R}(-1).
\end{equation}
Since $-1$ belongs to the resolvent sets of both non-negative
operators, $D_1$ is a bounded self-adjoint operator on
$\mathcal H$.

\begin{definition}\label{def:compact_perturbation}
    With the above notation, we say that $A$ is a compact perturbation of $\widetilde A$
    on $I$, if $\widetilde A$ and $A$ satisfy the assumptions (H1)-(H3) above and if the following two conditions hold:
    \begin{enumerate}
        \item[(i)] $D_1$ extends by density to a bounded
        operator from $\mathcal V_\uparrow$ to
        $\mathcal V_\downarrow$;
        \item[(ii)] for every $\lambda\in I$, both operators
        \begin{equation}\label{eq:compact_resolvent_product}
            K_1^\pm(\lambda)
            :=D_1\bigl[\iota+(\lambda+1)
                 \widetilde{\mathcal R}^\pm(\lambda)\bigr]
        \end{equation}
        are compact on $\mathcal V_\downarrow$.
    \end{enumerate}
    Inside the brackets in \eqref{eq:compact_resolvent_product},
    $\iota$ denotes the natural inclusion of
    $\mathcal V_\downarrow$ into $\mathcal V_\uparrow$.
\end{definition}

\begin{remark}\label{remark:sufficient_resolvent_compactness}
    A sufficient condition for Definition~\ref{def:compact_perturbation}
    is that $D_1$ extend to a compact operator from
    $\mathcal V_\uparrow$ to $\mathcal V_\downarrow$. Indeed, the
    operator in brackets in \eqref{eq:compact_resolvent_product}
    is bounded from $\mathcal V_\downarrow$ to
    $\mathcal V_\uparrow$, so its composition with $D_1$ is
    compact on $\mathcal V_\downarrow$. This sufficient condition
    is stronger than the definition. Changing the reference point would require additional weighted
    mapping properties of the resolvents; we use $-1$ throughout.
\end{remark}

\begin{remark}\label{remark:hazard_loret_equivalence}
    We explain precisely how this definition is related to
    \cite[Definition~2.9]{HazardLoret2007}. That definition concerns
    two bounded self-adjoint operators on the same Hilbert space. For
    unbounded operators on different spaces, the prescription in
    \cite[Section~2.5.3]{HazardLoret2007} is to apply it to the bounded
    operators $\mathcal R(-1)$ and
    $\widetilde{\mathcal R}(-1)$ on $\mathcal H$. Their difference
    is $D_1$, so the boundedness requirement in their
    definition is precisely condition~(i) above.

    To translate their compactness requirement, put
    $\mu=(\lambda+1)^{-1}$. The scalar identity
    \begin{equation}\label{eq:scalar_resolvent_change}
        \frac{1}{(\tau+1)^{-1}-(\lambda+1)^{-1}}
        =-(\lambda+1)
           \left(1+\frac{\lambda+1}{\tau-\lambda}\right)
    \end{equation}
    for $\tau\geqslant0$, $\tau\ne\lambda$, explains the change of
    spectral parameter in \cite[Section~2.5.3]{HazardLoret2007}. By the spectral
    calculus, the free boundary operator appearing in their
    definition, at the transformed parameter $\mu\pm i0$, is
    \begin{equation}\label{eq:translated_free_boundary_operator}
        -(\lambda+1)
          \bigl[\iota+(\lambda+1)
              \widetilde{\mathcal R}^\mp(\lambda)\bigr].
    \end{equation}
    The signs are interchanged because $z\mapsto(z+1)^{-1}$
    interchanges the upper and lower half-planes. On
    $\widetilde H_{\mathrm{ext}}$, the bounded comparison operator
    is zero and contributes $-\mu^{-1}\iota$ to
    \eqref{eq:translated_free_boundary_operator}. Thus the inclusion
    term acts on the whole common space; it cannot be omitted.
    Multiplication of \eqref{eq:translated_free_boundary_operator}
    by $D_1$ gives
    $-(\lambda+1)K_1^\mp(\lambda)$. Since the scalar factor
    is nonzero and both signs are required, their compactness
    condition is exactly condition~(ii). Hence
    Definition~\ref{def:compact_perturbation} is their definition
    expressed through the original resolvents, without introducing
    resolvents of bounded resolvent operators as new objects.
\end{remark}

We now give a self-contained version of the transfer argument of
\cite[Theorem~2.11]{HazardLoret2007}.
To make explicit the precise property of the free boundary values
used below, let $\widetilde P:\mathcal H\to\widetilde H$ be the
orthogonal projection associated with \eqref{eq:common_hilbert_space}.
By \eqref{eq:compatible_weighted_decompositions}, its restriction
maps $\mathcal V_\downarrow$ continuously into
$\widetilde V_\downarrow$. Besides the stable LAP, we require that
a free boundary value which belongs to the Hilbert space is an
ordinary solution of the corresponding equation. This compatibility
property is stated explicitly in the theorem; for the free operator $\Lambda_{\mathbb H}$ it is verified in
\cref{prop:stable_lap_free}. It replaces the appeal to a
generalized spectral representation in the formulation of the
transfer result.

\begin{theorem}
\label{thm:transfer_lap}
    Assume that $\widetilde A$ satisfies the stable LAP on $I$.
    Suppose also that, for every $\lambda\in I$, either sign, and
    $u\in\mathcal V_\downarrow$, the following implication holds:
    \begin{equation}\label{eq:free_boundary_consistency}
    \begin{split}
        v:=\widetilde{\mathcal R}^\pm(\lambda)u\in\mathcal H
        \quad\Longrightarrow\quad
        &v\in D(\widetilde A),\\
        &(\widetilde A-\lambda)v=\widetilde P u.
    \end{split}
    \end{equation}
    Here $D(\widetilde A)$ is identified with its copy in
    $\widetilde H\subset\mathcal H$. Suppose that $A$ is a compact
    perturbation of $\widetilde A$ in the sense of
    Definition~\ref{def:compact_perturbation} and has no eigenvalues
    in $I$. Then, for every $\lambda\in I$, the operators
    \begin{equation}\label{eq:fredholm_operator_original_resolvents}
        F_1^\pm(\lambda)
        :=\iota-(\lambda+1)K_1^\pm(\lambda)
    \end{equation}
    are invertible on $\mathcal V_\downarrow$. The perturbed
    extended resolvent has limits in
    $\mathcal L(\mathcal V_\downarrow,\mathcal V_\uparrow)$,
    locally uniformly for $\lambda\in I$, given by
    \begin{equation}\label{eq:transfer_original_resolvents}
        \mathcal R^\pm(\lambda)
        =\bigl[\widetilde{\mathcal R}^\pm(\lambda)
                    +K_1^\pm(\lambda)\bigr]
             \bigl[F_1^\pm(\lambda)\bigr]^{-1}.
    \end{equation}
    In this formula the range of $K_1^\pm(\lambda)$ is viewed
    in $\mathcal V_\uparrow$ through the natural embedding.
\end{theorem}

\begin{proof}
    All inner products and duality pairings are taken to be linear
    in the first argument. We take $\lambda>0$. The same argument applies to the other sign.

    First, we establish symmetry on the enlarged space. We denote by $D_1$ the difference of the resolvents in $-1$, including its bounded extension from
    $\mathcal V_\uparrow$ to $\mathcal V_\downarrow$.
    Self-adjointness of $D_1$ on $\mathcal H$ implies
    \begin{equation}\label{eq:extended_difference_real_pairing}
        \langle w,D_1 w\rangle_{\mathcal V_\uparrow,
                              \mathcal V_\downarrow}\in\mathbb R
        \qquad(w\in\mathcal V_\uparrow).
    \end{equation}
    Indeed, choosing $w_n\in\mathcal H$ such that
    $w_n\to w$ in $\mathcal V_\uparrow$, we have
    $D_1 w_n\to D_1 w$ in $\mathcal V_\downarrow$, so continuity of
    the duality pairing gives
    \[
        \langle w,D_1 w\rangle_{\mathcal V_\uparrow,
                              \mathcal V_\downarrow}
        =\lim_{n\to\infty}\langle w_n,D_1 w_n\rangle_{\mathcal H}.
    \]
    Each term on the right is real, proving
    \eqref{eq:extended_difference_real_pairing}.

    Next, we show the injectivity of the limiting Fredholm operator.
    Fix $\lambda\in I$, and abbreviate
    \[
        c=\lambda+1>0,\qquad
        S=\widetilde{\mathcal R}^\pm(\lambda),\qquad
        F=F_1^\pm(\lambda).
    \]
    Let $u\in\mathcal V_\downarrow$ satisfy $Fu=0$, and set
    \begin{equation}\label{eq:fredholm_kernel_auxiliary_vector}
        w=u+cSu\in\mathcal V_\uparrow.
    \end{equation}
    The definitions of $F$ and $K_1^\pm(\lambda)$ give
    \begin{equation}\label{eq:fredholm_kernel_relation}
        u=c D_1 w.
    \end{equation}
    Thus, by \eqref{eq:extended_difference_real_pairing},
    \begin{equation}\label{eq:kernel_zero_flux}
    \begin{split}
        0&=\operatorname{Im}\langle w,u\rangle_{\mathcal V_\uparrow,
                              \mathcal V_\downarrow}
          =\operatorname{Im}
            \bigl(\|u\|_{\mathcal H}^2+c\langle Su,u\rangle_{\mathcal V_\uparrow,
                              \mathcal V_\downarrow}\bigr)\\
         &=c\operatorname{Im}\langle Su,u\rangle_{\mathcal V_\uparrow,
                              \mathcal V_\downarrow}.
    \end{split}
    \end{equation}
    The extended resolvent $S$ is zero on
    $\widetilde H_{\mathrm{ext}}$. Hence its pairing with $u$
    is the native free-resolvent pairing with
    $\widetilde P u\in\widetilde V_\downarrow$.
    The stable LAP applied to \eqref{eq:kernel_zero_flux} gives
    \[
        v:=Su\in\widetilde H.
    \]
    By \eqref{eq:free_boundary_consistency}, we have
    \begin{equation}\label{eq:kernel_free_equation}
        v\in D(\widetilde A),\qquad
        (\widetilde A-\lambda)v=\widetilde P u.
    \end{equation}
    In particular, $w=u+cv$ belongs to $\mathcal H$. Moreover,
    \[
        \widetilde P w=\widetilde P u+cv
                     =(\widetilde A+1)v,
    \]
    so the ordinary free resolvent satisfies
    \begin{equation}\label{eq:kernel_free_reference_resolvent}
        \widetilde{\mathcal R}(-1)w=v.
    \end{equation}
    Combining \eqref{eq:fredholm_kernel_relation} with the definition of $D_1$ and
    \eqref{eq:kernel_free_reference_resolvent}, we obtain
    \begin{equation}\label{eq:kernel_perturbed_eigenvector}
        \mathcal R(-1)w
        =\widetilde{\mathcal R}(-1)w+D_1 w
        =v+c^{-1}u=c^{-1}w.
    \end{equation}
    Since the range of $\mathcal R(-1)$ is contained in
    $\mathcal D(A)\oplus\{0\}$, equation \eqref{eq:kernel_perturbed_eigenvector} implies that $w$
    belongs to this subspace and
    \[
        (A+1)w=cw,\qquad (A-\lambda)w=0.
    \]
    As $A$ has no eigenvalue at $\lambda$, it follows that
    $w=0$. Equation~\eqref{eq:fredholm_kernel_relation} then gives
    $u=0$, proving the injectivity of $F$.

    For completeness, we provide the compactness argument underlying
    the Fredholm alternative for the surjectivity and boundedness of the inverse. Setting
    $X=\mathcal V_\downarrow$ and
    $T=cK_1^\pm(\lambda)$, $T$ is compact on $X$ and
    $F=\iota-T$ is injective. There is a constant $m>0$ such
    that
    \begin{equation}\label{eq:fredholm_lower_bound}
        \|Fx\|_X\geqslant m\|x\|_X\qquad(x\in X).
    \end{equation}
    Otherwise, there would be unit vectors $x_n$ with $Fx_n\to0$.
    Compactness provides a subsequence for which $Tx_n$ converges.
    The identity $x_n=Fx_n+Tx_n$ then shows that this subsequence
    converges to a unit vector in $\ker F$, contradicting
    injectivity. Estimate~\eqref{eq:fredholm_lower_bound} also
    shows that $F(X)$, and more generally each $M_n=F^n(X)$,
    is closed.
    Suppose that $F(X)\ne X$. The closed subspaces $M_n$ would
    then satisfy $M_{n+1}\subsetneq M_n$ for every $n\geqslant0$:
    equality at any step, together with injectivity of $F^n$,
    would imply $F(X)=X$. Choosing a unit vector
    $x_n\in M_n\cap M_{n+1}^{\perp}$, where orthogonality is
    taken in the Hilbert space $X$, if $j>n$, then $Tx_j$ and
    $Fx_n$ belong to $M_{n+1}$. Since $Tx_n=x_n-Fx_n$, we have
    \[
        \|Tx_n-Tx_j\|_X
        =\|x_n-(Fx_n+Tx_j)\|_X\geqslant1.
    \]
    This contradicts the compactness of $T$. Thus $F$ is surjective,
    and \eqref{eq:fredholm_lower_bound} provides a bounded inverse.

    We now establish a resolvent identity at the fixed reference point. Let $z\in\mathbb C\setminus\mathbb R$ and put $t=z+1$.
    The ordinary resolvent identity remains valid after extension
    by zero. In particular,
    \begin{equation}\label{eq:reference_resolvent_factorisations}
    \begin{aligned}
        \bigl[\iota+t\mathcal R(z)\bigr]\mathcal R(-1)
          &=\mathcal R(z),\\
        \widetilde{\mathcal R}(-1)
          [\iota+t\widetilde{\mathcal R}(z)]
          &=\widetilde{\mathcal R}(z).
    \end{aligned}
    \end{equation}
    Substituting the difference defining $D_1$ and using
    \eqref{eq:reference_resolvent_factorisations} yields
    \begin{align*}
        &[\iota+t\mathcal R(z)]D_1
                     [\iota+t\widetilde{\mathcal R}(z)]\\
        &\quad=\mathcal R(z)[\iota+t\widetilde{\mathcal R}(z)]
              -[\iota+t\mathcal R(z)]
                   \widetilde{\mathcal R}(z)\\
        &\quad=\mathcal R(z)-\widetilde{\mathcal R}(z).
    \end{align*}
    Thus, we obtain the identity
    \begin{equation}\label{eq:resolvent_difference_fixed_point}
        \mathcal R(z)-\widetilde{\mathcal R}(z)
        ={}[\iota+(z+1)\mathcal R(z)]D_1
           [\iota+(z+1)\widetilde{\mathcal R}(z)].
    \end{equation}
    We can define bounded operators on $\mathcal V_\downarrow$ by
    \begin{equation}\label{eq:nonreal_fredholm_operators}
    \begin{aligned}
        K_1(z)
          &=D_1[\iota+(z+1)
                              \widetilde{\mathcal R}(z)],\\
        F_1(z)&=\iota-(z+1)K_1(z).
    \end{aligned}
    \end{equation}
    These mapping properties follow from condition~(i) of
    Definition~\ref{def:compact_perturbation} and the embeddings
    in \eqref{eq:common_weighted_scale}. Rearranging
    \eqref{eq:resolvent_difference_fixed_point} gives
    \begin{equation}\label{eq:resolvent_fredholm_identity}
        \mathcal R(z)F_1(z)
        =\widetilde{\mathcal R}(z)+K_1(z)
    \end{equation}
    as an identity from $\mathcal V_\downarrow$ to
    $\mathcal V_\uparrow$.

    Finally, we pass to the boundary in the operator norm.
    Let $J\subset I$ be a compact interval. The free LAP and
    boundedness of $D_1$ from $\mathcal V_\uparrow$ to
    $\mathcal V_\downarrow$ imply
    \begin{equation}\label{eq:fredholm_norm_convergence}
    \begin{split}
        \sup_{\lambda\in J}\bigl(&
          \|K_1(\lambda\pm i\varepsilon)
                       -K_1^\pm(\lambda)\|_{\mathcal L(X)}\\
          &+\|F_1(\lambda\pm i\varepsilon)
                       -F_1^\pm(\lambda)\|_{\mathcal L(X)}
                   \bigr)\longrightarrow0
        \quad(\varepsilon\to0^+).
    \end{split}
    \end{equation}
    Here $X=\mathcal V_\downarrow$ as above. For each
    $\varepsilon>0$, the free resolvent depends continuously on
    $\lambda$ in the weighted operator norm: this follows from
    its norm continuity on $\mathcal H$ and the continuous
    embeddings. Its locally uniform boundary limit is therefore
    continuous as well. Hence $K_1^\pm(\lambda)$ and
    $F_1^\pm(\lambda)$ are norm-continuous on $I$.

    As previously shown, each $F_1^\pm(\lambda)$ is
    invertible. The Neumann series gives continuity of their
    inverses: if $G$ is close enough to an invertible $F_0$, then
    \[
        G^{-1}
        =[\iota+F_0^{-1}(G-F_0)]^{-1}F_0^{-1},
    \]
    and the first inverse is the convergent geometric series
    whenever
    \[ \|F_0^{-1}(G-F_0)\|_{\mathcal{L}(X)}<1. \]
    Compactness of $J$ thus yields
    \[
        M_J:=\sup_{\lambda\in J}
             \|[F_1^\pm(\lambda)]^{-1}\|_{\mathcal L(X)}
             <\infty.
    \]
    By \eqref{eq:fredholm_norm_convergence}, for all sufficiently
    small $\varepsilon>0$ the operators
    $F_1(\lambda\pm i\varepsilon)$ are invertible,
    uniformly for $\lambda\in J$, and
    \begin{equation}\label{eq:fredholm_inverse_estimates}
    \begin{aligned}
        \|F_1(\lambda\pm i\varepsilon)^{-1}\|_{\mathcal L(X)}
             &\leqslant2M_J,\\
        \bigl\|F_1(\lambda\pm i\varepsilon)^{-1}
             -[F_1^\pm(\lambda)]^{-1}\bigr\|_{\mathcal L(X)}
             &\leqslant2M_J^2
               \bigl\|F_1(\lambda\pm i\varepsilon)
                       -F_1^\pm(\lambda)\bigr\|_{\mathcal L(X)}.
    \end{aligned}
    \end{equation}
    The second inequality also follows directly by subtracting
    the two inverses. We can therefore solve
    \eqref{eq:resolvent_fredholm_identity} to obtain
    \[
        \mathcal R(\lambda\pm i\varepsilon)
        =[\widetilde{\mathcal R}(\lambda\pm i\varepsilon)
                         +K_1(\lambda\pm i\varepsilon)]
                     F_1(\lambda\pm i\varepsilon)^{-1}.
    \]
    The first factor converges in
    $\mathcal L(\mathcal V_\downarrow,\mathcal V_\uparrow)$,
    uniformly for $\lambda\in J$, and the second converges in
    $\mathcal L(\mathcal V_\downarrow)$ by
    \eqref{eq:fredholm_inverse_estimates}. Both are uniformly
    bounded there. Continuity of the multiplication of bounded
    operators now gives \eqref{eq:transfer_original_resolvents}
    in the asserted operator norm, locally uniformly on $I$.
\end{proof}

\subsection{Local energy decay from absolute continuity and compactness}

The role of the LAP in our decay argument is to establish the absolute continuity of the spectral measures. Once this spectral property is known, the decay follows from an abstract localisation argument. We provide its proof in detail, keeping the spectral assumption separate from the compactness property of the observations. Related spectral decay arguments are discussed in Teschl \cite{teschl2009mathematical}. Here, compactness is required only after restriction to a bounded spectral interval separated from zero.

\begin{proposition}
\label{prop:abstract_local_energy_decay}
Let $A$ be a non-negative self-adjoint operator on $H$ with purely absolutely continuous spectrum, and let $E_A$ denote its projection-valued spectral measure. Let $Y_0,Y_1$ be Hilbert spaces and let $C_j\in\mathcal L(H,Y_j)$, $j=0,1$. Suppose that, for every $0<a<b<\infty$ and $j=0,1$,
\begin{equation}\label{eq:compact_spectral_observations}
    C_j E_A([a,b]):H\longrightarrow Y_j
    \quad\hbox{is compact}.
\end{equation}
For $w_0\in \mathcal D(A^{\frac12})$ and $w_1\in H$, let $w$ be the solution of $\ddot w+Aw=0$, $w(0)=w_0$, $\dot w(0)=w_1$. Then, as $|t|\to\infty$,
\begin{equation}\label{eq:abstract_weak_energy_decay}
    A^{\frac12}w(t) \to 0 \quad \text{and} \quad \dot w(t) \to 0 \quad \text{weakly in } H,
\end{equation}
and
\begin{equation}\label{eq:abstract_local_energy_decay}
    \|C_0 A^{\frac12}w(t)\|_{Y_0}
        +\|C_1\dot w(t)\|_{Y_1}\longrightarrow0.
\end{equation}
\end{proposition}

\begin{proof}
    First, we prove the weak decay of the half-wave groups. For $f,h\in H$, define the complex measure
    \[
        \mu_{f,h}(B)=\langle E_A(B)f,h\rangle_H
           \qquad(B\subset\mathbb R\hbox{ Borel}).
    \]
    This measure has finite total variation. Indeed, for any finite family of disjoint Borel sets $B_k$, orthogonality of the spectral projections and the Cauchy--Schwarz inequality give
    \[
    \begin{split}
        \sum_k|\mu_{f,h}(B_k)|
           &\leqslant\sum_k
                \|E_A(B_k)f\|_H\,\|E_A(B_k)h\|_H\\
           &\leqslant
              \left(\sum_k\|E_A(B_k)f\|_H^2\right)^{\frac12}
              \left(\sum_k\|E_A(B_k)h\|_H^2\right)^{\frac12}
            \leqslant\|f\|_H\|h\|_H.
    \end{split}
    \]
    Pure absolute continuity and the Radon--Nikodym theorem therefore provide $\rho_{f,h}\in L^1(0,\infty)$ such that ${\rm d}\mu_{f,h}(\lambda)=\rho_{f,h}(\lambda)\,{\rm d}\lambda$. Since there is no mass at zero, the spectral theorem, followed by the change of variable $\lambda=r^2$, gives
    \[
    \begin{split}
        \langle e^{it\sqrt A}f,h\rangle_H
           &=\int_0^\infty e^{it\sqrt\lambda}
                        \rho_{f,h}(\lambda)\,{\rm d}\lambda\\
           &=\int_0^\infty e^{itr}
                        2r\rho_{f,h}(r^2)\,{\rm d}r.
    \end{split}
    \]
    The amplitude of this final integral belongs to $L^1(0,\infty)$, because
    \[
        \int_0^\infty 2r|\rho_{f,h}(r^2)|\,{\rm d}r
           =\int_0^\infty|\rho_{f,h}(\lambda)|\,{\rm d}\lambda
           <\infty.
    \]
    The scalar Riemann--Lebesgue lemma now guarantees that
    \begin{equation}\label{eq:weak_half_wave_decay}
        e^{it\sqrt A}f \to 0
           \quad \text{weakly in } H \text{ as }|t|\to\infty,
           \qquad f\in H.
    \end{equation}

    Next, we lift this to strong decay after spectral localisation. Let $C$ be either $C_0$ or $C_1$, with target space $Y$, and set $P_{a,b}=E_A([a,b])$. For fixed $0<a<b<\infty$, the compactness of $CP_{a,b}$ and the weak convergence \eqref{eq:weak_half_wave_decay} imply
    \[
        \|CP_{a,b}e^{it\sqrt A}f\|_Y\longrightarrow0 \quad \text{as }|t|\to\infty.
    \]
    To justify this implication explicitly, note that any sequence of times $|t_n|\to\infty$ yields a bounded, weakly null sequence $e^{it_n\sqrt A}f$. Its image under a compact operator has norm-convergent subsequences, and every such limit must be zero by weak convergence. Thus the whole image sequence tends to zero.
    Since the spectral projection commutes with the unitary group, we can write
    \[
        \|Ce^{it\sqrt A}f\|_Y
           \leqslant\|CP_{a,b}e^{it\sqrt A}f\|_Y
                 +\|C\|\,\|(I-P_{a,b})f\|_H.
    \]
    It follows that
    \[
        \limsup_{|t|\to\infty}\|Ce^{it\sqrt A}f\|_Y
           \leqslant\|C\|\,\|(I-P_{a,b})f\|_H.
    \]
    Since $A\geqslant0$ and $E_A(\{0\})=0$, we have $E_A([1/n,n])f\to f$ in $H$. Taking $a=1/n$ and $b=n$, and letting $n\to\infty$, yields
    \begin{equation}\label{eq:localized_half_wave_decay}
        \|Ce^{it\sqrt A}f\|_Y\longrightarrow0
              \text{ as }|t|\to\infty,\qquad f\in H.
    \end{equation}
    Notice that the approximation error is uniform in time; no estimate on the resolvent at zero or at infinity is required.

    Finally, we apply this to the wave equation. Set $f_0=A^{\frac12}w_0$ and $f_1=w_1$. The spectral solution formula gives, in $H$,
    \begin{equation}\label{eq:wave_energy_half_wave_formula}
    \begin{aligned}
        A^{\frac12}w(t)
           &=\cos(t\sqrt A)f_0+\sin(t\sqrt A)f_1,\\
        \dot w(t)
           &=-\sin(t\sqrt A)f_0+\cos(t\sqrt A)f_1.
    \end{aligned}
    \end{equation}
    These identities hold for the stated energy data by the functional calculus, since all four multipliers on the right are bounded. Writing the sine and cosine as linear combinations of $e^{it\sqrt A}$ and $e^{-it\sqrt A}$, the first step proves \eqref{eq:abstract_weak_energy_decay}. The second step, applied with $C_0$ to the first identity and with $C_1$ to the second, proves the strong local energy decay \eqref{eq:abstract_local_energy_decay}.
\end{proof}

\begin{remark}
\label{rem:weak_energy_zero_threshold}
The weak convergence in \eqref{eq:abstract_weak_energy_decay} concerns the energy variables $(A^{\frac12}w,\dot w)$. It does not in general imply weak convergence of $(w,\dot w)$ in $\mathcal D(A^{\frac12})\times H$ with its inhomogeneous graph norm: conservation of $\|A^{\frac12}w(t)\|_H^2+\|\dot w(t)\|_H^2$ need not bound $\|w(t)\|_H$ when the spectrum reaches zero.

For example, consider the operator $A=\Lambda_{\mathbb{H}}$ on $L^2(\mathbb R)$ (see \cref{sec:free_waves}), with initial data $w_0=0$, and $w_1\in\mathscr S(\mathbb R)$ satisfying $\widehat w_1(0)\ne0$. By Plancherel's theorem, we have
\[
    \|w(t)\|_{L^2(\mathbb{R})}^2
       =\int_{\mathbb R}
          \frac{\sin^2(t\sqrt{|\xi|})}{|\xi|}
                    |\widehat w_1(\xi)|^2\,{\rm d}\xi.
\]
The continuity of $\widehat w_1$ at zero provides $c,\delta>0$ such that, for $t>0$,
\[
    \|w(t)\|_{L^2(\mathbb{R})}^2
       \geqslant c\int_0^\delta
                  \frac{\sin^2(t\sqrt\xi)}{\xi}\,{\rm d}\xi
       =2c\int_0^{t\sqrt\delta}\frac{\sin^2 r}{r}\,{\rm d}r.
\]
This last expression grows logarithmically as $t\to\infty$. This can be seen by using $\sin^2r=(1-\cos(2r))/2$ on $[1,t\sqrt\delta]$ and integrating by parts to bound the cosine integral. Thus, $w(t)$ is not even bounded in $L^2(\mathbb{R})$, although the initial data are smooth and rapidly decreasing. A sufficient additional condition for weak decay in the inhomogeneous space is $w_1\in \mathcal D(A^{-\frac12})$. Indeed, in this case,
\[
    w(t)=\cos(t\sqrt A)w_0+
                   \sin(t\sqrt A)A^{-\frac12}w_1
\]
is weakly null in $H$ by \eqref{eq:weak_half_wave_decay}. Together with \eqref{eq:abstract_weak_energy_decay}, this yields weak convergence in $\mathcal D(A^{\frac12})\times H$.
\end{remark}

\clearpage
\printbibliography

@book {MR1925354,
    AUTHOR = {Kuznetsov, N. and Maz'ya, V. and Vainberg, B.},
     TITLE = {Linear water waves},
      NOTE = {A mathematical approach},
 PUBLISHER = {Cambridge University Press, Cambridge},
      YEAR = {2002},
     PAGES = {xviii+513},
      ISBN = {0-521-80853-7},
   MRCLASS = {76B15 (35Q35 76B03)},
  MRNUMBER = {1925354},
MRREVIEWER = {Mark\ D.\ Groves},
       DOI = {10.1017/CBO9780511546778},
       URL = {https://doi.org/10.1017/CBO9780511546778},
}

@article{KM13,
  title={Estimate for a solution to the water wave problem in the presence of a submerged body},
  author={Kamotski, Ilia and Maz'ya, V},
  journal={Russian Journal of Mathematical Physics},
  volume={20},
  number={4},
  pages={453--467},
  year={2013},
  publisher={Springer}
}

@book{agmon1975spectral,
  author    = {Agmon, S.},
  title     = {Spectral Properties of {S}chr{\"o}dinger Operators and Scattering Theory},
  publisher = {Scuola Normale Superiore},
  year      = {1975},
  lccn      = {77378038}
}

@book{ahlfors1979complex,
  title     = {Complex Analysis},
  author    = {Ahlfors, Lars V},
  edition   = {3rd},
  year      = {1979},
  publisher = {McGraw-Hill}
}

@book {Hardy,
    AUTHOR = {Hardy, G. H. and Littlewood, J. E. and P\'{o}lya, G.},
     TITLE = {Inequalities},
      NOTE = {2d ed},
 PUBLISHER = {Cambridge, at the University Press},
      YEAR = {1952},
     PAGES = {xii+324},
   MRCLASS = {27.0X},
  MRNUMBER = {46395},
}

@book {TW05,
    AUTHOR = {Toselli, Andrea and Widlund, Olof},
     TITLE = {Domain decomposition methods---algorithms and theory},
    SERIES = {Springer Series in Computational Mathematics},
    VOLUME = {34},
 PUBLISHER = {Springer-Verlag, Berlin},
      YEAR = {2005},
     PAGES = {xvi+450},
      ISBN = {3-540-20696-5},
   MRCLASS = {65-02 (65N55 74S05 76M10)},
  MRNUMBER = {2104179},
MRREVIEWER = {R\'{e}mi\ Vaillancourt},
       DOI = {10.1007/b137868},
       URL = {https://doi.org/10.1007/b137868},
}

@book {dlv2,
    AUTHOR = {Dautray, Robert and Lions, Jacques-Louis},
     TITLE = {Mathematical analysis and numerical methods for science and
              technology. {V}ol. 2},
      NOTE = {Functional and variational methods,
              With the collaboration of Michel Artola, Marc Authier,
              Philippe B\'{e}nilan, Michel Cessenat, Jean Michel Combes,
              H\'{e}l\`ene Lanchon, Bertrand Mercier, Claude Wild and Claude
              Zuily,
              Translated from the French by Ian N. Sneddon},
 PUBLISHER = {Springer-Verlag, Berlin},
      YEAR = {1988},
     PAGES = {xvi+561},
      ISBN = {3-540-19045-7; 3-540-66098-4},
   MRCLASS = {00A05},
  MRNUMBER = {969367},
       DOI = {10.1007/978-3-642-61566-5},
       URL = {https://doi.org/10.1007/978-3-642-61566-5},
}

@book {galdinsbook,
    AUTHOR = {Galdi, G. P.},
     TITLE = {An introduction to the mathematical theory of the
              {N}avier-{S}tokes equations},
    SERIES = {Springer Monographs in Mathematics},
   EDITION = {Second},
      NOTE = {Steady-state problems},
 PUBLISHER = {Springer, New York},
      YEAR = {2011},
     PAGES = {xiv+1018},
      ISBN = {978-0-387-09619-3},
   MRCLASS = {35Q30 (35-02 76D03 76D05 76D07)},
  MRNUMBER = {2808162},
       DOI = {10.1007/978-0-387-09620-9},
       URL = {https://doi.org/10.1007/978-0-387-09620-9},
}

@book{amrein1996c0,
  title     = {$C_0$-Groups, Commutator Methods and Spectral Theory of $N$-Body Hamiltonians},
  author    = {Amrein, Werner O. and Boutet de Monvel, Anne and Georgescu, Vladimir},
  series    = {Progress in Mathematics},
  volume    = {135},
  year      = {1996},
  publisher = {Birkh{\"a}user Verlag},
  address   = {Basel},
  doi       = {10.1007/978-3-0348-0733-3}
}

@book{axler2013harmonic,
  author    = {Axler, S. and Bourdon, P. and Ramey, W.},
  title     = {Harmonic Function Theory},
  series    = {Graduate Texts in Mathematics},
  publisher = {Springer New York},
  year      = {2013},
  isbn      = {9781475781373},
  doi       = {10.1007/978-1-4757-8137-3}
}

@article{Brown1994,
  author    = {Russell Brown},
  title     = {The mixed problem for {L}aplace's equation in a class of {L}ipschitz domains},
  journal   = {Communications in Partial Differential Equations},
  volume    = {19},
  number    = {7-8},
  pages     = {1217--1233},
  year      = {1994},
  publisher = {Taylor \& Francis},
  doi       = {10.1080/03605309408821052}
}

@book{Dauge1988,
  author    = {Dauge, Monique},
  title     = {Elliptic Boundary Value Problems on Corner Domains},
  series    = {Lecture Notes in Mathematics},
  volume    = {1341},
  publisher = {Springer},
  address   = {Berlin, Heidelberg},
  year      = {1988},
  doi       = {10.1007/BFb0086682}
}

@article{EichmairMetzger2012,
  author  = {Eichmair, Michael and Metzger, Jan},
  title   = {On large volume preserving stable {CMC} surfaces
             in initial data sets},
  journal = {Journal of Differential Geometry},
  volume  = {91},
  number  = {1},
  pages   = {81--102},
  year    = {2012}
}

@article{HazardLoret2007,
  author    = {Hazard, Christophe and Loret, Fran{\c{c}}ois},
  title     = {Generalized eigenfunction expansions for conservative scattering problems with an application to water waves},
  journal   = {Proceedings of the Royal Society of Edinburgh Section A: Mathematics},
  volume    = {137},
  number    = {5},
  pages     = {995--1035},
  year      = {2007},
  publisher = {Royal Society of Edinburgh},
  doi       = {10.1017/S0308210506000138}
}

@article{John1949OnTM,
  author    = {John, Fritz},
  title     = {On the motion of floating bodies. {I}},
  journal   = {Communications on Pure and Applied Mathematics},
  volume    = {2},
  pages     = {13--57},
  year      = {1949},
  doi       = {10.1002/cpa.3160020102}
}

@article{John1950OnTM,
  author    = {John, Fritz},
  title     = {On the motion of floating bodies {II}. Simple harmonic motions},
  journal   = {Communications on Pure and Applied Mathematics},
  volume    = {3},
  pages     = {45--101},
  year      = {1950},
  doi       = {10.1002/cpa.3160030106}
}

@book{lannes2013water,
  author    = {Lannes, David},
  title     = {The Water Waves Problem: Mathematical Analysis and Asymptotics},
  volume    = {188},
  publisher = {American Mathematical Society},
  year      = {2013},
  isbn      = {9780821894705},
  doi       = {10.1090/surv/188}
}

@misc{lannes2025fjohnmodelcummins,
  author        = {Lannes, David and Paulsen, Martin Oen},
  title         = {The {F}. {J}ohn model and {C}ummins' equations for freely floating objects},
  year          = {2025},
  eprint        = {2511.16826},
  howpublished  = {arXiv:2511.16826},
  archivePrefix = {arXiv},
  primaryClass  = {math.AP}
}

@book{lunardi2018interpolation,
  title={Interpolation Theory},
  author={Lunardi, A.},
  isbn={9788876426391},
  series={Publications of the Scuola Normale Superiore},
  url={https://books.google.fr/books?id=a7UStQEACAAJ},
  year={2018},
  publisher={Scuola Normale Superiore}
}

@book{mclean2000strongly,
  title={Strongly Elliptic Systems and Boundary Integral Equations},
  author={McLean, W.C.H.},
  isbn={9780521663755},
  lccn={99030938},
  url={https://books.google.fr/books?id=RILqjEeMfK0C},
  year={2000},
  publisher={Cambridge University Press}
}

@article{lannes_ming_2026well,
  title={Well-posedness of F. John's floating body problem for a fixed object},
  author={Lannes, David and Ming, Mei},
  journal={Revista Matem{\'a}tica Iberoamericana},
  volume={42},
  number={2},
  pages={665--746},
  year={2026},
  doi={10.4171/RMI/1592}
}

@book{Obs_Book,
  author    = {Tucsnak, Marius and Weiss, George},
  title     = {Observation and Control for Operator Semigroups},
  publisher = {Birkh\"auser Verlag},
  address   = {Basel},
  year      = {2009},
  doi       = {10.1007/978-3-7643-8994-9}
}

@article{Ocqueteau2026,
  author  = {Ocqueteau, Vicente and Tucsnak, Marius},
  title   = {On an Initial Value Problem Describing the Small Oscillations of a Floating Cylinder},
  journal = {Journal of Mathematical Fluid Mechanics},
  year    = {2026},
  volume  = {28},
  eid     = {66},
  doi     = {10.1007/s00021-026-01047-0},
  issn    = {1422-6952}
}

@book{pommerenke1992boundary,
  title     = {Boundary Behaviour of Conformal Maps},
  author    = {Pommerenke, Christian},
  volume    = {299},
  year      = {1992},
  publisher = {Springer Berlin, Heidelberg}
}

@book{rudin1991functional,
  author    = {Rudin, Walter},
  title     = {Functional Analysis},
  edition   = {2nd},
  publisher = {McGraw-Hill},
  year      = {1991},
  isbn      = {978-0-07-054236-5}
}

@book{stoker1957water,
  author    = {Stoker, J. J.},
  title     = {Water Waves: The Mathematical Theory with Applications},
  series    = {Wiley Classics Library},
  publisher = {Interscience Publishers},
  year      = {1957},
  isbn      = {978-0-471-57034-9},
  note      = {Reprint of the 1957 original}
}

@book{teschl2009mathematical,
  author    = {Teschl, Gerald},
  title     = {Mathematical Methods in Quantum Mechanics: With Applications to {S}chr{\"o}dinger Operators},
  series    = {Graduate Studies in Mathematics},
  volume    = {99},
  publisher = {American Mathematical Society},
  year      = {2009},
  doi       = {10.1090/gsm/099}
}

@article{Ursell1949ONTH,
  author    = {Ursell, Fritz},
  title     = {On the heaving motion of a circular cylinder on the surface of a fluid},
  journal   = {Quarterly Journal of Mechanics and Applied Mathematics},
  volume    = {2},
  pages     = {218--231},
  year      = {1949},
  doi       = {10.1093/qjmam/2.2.218}
}

@phdthesis{Weber_1976,
  author      = {Weber, Charles F.},
  title       = {The existence and decay of water waves in the presence of a fixed obstacle},
  school      = {New York University, Courant Institute of Mathematical Sciences},
  year        = {1976},
  address     = {New York, NY},
  note        = {Available from University Microfilms International}
}
\end{document}